\documentclass[12pt,reqno]{amsart}
\usepackage{amssymb}
\usepackage{amsfonts}
\usepackage{indentfirst}
\usepackage{amsopn}
\usepackage{amsmath}
\usepackage{amsthm}
\usepackage{amscd}
\usepackage{graphicx}
\usepackage{epstopdf}

\usepackage[english]{babel}
\usepackage{amsmath}
\usepackage{subfigure}
\usepackage{hyperref}
\usepackage{enumerate}

\newtheorem{lemma}{Lemma}[section]

\newtheorem{proposition}[lemma]{Proposition}

\newtheorem{definition}[lemma]{Definition}
\newtheorem{remark}[lemma]{Remark}

\numberwithin{equation}{section}

\newcommand{\pd}[2]{\frac{\partial {#1}}{\partial {#2}}}
\newcommand{\beq}{\begin{equation}}
\newcommand{\eeq}{\end{equation}}
\newcommand{\be}{\begin{equation*}}
\newcommand{\ee}{\end{equation*}}

\newcommand{\RE}{\mathbb R}
\newcommand{\CO}{\mathbb C}

\newcommand{\DD}{\mathsf D}

\newcommand{\GG}{\mathcal{G}}

\newcommand{\Ker}{\operatorname{Ker}\,}

\newcommand{\al}{\alpha}

\newcommand{\sech}{\operatorname{sech}\,}
\newcommand{\f}{\frac}

\newcommand{\K}{\mathcal K}
\newcommand{\x}{\underline{x}}
\newcommand{\y}{\underline{y}}

\newcommand{\qmax}{q_{\rm max}}

\newcommand{\der}[2]{\frac{\mathrm{d}^{#1}}{\mathrm{d}#2^{#1}}}
\usepackage{latexsym,epsfig,bm,xcolor}

\DeclareMathOperator{\rank}{rank}

\newcommand{\blue}[1]{{\color{blue}#1}}

\newcommand{\black}[1]{{\color{black}#1}}

\usepackage{tikz}
\tikzstyle{nodino}=[circle,draw,fill,inner sep=0pt]
\tikzstyle{infinito}=[circle,inner sep=0pt,minimum size=0mm]
\tikzstyle{nodo}=[circle,draw,fill,inner sep=0pt,minimum size=0.5* width("k")]

\begin{document}

\title[Standing waves on quantum graphs]{\bf Standing waves on quantum graphs}

\author[A.Kairzhan]{Adilbek Kairzhan}
\address[A.Kairzhan]{Department of Mathematics, University of Toronto, 40 St. George St., Room 6290
Toronto, Ontario, Canada}
\email{kairzhan@math.toronto.edu}

\author[D. Noja]{Diego Noja}
\address[D. Noja]{Dipartimento di Matematica e Applicazioni, Universit\`a di Milano Bicocca, via R. Cozzi 55, 20125 Milano, Italy}
\email{diego.noja@unimib.it}

\author[D.E. Pelinovsky]{Dmitry E. Pelinovsky}
\address[D.E. Pelinovsky]{Department of Mathematics and Statistics, McMaster University,
Hamilton, Ontario, Canada, L8S 4K1}
\email{dmpeli@math.mcmaster.ca}

\date{\today}
\maketitle

\begin{abstract}
We review evolutionary models on quantum graphs expressed by linear and nonlinear partial differential equations. Existence and stability of the standing waves trapped on quantum graphs are studied by using methods of the variational theory, dynamical systems on a phase plane, and the Dirichlet--to--Neumann mappings.
\end{abstract}

\tableofcontents

\begin{footnotesize}
	
\end{footnotesize}

\section{Introduction}

The goal of this review is to introduce readers to the subject of standing waves on quantum graphs described by nonlinear partial differential equations. The linear evolutionary models and the spectral theory for quantum graphs are well covered by several monographs \cite{alimehmeti,BK13,Exner} and many recent publications. The first review of the nonlinear 
evolutionary models was published some time ago \cite{N14} and was complemented by the recent reviews \cite{Adami,ASTrev,Dovetta} which explained various technical aspects of mathematical analysis of the ground state on the quantum graphs. Compared to these publications, we would like to present a general overview of how the standing waves arise in the nonlinear models and how their existence and stability can be analyzed  with different analytical methods.

\subsection{Main definitions} 

\begin{figure}[htb]
\begin{center}
\begin{tikzpicture}[xscale= .60,yscale=.4]
\node at (-11.5,0) [infinito] (GS) {\scriptsize $\infty$};
\node at (-9,0) [infinito] (GS) {};
\node at (-1,0) [nodo] (G1) {};
\node at (1,1) [nodo] (G2) {};
\node at (2,-1) [nodo](G3) {};
\node at (4,1) [nodo](G4) {};
\node at (5,-1)[nodo] (G5) {};
\node at (15.5,1) [infinito] (GD2) {\scriptsize $\infty$};
\node at (16.5,-1)  [infinito] (GD1) {\scriptsize $\infty$};

\node at (14,-1)  [infinito] (GD1) {};
\node at (13,1)  [infinito] (GD2) {};

\node at (2.7,-2)  [infinito] (G) {$ \GG$};

\draw[-] (GS)--(G1);
\draw[-] (G1)--(G2);
\draw[-] (G1)--(G3);
\draw[-] (G2)--(G3);
\draw[-] (G3)--(G5);
\draw[-] (G2)--(G4);
\draw[-] (G4)--(G5);
\draw[-] (G4)--(GD2);
\draw[-] (G5)--(GD1);
\draw[dashed] (GD1)--(16,-1);
\draw[dashed] (GD2)--(15,1);
\draw[dashed] (GS)--(-11,0);

\end{tikzpicture}
\end{center}
\caption{A metric graph with six bounded edges, three unbounded edges and five vertices}
\label{fig-1}
\end{figure}
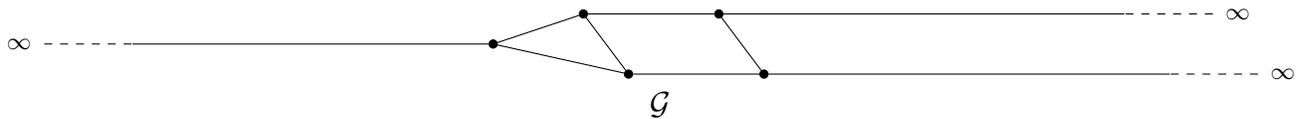

A metric quantum graph $\GG$ is a connected network made up of edges, i.e. bounded or unbounded segments joined together at their endpoints, named vertices. On this structure of edges and vertices, we will give a distance, so that the graph becomes a metric space, which explains the name of \emph{metric graph}. Thus, a metric graph is a one dimensional structure (see Figure \ref{fig-1}), which however should not thought to be embedded in the plane as the angles between edges do not play any role in the theory. 

Let us indicate the set of edges of the graph with $E=\{e_j\}$ and the set of vertices of the graph with $V=\{v_k\}$. We will assume throughout this review that the cardinalities $|E|$ of $E$ and $|V|$ of $V$  are finite. 

We assume that any vertex is a finite endpoint of a finite number of edges. This number is called the degree of the vertex and is denoted by $d_v$. With the symbol $e \prec v$ we mean that edge $e$ is ingoing to $v$ or outgoing from $v$ according to the edge orientation, and $E_v=\{e : \; e\prec v \}$ is the set of all edges ingoing to or outgoing from $v$. The special situation with $d_v = 1$ gives the pendant vertex $v$, while vertices $v$ with $d_v = 2$ can often be treated as the dummy vertices (interior points) with two edges adjacent to $v$ concatenated together. 

An edge $e\in E$ has a length $\l_e$ that can be finite or infinite. The coordinate system on the metric graph is largely arbitrary. When the length $\l_e$ is finite, the edge $e \in E$ can be identified a different choice for the interval $I_e$. For example $[0, \l_e]$ or $[-\l_e/2, +\l_e/2]$ are both useful identifications for the edge $e \in E$. When the length $\l_e$ is infinite, the edge can be identified with any half-line $I_e=[a,\infty)$, where $a \in \mathbb{R}$ is at our disposal. 

A choice of the edge orientation is conventional and all choices are equivalent. A graph is undirected when no choice of orientation is made and directed in the opposite case. We always consider directed graphs.

A point $\x$ on the graph can be identified giving the edge $e$ and the coordinate $x$ on the edge: $\x=(e,x).$
Having identified edges with given intervals and chosen coordinates on these intervals, the length of a path on the graph has a well defined meaning. If $\x,  \y \in \GG$ the distance $d(\x, \y) $ is the infimum of the length of the paths connecting the two points.
This makes $(\GG, d)$ a metric space. The metric graph $\GG$ is compact if and only if $l_e <\infty$  for all $ e \in E$.

\begin{figure}[htb]
	\begin{center}
		\begin{tikzpicture}[xscale= .60,yscale=.4]		
		\node at (-11.5,0) [infinito] (GS) {\scriptsize $\infty$};
		\node at (-9,0) [infinito] (GS) {};
		\node at (-1,0) [nodo] (G1) {};
		\node at (1,1) [nodo] (G2) {};
		\node at (2,-1) [nodo](G3) {};
		\node at (4,1) [nodo](G4) {};
		\node at (5,-1)[nodo] (G5) {};
		\node at (15.5,1) [infinito] (GD2) {\scriptsize $\infty$};
		\node at (16.5,-1)  [infinito] (GD1) {\scriptsize $\infty$};
		
		\node at (14,-1)  [infinito] (GD1) {};
		\node at (13,1)  [infinito] (GD2) {};
		
		\node at (2.7,-2)  [infinito] (G) {$ \GG$};
		\node at (3.2,4.6)  [infinito] (G) {$\blue \Psi$};
		
		\draw[-] (GS)--(G1);
		\draw[-] (G1)--(G2);
		\draw[-] (G1)--(G3);
		\draw[-] (G2)--(G3);
		\draw[-] (G3)--(G5);
		\draw[-] (G2)--(G4);
		\draw[-] (G4)--(G5);
		\draw[-] (G4)--(GD2);
		\draw[-] (G5)--(GD1);
		\draw[dashed] (GD1)--(16,-1);
		\draw[dashed] (GD2)--(15,1);
		\draw[dashed] (GS)--(-11,0);

		\node[-] at (-1,4) [infinito] (F1) {};
		\node[-] at (1,4) [infinito] (F2) {};
		\node[-] at (2,2.5) [infinito] (F3) {};
		\node[-] at (4,3.5) [infinito] (F4) {};
		\node[-] at (5,1.5) [infinito] (F5) {};
		\draw[-,blue,thick] (-8.55,.2) to [out= 5, in = 130] (F1);
		\draw[-,blue,thick] (F1)  to [out= 0, in = 160]  (F2);
		\draw[-,blue,thick] (F1)  to [out= -30, in = 160] (F3);
		\draw[-,blue,thick] (F2)  to [out= -30, in = 120] (F3);
		\draw[-,blue,thick] (F2)  to [out= 40, in = 200] (F4);
		\draw[-,blue,thick] (F3)  to [out= -30, in = 215] (F5);
		\draw[-,blue,thick] (F4)  to [out= -30, in = 195] (F5);
		\draw[-,blue,thick] (F4)  to [out= -25, in = 180] (12.55,1.2);
		\draw[-,blue,thick] (F5)  to [out= -45, in = 210] (13.55,-.8);
		
		\end{tikzpicture}
\end{center}
		\caption{A real function on a metric graph}
	\label{fig-2}
\end{figure}
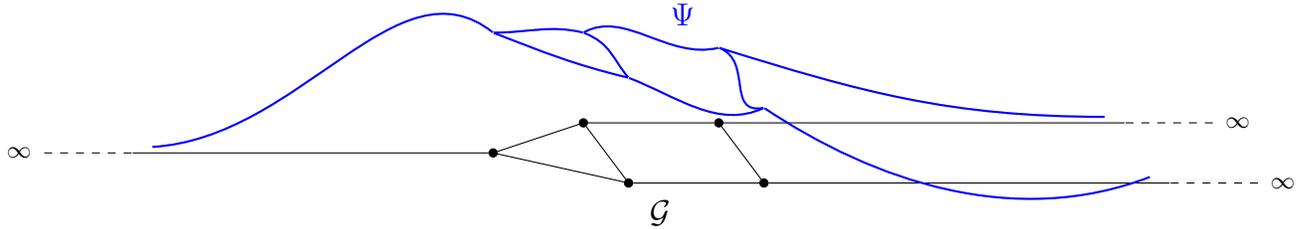

\subsection{Functions and functional spaces}

A function on the metric graph is defined component-wise, giving its value on any edge (see Figure \ref{fig-2} for a real function on the same metric graph 
as on Figure \ref{fig-1}). To assign $\Psi$ on $\GG$ means to assign its edge components $\{ \psi_e \}_{e\in E}$ with $\psi_e : I_e \to \CO$. Equivalently, if $\x = (e,x)$ we set $\Psi(\x) = \psi_e(x)$.

Since a generic metric graph $\GG$ is identified with families of points and intervals, a natural choice for a metric space on $\GG$ is the Lebesgue spaces on the given intervals. In particular, if a function $\Psi : \GG \to \mathbb C$ is integrable, then every component $\psi_e$ of $\Psi$ is integrable on $I_e$ of the edge $e$, and
$$ 
\int_\GG \Psi  \ : = \ \sum_{e \in E} \int_{I_e} \psi_e (x) \, dx.
$$
The space of the $p$-integrable functions $L^p (\GG)$, where $1\leq p\leq \infty$, is the set of the functions $\Psi : \GG \to \mathbb C$ such that $\psi_e \in L^p (I_e)$ for all $e\in E$. Moreover, the norm in $L^p(\GG)$ is defined by 
\begin{align*}
\| \Psi \|_{L^p(\GG)} = \left( \sum_{e\in E} \|\psi_e \|^p_{L^p (I_e) } \right)^{\frac{1}{p}}, \ \ \ \ \  1\leq p<\infty
\end{align*}
and
\begin{align*}
\|\Psi \|_{L^{\infty}(\GG)} = \max_{e\in E}  \|\psi_e \|_{L^\infty (I_e) }, \ \ \ \ \ \ p=\infty.
\end{align*}
These Lebesgue spaces are all Banach spaces (complete normed vector spaces) and they retain the classical properties of these spaces. In particular, $L^2(\GG)$ is a Hilbert space and we will denote by $\langle \cdot ,\, \cdot \rangle$ its inner product
\begin{align*}
\langle \Psi_1,\Psi_2 \rangle =\int_\GG \overline \Psi_1\Psi_2,
\end{align*}
such that $\| \Psi \|_{L^2(\GG)} = \sqrt{\langle \Psi,\Psi \rangle}$.

We denote by $C(\GG)$ the set of continuous functions on $\GG$. This means in particular that the values of the functions in $C(\GG)$ are well defined at the vertex and they coincide with the limit of the functions defined on the edges adjacent to that vertex. It is not especially useful to define spaces $C^k(\GG)$ of the $k$-times continuously differentiable functions because if the functions are continuous across the vertices of $\GG$, their derivatives are generally discontinuous. Instead, Sobolev spaces of weekly differentiable functions that are square integrable with their derivatives play a more important role. The first Sobolev space is especially important:
\begin{align*} 
H^1(\GG) := \left\{\Psi\in L^2(\GG) : \quad \psi_e \in H^1(I_e),  \quad \forall e \in E\right\}
\end{align*}
endowed with the norm 
\begin{align*}
\| \Psi \|_{H^1(\GG)} \ = \ \sqrt{\sum_{e\in E} \| \psi_e \|_{H^1(I_e)}^2}
\end{align*}
In order to incorporate the continuity condition at the vertices, 
we also use 
$$
H_C^1(\GG) := H^1(\GG) \cap C(\GG). 
$$
Note that continuity in the interior of edges of the graph is guaranteed by the one-dimensional Sobolev embedding of $H^1(I_e)$ into $C(I_e) \cap L^{\infty}(I_e)$ expressed by the Sobolev inequality
\begin{align}\label{immersion}
\| \Psi \|_{L^\infty(\GG)} \ \leq \ C_{\infty} 
\| \Psi \|_{H^1(\GG)} \qquad \forall \Psi \in H^1(\GG).
\end{align} 

Analogously we introduce the Sobolev spaces of higher regularity for $k \in \mathbb{N}$:
\begin{align*} 
H^k(\GG) := \left\{\Psi\in L^2(\GG) : \quad  \psi_e \in H^k(I_e), \quad \forall e \in E\right\}
\end{align*}
equipped with the norm 
\begin{align*} 
\| \Psi \|_{H^k(\GG)} \ = \ \sqrt{\sum_{e\in E} \| \psi_e \|_{H^k(I_e)}^2}.
\end{align*}
We will make use of Sobolev spaces in all examples throughout this review. Properties of Sobolev spaces can often expressed by 
inequalities. The following proposition gives the Gagliardo--Nirenberg inequality which relates the $L^q(\GG)$ norm for $q > 2$ to the $H^1(\GG)$ norm.
 
 \begin{proposition}
 	Let $  {\GG}$ be connected and non-compact with finitely many edges. For every $2 \leq q \leq \infty$ there exists a constant $C_q >0$ that depends on $q$ such that
\begin{align}\label{gagliardo2}
\| \Psi \|_{L^q(\GG)} \ \leq \ C_q \| \Psi \|_{L^2(\mathcal G)}^{\f 1 2 + \f 1 q}
\| \Psi' \|_{L^2(\GG)}^{\f 1 2 - \f 1 q} \qquad \forall \Psi \in H^1(\GG).
\end{align}
\label{prop-GN}
\end{proposition}

\begin{remark}
Inequality (\ref{gagliardo2}) cannot hold for a compact graph (take a constant function as a counterexample) and in that case it can be replaced by
\begin{align}\label{gagliardo3}
\| \Psi \|_{L^q(\GG)} \ \leq \ C_q \| \Psi \|_{L^2(\mathcal G)}^{\f 1 2 + \f 1 q}
\| \Psi \|_{H^1(\GG)}^{\f 1 2 - \f 1 q} \qquad \forall \Psi \in H^1(\GG).
\end{align}
which is weaker than \eqref{gagliardo2} but it holds for any metric graph, compact or non-compact. A proof of inequality \eqref{gagliardo2} is in \cite{AST16}, while the proof of \eqref{gagliardo3} for compact graphs is in \cite{Mugnolo}. A generalization of both inequalities is discussed in \cite{CFN17}.
\end{remark}

The kind of property expressed by an inequality as $\|\Psi\|_Y\leq c\|\Psi\|_X \  \Psi\in X $ is named in functional analysis as \emph{continuous embedding} of the space $X$ into the space $Y$. It implies that any convergent sequence of functions of $X$ actually converges in the norm of $Y$ as well.

Another kind of embedding of the space $X$ into the space $Y$ is the \emph{compact embedding}. This means that a sequence of elements of $X$ which is uniformly bounded in the norm of $X$ admits a subsequence that converges in the norm of $Y$.
In other words, if $X$ is compactly embedded into $Y$, then 
\begin{align*}
\{\Psi_n\}_{n \in \mathbb{N}} \subset X,\  \|\Psi_n\|_X\leq C \; \Longrightarrow \; \exists \{\Psi_{n_k}\}_{k \in \mathbb{N}} \subset\{\Psi_n\}_{n \in \mathbb{N}},\  \Psi_\infty\in X : \; \lim_{k\to\infty} \|\Psi_{n_k}-\Psi_\infty\|_Y=0.
\end{align*} 
For our purposes we only need to know that $H^1(\GG)$ is compactly embedded in $L^p(\GG)$ for any $p$ if $\GG$ is a compact metric graph, but that the compact embedding fails if the metric graph is not compact. The compact embedding is very useful when searching minima of functionals in variational calculus. It is often easy to bound the $H^1(\GG)$ norm of a minimizing sequence of functions. This does not provide a converging subsequence as in the finite dimensional setting, but at least gives a converging subsequence in the weaker norm of $L^p(\GG)$ which is often a fundamental step in showing convergence of a minimizing sequence to the minimum of the given functional. On the other hand, when the compact embedding fails, e.g., in the case of unbounded metric graphs, it is a much more difficult issue to manage convergence of a minimizing sequence.

\subsection{Differential operators on metric graphs}

We can now define the differential operators on metric graphs, an essential step to set up mathematical models expressed by differential equations. While there is no difficulty in transferring an operator defined on an interval for any single edge of the metric graph, it is clear that an ambiguity arises as regards the matching conditions at vertices for the functions from different edges connected to the same vertex. Such matching conditions are essential in defining the operator domain, and the choice may depend on physical properties of the modeling problem.

Driven by applications of quantum graphs, the most frequently studied operators are given in the following list:
\begin{itemize}
	\item $\mathcal{H}=-\der{2}{x}$ (Laplace operator);
	\item $\mathcal{H}=-\der{2}{x}+V$ (Schr\"odinger operator);
	\item $\mathcal{H}=-\left( \der{}{x}-iA \right)^2+V$ (Magnetic Schr\"odinger operator).
\end{itemize}
The functions $V$ and $A$ are called respectively electric and magnetic potentials and they act on the domain of $\mathcal H$ by multiplication. In what follows, we will not consider the magnetic Schr\"odinger operator, because it played only a very minor role until now (the only references known to the authors are \cite{cfn15,Hof21}). 

In regards to the Laplace or Schr\"{o}dinger operators, the standard requirement is their self-adjointness in Hilbert space $L^2(\GG)$. As it is well known, 
if $\mathcal{H}$ is self-adjoint, then the  dynamics generated by 
the time-dependent Schr\"{o}dinger equation  $i \partial_t \Psi = \mathcal{H} \Psi$ is well defined and unitary, according to the Stone theorem. This is also an important step to guarantee well-posedness of the dynamics of the nonlinear Schr\"{o}dinger equation, as we will see in Section \ref{sec-2}.

Let us consider the Laplace operator on a metric graph $\GG$ denoted by  $\Delta_\GG$.  The symmetry of $\Delta_{\GG}$, which is expressed by the condition 
\[
\langle\Phi,\Delta_\GG \Psi\rangle - \langle\Delta_\GG\Phi, \Psi\rangle = 0,
\]
translates into constraints on the values of the functions $\Psi$, $\Phi$ and their derivatives at vertices. Integrating by parts and taking into account the orientation of edges at vertices yield the constraints in the form:
\begin{equation*}
\sum_{b_j}\left(\overline{\phi_{e_j}(b_j)}\psi_{e_j}'(b_j)- \overline{\phi_{e_j}'(b_j)}\psi_{e_j}(b_j)\right) - 
\sum_{a_j}\left(\overline{\phi_{e_j}(a_j)}\psi_{e_j}'(a_j)- \overline{\phi_{e_j}'(a_j)}\psi_{e_j}(a_j)\right) =0,
\end{equation*}
where the edge $e_j$ is parameterized by $[a_j,b_j]$ and the summation is taken over all possible $a_j$ and $b_j$. This constraint is not easily readable and does not specify a practical condition to check if a given boundary condition is self-adjoint or not. 

An explicit classification of self-adjoint boundary conditions is well known under various forms. For separated boundary conditions at each vertex, the self-adjoint boundary conditions are identified in Theorem 1.4.4 of \cite{BK13} which is reproduced here.

\begin{proposition}
	\label{prop-self-adjoint}
For any vertex $v$, we assume that 
$$ 
A_v F(v) + B_v F'(v)=0,
$$
where $A_v$ and $B_v$ are complex-valued $d_v\times d_v$ matrices, 
$$
F(v) :=(f_1(v),\dots,f_{d_v}(v)), \quad F'(v):=(f_1'(v),\dots,f_{d_v}'(v)),
$$
and we have used the convention that the derivatives are taken in the outgoing direction from the vertex $v$. The Laplacian $\Delta_{\GG}$ is self-adjoint in $L^2(\GG)$ if 
\begin{enumerate}
\item $\rank (A_v,B_v)=d_v$; 
\item the matrix $A_vB_v^*$ is self-adjoint such that $A_vB_v^*=B_vA_v^*$;
\item $\psi_e\in H^2(I_e)\  \forall e\in \GG$.
\end{enumerate}
\end{proposition}

For every fixed choice of $A_v$ and $B_v$ at $v \in \GG$, the conditions 
of Proposition \ref{prop-self-adjoint} define the domain $\mathcal D (\Delta_\GG)$ of the unique self-adjoint operator $\Delta_\GG$ in $L^2(\GG)$. If $\Delta_{\GG}$ with the dense domain $\mathcal D (\Delta_\GG)$ is self-adjoint in $L^2(\GG)$, then  the dynamics of the time-dependent Schr\"odinger equation 
\[
i \frac d {dt} \Psi_t = - \Delta_{\GG}\Psi_t 
\qquad\qquad
\]
is well defined and the evolution operator $U(t)=e^{it\Delta_{\GG}}$ can be extended to a unitary operator in $L^2(\GG)$.
From the physical point of view, the self-adjoint boundary conditions guarantee  conservation of the probability current through the vertex. On the other hand, it is not easy and is not fully understood how to select the most suitable boundary conditions in specific physically relevant cases. 

Let us now review the most common examples of the self-adjoint boundary conditions.  

We say that a function $\Psi\in H^2(\GG)$ satisfies the \emph{Neumann--Kirchhoff conditions} at vertex $v$ if 
\begin{equation}
\label{NK}
	\left\{ \begin{array}{l}
	\Psi \text{ is continuous at }v \ \ \ (\iff \psi_e(v)=\psi_{e'}(v),\;\; \forall\ e,e' \prec v)\\
	\sum\limits_{e\prec v}\der{}{x_e}\Psi(v)=0 \ \ \ \ \ \ \ \ \ (\iff \sum\limits_{e \prec v} \frac{d \psi_e}{d x_e}(v)=0)	\end{array} \right. 
\end{equation}
	
\begin{remark} 
We always assume that if the edge is outgoing from the vertex, then the derivative is taken with the positive sign, and if the edge is incoming to the vertex, then it is taken with the negative sign. Thus, the second Neumann--Kirchhoff condition is written more explicitly as
$$ 
\sum_{e \prec v} \frac{d \psi_e}{d x_e}(v) =  \sum_{e \leftarrow v} \psi_e' (v)  - \sum_{e \rightarrow v}
 \psi_e' (v) = 0,  
$$
so that the global outgoing derivative vanishes at every vertex.
\end{remark}

A generalization of the Neumann-Kirchhoff condition is given by the so called  \emph{$\delta$-type conditions} \cite{[ACFN14],acfn-aihp,ACFN16,PG19}, for which the function $\Psi\in H^2(\GG)$ at vertex $v$ satisfies 
\begin{equation}
\label{BC-delta} 
\left\{ \begin{array}{l}
\psi_e(v)=\psi_{e'}(v),\ \quad \forall\ e,e' \prec v, \\
\sum\limits_{e\prec v}\der{}{x_e}\psi_e(v)=\alpha_v \psi_{e'}(v), \end{array} \right. 
\end{equation}
where the number $\alpha_v\in\mathbb{C}$ is called the \emph{strength} of the vertex. When the degree of a vertex $v$ is $d_v = 2$, these conditions correspond to the well-known $\delta$-interactions in quantum mechanics. If $\alpha_v=0$ we recover the Neumann-Kirchhoff conditions (\ref{NK}).

There exist boundary conditions such that the function in the operator domain is not continuous. A first example is given by the \emph{$\delta'$-type conditions} \cite{PG20,G22}, for which $\Psi \in H^2(\GG)$ at vertex $v$ satisfies
\begin{equation}
\label{BC-delta-prime}
\left\{\begin{array}{l}
 \der{}{x_e}\psi_e(v) =  \der{}{x_e'} \psi_{e'}(v),\ \quad \forall\ e,e' \prec v, \\
\sum\limits_{e\prec v}\psi_e(v)=\beta_v\der{}{x_e}\psi_{e}(v).
\end{array}\right. 
\end{equation}
Another example is given by the ``generalized Kirchhoff conditions" 
\cite{KP2,KGP,sabirov_ea13,sobirov_ea10}, for which $\Psi\in H^2(\GG)$ at vertex $v$ satisfies
\begin{equation}
\label{NK-generalized}
	\left\{ \begin{array}{l}
\alpha_e \psi_e(v) = \alpha_{e'} \psi_{e'}(v),\ \quad \forall\ e,e' \prec v, \\
	\sum\limits_{e \prec v}  \frac{1}{\alpha_e} \frac{d \psi_e}{d x_e}(v)=0,
	\end{array} \right. 
\end{equation}
where the set of positive constants $\{ \alpha_e \}_{e \in E}$ is given.

We end this section by introducing the energy associated with each boundary condition. Recall that if a self-adjoint operator $A$ in $L^2(\GG)$ is bounded from below, then there exists a quadratic form given by the map $\psi\mapsto \langle \psi, A\psi\rangle$, defined on the dense domain $\mathcal E(A) \subset L^2(\GG)$. This quadratic form is also bounded from below, and there exists $a\geq 0$ such that $\langle \psi, (A +a)\psi\rangle $ is positive. As a result, $\langle\psi,\psi\rangle_A:=\langle \psi, (A +a)\psi\rangle $ is a scalar product on $\mathcal E(A)$. The completion of $\mathcal E(A)$ in the norm induced by $\langle\cdot,\cdot\rangle_A$ is a Hilbert space which is by definition the form domain of the given operator $A$ that we continue to indicate as $\mathcal E(A)$ with an abuse of notation. The quadratic form $E(\psi) := \langle \psi, A \psi \rangle$ is called energy and the form domain $\mathcal E(A)$ is called  the energy space. The explicit expression of the energy $E(\psi)$ is obtained through integration by parts. Specifically for the Laplacian $\Delta_{\GG}$, we have
\begin{itemize}
\item
Neumann--Kirchhoff conditions:
\begin{equation}
E(\Psi)=\|\Psi'\|^2_{L^2(\GG)}, \qquad \qquad \mathcal{E}(-\Delta_{\GG}) = H^1_C(\GG),
\end{equation}
\item $\delta$-type conditions: 
\begin{equation}E(\Psi) = \| \Psi ' \|^2_{L^2(\GG)} + \sum_{v\in V} \al(v) |\Psi(v)|^2, \qquad \mathcal{E}(-\Delta_{\GG}) =H^1_C(\GG),
\end{equation}
\item $\delta'$-type conditions (consider for simplicity a single vertex)
\begin{equation}
E(\Psi)=\| \Psi ' \|^2_{L^2(\GG)} + \frac{1}{\beta}\left|\sum_{e\prec v}\psi'_e(v)\right|, \qquad \mathcal{E}(-\Delta_{\GG}) =H^1(\GG),
\end{equation}
\item generalized Kirchhoff conditions: 
\begin{equation}E(\Psi)=\|\Psi'\|^2_{L^2(\GG)}, \qquad \mathcal{E}(-\Delta_{\GG}) = \{\Psi\in H^1(\Psi) : \;\; \alpha_e \psi_e(v) = \alpha_{e'} \psi_{e'}(v),\ \quad \forall\ e,e' \prec v \}.
\end{equation}\end{itemize}
Different boundary conditions can be posed at different vertices, with the corresponding changes in the definitions of the operator domain $\mathcal{D}(\Delta_{\GG})$ and the energy space $\mathcal{E}(-\Delta_{\GG})$. 

\subsection{Other differential operators} 

Other differential operators can be defined on a metric graph. Many examples such as the Dirac operator and the Airy operator have been considered in the literature due to their importance for physical applications. As a case study, we consider here the Airy operator 
\begin{equation}
A := \alpha \frac{\partial^3 u}{\partial
	x^3}+\beta\frac{\partial u}{\partial x},
\label{airy}
\end{equation}
where $\alpha\in \mathbb R\setminus\{0\}$ and $\beta \in \mathbb R$ are certain physical constants. This third-order differential operator defines the linear part of the Korteweg--De Vries (KdV) equation, which is the most studied dispersive equation after the NLS equation. It describes the unidirectional propagation in shallow water channels and more generally dynamics of long waves in weakly dispersive systems. The linear KdV equation is $\partial_t u = A u$ can be considered 
on a metric graph $\GG$ with the obvious motivation to understand propagation of shallow water waves in branching channels. 

A first relevant remark that should be taken into account is that if we replace the second-order Schr\"odinger operator with the third-order Airy operator, then the wave propagation becomes sensible to the direction. This forces us to consider directed metric graphs. 

This fact is well-known already on the half-line or on a finite interval \cite{bona-half}. Without loss of generality, we can choose $\alpha<0$. Then, the initial-value problem for the linear KdV equation $\partial_t u = A u$ is well-posed 
on $(0,+\infty)$ if a single boundary condition is imposed at $x = 0$ 
and is well-posed on $(-\infty,0)$ if two boundary conditions are imposed at $x = 0$ \cite{bona-half}. More generally for a star graph, the number and type of boundary conditions at the only vertex required for the well posedness of the linear KdV equation on the star graph are classified in \cite{MNS18}. 

Without reproducing the details of \cite{MNS18}, we mention that the Airy operator $A$ is formally antisymmetric, being a combination of odd derivatives. A proper choice of the domain $\mathcal D(A)$ will make it skew-adjoint and not self-adjoint. This is not an obstacle because the skew-adjoint operators $A$ in $\partial_t u = A u$ generate unitary operators due to the Stone theorem. 

It turns out that only on {\em balanced} metric graphs (i.e. star graphs with an equal number of ingoing and outgoing edges) appropriate boundary conditions exist that make the Airy operator skew-adjoint. The boundary conditions can be expressed requiring that a finite-dimensional map from the vector space of the ingoing boundary values $u(0-),u'(0-),u''(0-)$ to the vector space of the outgoing boundary values $u(0+),u'(0+),u''(0+)$ exists and has to be unitary. Such conditions are completely different from the ones described in Proposition \ref{prop-self-adjoint} for $\Delta_{\GG}$. 

If the initial-value problem for the linear KdV equation is well-posed on a balanced star graph with the appropriate boundary conditions, then the time evolution is unitary and preserves the squared $L^2$-norm $P(u) := \int_{\mathcal G} |u|^2$ of the solution $u$, which in the context of fluid flow is to be interpreted as the momentum. On the other hand, only a subset of the admitted boundary conditions preserve the mass, which is represented by $M(u) := \int_{\mathcal G} u$. 

If the star graph is unbalanced (representing, for example, the confluence of rivers or a branching of the water channels), then the time evolution is not unitary and does not result in the conservation of the momentum $P(u)$. In this case, the time evolution 
of the linear KdV equation is well-posed if and only if the Airy operator $A$ 
is {\em dissipative} or equivalently if it generates a  contraction semigroup, instead of a unitary group. We recall that $A$ is dissipative if and only if 
$$
{\rm Re} \langle Au, u\rangle \leq 0, \quad \mbox{\rm for all } u\in \mathcal{D}(A).
$$ 
Again, this translates into a condition on a finite-dimensional map relating the ingoing and outgoing boundary values of $u$ at the vertex, described explicitly in \cite{MNS18}. As in the conservative case, only a proper subset if these boundary conditions preserves the mass $M(u)$. 

While these general and qualitative results appear suggestive and physically sensible, one should also admit that it is unclear which boundary conditions allowed by the theory have to be used in concrete problems. A reduction from more realistic, two- or three-dimensional models to the problems on the one-dimensional metric graphs is a fundamental step in clarifying this essential question \cite{Beck,Costa,Kuch,Post1,Post2}.

\section{Nonlinear partial differential equations on metric graphs}
\label{sec-2}

Thanks to the possibility of defining operators on metric graphs, we can consider different linear and nonlinear partial differential equations on metric graphs. Among these are the wave, Klein-Gordon, Schr\"odinger, Dirac, and Korteweg--de Vries equations, which are classical examples of the dispersive wave equations of the form
\begin{equation}
\label{dispersive-PDE}
\partial_t u(t)=Lu(t) +  N(u),
\end{equation}
where $L$ is a  differential operator with constant coefficients and $N$ is the nonlinear part. The three examples addressed here are given by  
\begin{align*}
\partial_t u(x,t)&=\ \ \partial_x u(x,t) - \partial^3_x u(x,t) + u(x,t)\partial_x u(x,t)\ \ \ \ \ \ \ \ \ \text{KdV equation}\\
i\partial_t u(x,t)&=-\partial^2_x u(x,t) \pm |u(x,t)|^{2p} u(x,t)\ \ \ \ \ \ \ \ \ \ \ \ \ \ \ \ \text{NLS equation}\\
\partial^2_t u(x,t)&=\ \ \partial^2_x u(x,t) -m^2u(x,t) \pm |u(x,t)|^{2p} u(x,t)\ \ \ \ \ \ \ \ \ \ \text{NLKG equation}
\end{align*}
whereas the Dirac equations can be found in Section \ref{sec-7}.

\subsection{Linear dispersion}

We briefly recall what dispersion means in this context by considering the linear part of the general equation (\ref{dispersive-PDE}) given by 
$\partial_t u(t)=Lu(t)$. If the linear equation is considered on the real line $\mathbb{R}$, we can search for the plane wave solutions, $u(x,t)= e^{i(\xi_0x-\omega_0 t)}$, or more generally for solutions in the Fourier form,
$$
u(x,t)= \int_\RE c(\xi)e^{i(\xi x-\omega(\xi) t)}\ d\xi,
$$
where $c(\xi)$ is the Fourier transform of the initial datum 
$u(x,0) = u_0(x)$. 

If $\widehat{L u}(\xi,t) := -ip(\xi)\hat u(\xi,t)$ under the Fourier transform,
then $\omega(\xi) = p(\xi)$ is called {\em the dispersion relation}. The linear equation $\partial_t u(t)=Lu(t)$ is called dispersive if $\omega''(\xi)\neq 0 $.
For the dispersive equations, the "group velocity" $v_g(\xi)=\omega'(\xi)$ depends on $\xi$, and in general differs form the "phase velocity" of the wave, $v_p(\xi)={\omega(\xi)}/{\xi}$. As a consequence, an initially concentrated wave-packet "disintegrates" or "disperses" with time. Let us inspect the same examples above. 
\begin{itemize}
	\item The linear KdV equation $\partial_t u=\partial_x u - \partial^3_x u$ is dispersive due to $\omega(\xi)=-\xi-\xi^3$ and $\omega''(\xi)=-6\xi$.
	\vskip5pt
	\item The linear Schr\"{o}dinger equation $i\partial_t u =-\partial^2_x u$ is dispersive due to $\omega(\xi)=\xi^2$ and $\omega''(\xi)=2$.
	\vskip5pt
	\item The linear  Klein-Gordon equation $\partial^2_t u =\partial^2_x u - m^2 u$ is dispersive due to $\omega(\xi)=\pm\sqrt{m^2+\xi^2}$ and $\omega''(\xi)\neq 0$. However, $\omega''(\xi)\to 0$ as $|\xi|\to \infty$, hence the short waves are less dispersive compared to the long waves. If $m = 0$, then there is no dispersion for long waves and the equation becomes the one-dimensional wave equation. 
\end{itemize}

Using the Fourier transform and the dispersion relation, one can obtain very precise results about the decay of solutions in time by means of the stationary phase method. 
A different and very powerful tool suitable to measure the decay of solutions in time is given by the so-called dispersive estimates.

Consider for example the linear Schr\"odinger equation $i\partial_t u = -\partial^2_x u$ on the line $\mathbb{R}$. If the initial datum $u_0$ is sufficiently regular and decaying at infinity, then the unique solution is given in the convolution form:
\begin{equation*}
u = U(t)u_0 := \frac{1}{(4\pi i t)^{1/2}}\int_{\RE} e^\frac{{i|x-y|^2}}{4t}u_0(y)\ dy.
\end{equation*}
It follows from this representation that 
\begin{itemize}
	\item The $L^2$ norm conserves in time:
	$$
	\| U(t)u_0 \|_{L^2(\mathbb{R})} = \|u_0 \|_{L^2(\mathbb{R})}.
	$$
	\item The $L^{\infty}$ norm decays in time: 
	\begin{equation*}
	\|U(t) u_0 \|_{L^\infty(\mathbb{R})} := \sup_{x\in \RE}|U(t)u_0| \leq \frac{1}{(4\pi t)^{1/2}}||u_0||_{L^1(\mathbb{R})}.
	\end{equation*}
\end{itemize}
With some additional efforts, one can also 
interpolate the two estimates above to the more general estimate:
\begin{equation*}
\| U(t)u_0 \|_{L^p(\mathbb{R})} \leq \frac{1}{(4\pi t)^{\frac{1}{2}(\frac{1}{q}-\frac{1}{p})}}||u_0||_{L^q(\mathbb{R})} , \ \ \ \ \ \ p\in [2,\infty],\ \ \ \frac{1}{p}+ \frac{1}{q}=1,
\end{equation*}
and to the so-called Strichartz estimates, which involve the space-time norms \cite{Cazenave}. 

A similar toolbox of dispersive estimates exists for the linear KdV and Klein-Gordon equations. The important point is that dispersion, or the decay just illustrated, conflicts with the tendency to concentrate solutions due to the nonlinear part of the general equation (\ref{dispersive-PDE}). This competition is exactly what makes the analysis of nonlinear dispersive equations so interesting, because it makes possible the existence of coherent structures such as standing waves.

The analysis of dispersion on metric graphs is not fully developed until now and it is mostly an open problem. We only mention that dispersive and Strichartz estimates for the Schr\"odinger equations on star graphs and trees were studied in \cite{[ACFN1],Banica2,Grecu18}, while those on the tadpole graph were studied in \cite{AMAN17}. 

\subsection{The NLS equation} 

The most studied nonlinear dispersive wave equation on metric graphs is without any doubt the focusing version of the nonlinear Schr\"{o}dinger (NLS) equation with power nonlinearity, 
\begin{equation}\label{tNLS}
i \frac d {dt} \Psi =  -\Delta_{\GG} \Psi + V \Psi - (p+1)| \Psi |^{2p} \Psi,
\end{equation}
where $V$ is the external potential and $p > 0$ is the power. 
The focusing version of the NLS equation is defined by the minus sign in front of the power nonlinearity, in contrast to the defocusing version defined by the plus sign. It appears that the defocusing NLS equation on metric graphs is less studied than the focusing NLS equation.

One of the main physical applications of the focusing NLS equation is propagation of matter waves (attractive Bose-Einstein condensates) or optical pulses in nonlinear optics in $T$-junctions, $X$-junctions and other similar geometrical settings (see \cite{lorenzo,N14} and references therein).

The power nonlinearity of the NLS equation (\ref{tNLS}) 
can be defined edge by edge. Each component on the metric graph satisfies the following one-dimensional NLS equation on edges of the graph $\GG$,
\begin{align*}
&i \pd{}{t} \psi_{e}(x,t) = -\pd{^2}{x^2}  \psi_e(x,t) + V(x)\psi_e(x,t) - (p+1)| \psi_e(x,t)|^{2p}  \psi_e(x,t), \quad
x \in e,
\end{align*}
closed with the corresponding boundary conditions at vertices. 
The problem consists in a system of identical equations, as many as there exist edges of the graph, coupled through the boundary conditions at the vertices $\{ v \}_{v \in V}$, where the values of the wave functions $\{ \psi_e \}_{e \in E}$ and their first derivatives are related.

The first issue to deal with is the time evolution of the NLS equation (\ref{tNLS}). As it is well known, the classical analysis of the existence and uniqueness of the initial-value problem is difficult due to the presence of the unbounded operator $\Delta_{\GG}$. To overcome the difficulty, one can use good mapping properties of the propagator $U(t)=e^{it\Delta_\GG}$, which is a unitary operator in $L^2(\GG)$. Using Duhamel's formula, one can rewrite 
the initial-value problem for the NLS equation (\ref{tNLS}) 
as the integral equation
\begin{equation}\label{duhamel}
\Psi(t) = U(t)\Psi_0 - i \int_0^t U(t-s) V \Psi(s)\ ds + i(p+1) \int_0^t U(t-s)| \Psi(s) |^{2p} \Psi(s)\ ds,
\end{equation}
where $\Psi_0 = \Psi(0)$. The NLS equation (\ref{tNLS}) follows from the integral equation (\ref{duhamel}) by formal differentiation in $t$. The validity of the Duhamel's formula can be  shown exactly as in the standard case of the NLS equation on $\mathbb{R}$ (see \cite{Cazenave}).

Since $U(t)$ is unitary in $L^2(\GG)$, one can use iterations of the Banach fixed point theorem and show local well posedness of the initial-value problem in different functional spaces: $L^2(\GG)$, energy space $\mathcal{E}(-\Delta_{\GG})$, and the operator domain $\mathcal{D}(\Delta_\GG)$.

\begin{figure}[htb]
\begin{center}
\begin{tikzpicture}
\node at (-6.5,0) [infinito](0) {${\infty}$};
\node at (-2,0) [nodo,black] (5) {};
\node at (0,0) [nodo,black] (4) {};
\node at (-2,2) [nodo,black] (2) {};
\node at (0,2) [nodo,black] (3) {};
\node at (4,0) [infinito] (8) {${\infty}$};
\node at (4,2) [infinito] (7) {${\infty}$};
\node at (2,2) [nodo,black] (6) {};
\node at (-3.8,2) [nodo,blue] (1) {};
\node at (-2.2,2.2){\black{$v_2$}};
\node at (0.1,2.2){\black{$v_3$}};
\node at (0.2,-0.2){\black{$v_4$}};
\node at (-2.2,-0.2){\black{$v_5$}};
\node at (2,2.2){\black{$v_6$}};

\node at (-3.6,2.2){\black{$v_1$}};

\draw [-,black] (-4.2,2) circle (0.4cm) ;
\draw [-,black] (5) -- (2);
\draw [-,black] (0) -- (5);
\draw [-,black] (5) -- (4);
\draw [-,black] (5) to [out=-40,in=-140] (4);
\draw [-,black] (5) to [out=130,in=-130] (2);
\draw [-,black] (5) -- (3);
\draw [-,black] (1) -- (2);
\draw [-,black] (2) -- (3);
\draw [-,black] (3) -- (6);
\draw [-,black] (4) -- (3);
\draw [-,black] (6) -- (7);
\draw [-,black] (4) -- (8);

\node at (2,0.5) [nodo,black] (10) {};
\node at (2,0.3){\black{$v_7$}};
\draw [-,black] (3) -- (10);

\end{tikzpicture}
\end{center}
\caption{A metric graph with a compact core and a finite number of half-lines.}
\label{starlike}
\end{figure}
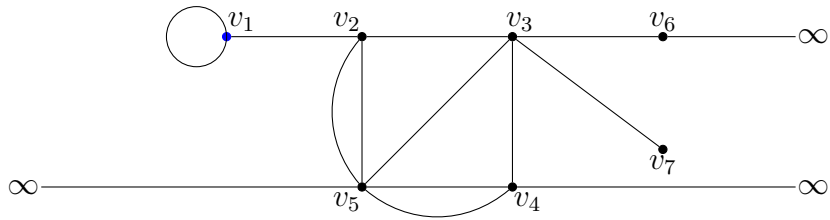

A prototypical and fairly general case where this strategy has been pursued is the case of NLS on a metric graph with a compact core and a finite set of half-lines (see Figure \ref{starlike}). Moreover, vertices with $\delta$-type conditions and the external potentials $V$ on the edges can be added to the consideration. The following proposition gives local well-posedness of the corresponding initial-value problem in the energy space $\mathcal{E}(-\Delta_{\GG}) = H^1_C(\GG)$ (see \cite{CFN17} for the proof).

\begin{proposition}
\label{prop-local}
Let $\GG$ be a graph with a compact core and $V = V_+ - V_-$, where $V_\pm\geq 0$ satisfy $V_+\in L^1(\GG) + L^{\infty}(\GG)$ and $V_-\in L^r(\GG)$ for some $r\in [1,1+1/p]$ and $p>0$. Then, 
\begin{itemize} 
	\item For any $\Psi_0 \in H^1_C(\GG)$, there exists $T > 0$ such that 
the integral equation \eqref{duhamel}  has a unique solution $\Psi \in C ([0,T),
H^1_C(\GG)) \cap C^1 ([0,T), H^1_C(\GG)^\star)$, where $H^1_C(\GG)^\star$ is the dual of $H^1_C(\GG)$.
\item Moreover, if $T^\star$ is the maximal existence time, then the following ``blow-up alternative'' holds: either $T^\star = \infty$ or
$\lim\limits_{t \to T^\star} \| \Psi(t) \|_{H^1(\GG)} = \infty$.
\end{itemize}
\end{proposition}

\begin{remark} The meaning of the blow-up alternative is that there are only two possibilities: existence for all times or explosion of the $H^1$ norm of the solution in a finite blow-up time $T^\star$. 
\end{remark}

Thanks to the blow-up alternative, one can obtain global well posedness 
of the NLS equation (\ref{tNLS}) in $H^1_C(\GG)$ if the $H^1(\GG)$ norm of the solution can be controlled in the time evolution. To this end, a typical strategy is to exploit conservation laws,  since the NLS equation (\ref{tNLS}) admits the following conservation laws of mass and energy.

\begin{proposition}
	\label{prop-laws}
For any solution $\Psi \in C([0,T), H^1_C(\GG))
\cap C^1 ([0,T), H^1_C(\GG)^\star)$ to the NLS equation (\ref{tNLS}) with any $p > 0$ the following conservation laws hold at
any time $t \in [0,T)$:
\begin{align*}
M(\Psi(t)) \ &=\left\|\Psi(t) \right\|_{L^2(\GG)}^2 = M(\Psi_0), \\
E(\Psi(t)) \ &=  \| \Psi '(t) \|_{L^2(\GG)}^2 + \langle \Psi(t),V\Psi(t) \rangle+  \sum_{v\in V} \al(v) |\Psi(v,t)|^2  -\|\Psi(t) \|_{L^{2p+2}(\GG)}^{2p + 2} = E(\Psi_0).
\end{align*}
\end{proposition}

By using Gagliardo--Nirenberg inequality of Proposition \ref{prop-GN} 
and the mass and energy conservation of Proposition \ref{prop-laws}, 
one can obtain \emph{a priori} bounds on the $H^1_C(\GG)$ norm of the 
local solution of Proposition \ref{prop-local} 
assuring global well posedness of the initial-value problem. 

\begin{proposition}
	\label{t:wp}
For any $\Psi_0 \in H^1_C(\GG)$,  the
NLS-equation \eqref{tNLS} with $0 < p < 2$ 
has a unique solution $\Psi \in C([0,\infty),H^1_C(\GG))
\cap C^1 ([0,\infty), H^1_C(\GG)^\star)$. 
\end{proposition}

The range of power nonlinearities $p\in(0,2)$ where the global well-posedness 
of Theorem \ref{t:wp} holds is called \emph{sub-critical}. When $p=2$, solutions can again exist for every time but only for the initial data with small mass $M(\Psi_0)$. The case $p=2$ where the global existence is subordinated by the mass of the initial datum is called \emph{critical}. The case $p > 2$ is \emph{supercritical} and in that case one can expect explosion of the $H^1(\GG)$ norm of the solution in a finite time even for small mass $M(\Psi_0)$ and for negative energy $E(\Psi_0)$ \cite{Cazenave}. The finite-time blow-up has not been proved for the NLS equation on a general metric graph $\GG$ but only for the simplified case of star graphs (see \cite{GO19}).

\begin{remark}
A similar analysis can be carried out for more regular initial data in the operator domain $\mathcal{D}(\Delta_{\GG})$ (see again \cite{CFN17,GO19}) obtaining solutions with values in the operator domain. Global well posedness is identical to the case of data in the energy space $\mathcal{E}(-\Delta_{\GG}) = H^1_C(\GG)$.
\end{remark}

\begin{remark} There are fewer works on the well posedness of the NLS equation \eqref{tNLS} with boundary conditions different from the $\delta$-type conditions and in $L^2(\GG)$. For the special case of star graphs, well-posedness of the NLS equation in the case of $\delta'$-type conditions was considered in  \cite{G22}, while dispersive and Strichartz estimates have been used in \cite{Grecu18} to prove well-posedness of the NLS equation 
with some general boundary conditions in $L^2(\GG)$.
\end{remark}

\section{Standing waves of the NLS equation}

Consider the focusing NLS equation (\ref{tNLS}) on the metric graph $\GG$.
A standing wave of the NLS equation is a solution of the form 
$\Psi(x,t) = e^{-i\omega t}\Phi_\omega(x)$, 
where the parameter $\omega$ is often called frequency of the standing wave. The profile $\Phi_\omega$ satisfies the stationary NLS equation
\begin{equation}\label{stationary}
 - \Delta_{\GG} \Phi_\omega   + V\Phi_\omega- (p+1)| \Phi_\omega  |^{2p} \Phi_\omega=\omega \Phi_\omega.
 \end{equation}
The stationary NLS equation is a system of ordinary differential equations coupled through boundary conditions at the vertices of $\GG$. 

\subsection{Ground state and bound states}

Recall the conservation laws of mass and energy of the NLS equation 
in Proposition \ref{prop-laws}. It is obvious that these quantities computed 
at the standing waves are independent of time. 
Correspondingly, a class of standing waves can be defined by means of a variational principle.

A {\em ground state} of the NLS equation is defined as the minimizer of the energy $E(\Psi)$ at fixed mass $M(\Psi) = \mu$. More formally, $\Phi \in H^1_C(\GG)$ is a ground state if $E(\Phi) = \mathcal{E}_{\mu}$, where 
\begin{equation}\label{minimization}
\mathcal{E}_{\mu} = \inf_{\Psi \in H^1_C(\GG)} \{ E(\Psi) : \quad M(\Psi) =\mu \}.
\end{equation}
This definition assumes that the following two requirements are satisfied:
\begin{itemize}
\item $\mathcal{E}_{\mu} > - \infty$
\item $\mathcal{E}_{\mu}$ is attained at some $\Phi\in H^1_C(\GG)$.
\end{itemize}
The definition of the ground state is completely analogous to the one used for quantum systems in the framework of the Schr\"{o}dinger theory of Quantum Mechanics.

A \emph{minimizing sequence} is any sequence $\{\Psi_n\}_{n \in \mathbb{N}} \in H^1_C(\GG)$ such that $\lim\limits_{n\to \infty}{E(\Psi_n)} = \mathcal{E}_{\mu}$. 
If  $\mathcal{E}_{\mu} > - \infty$ and the infimum is attained at the minimizer 
$\Phi \in H^1_C(\GG)$, then the minimizing sequence converges to the minimizer. However, if  $\mathcal{E}_{\mu} > - \infty$  but the infimum is not attained, then the minimizing sequence may converge to some trivial function in $H^1_C(\GG)$ due to either splitting 
or vanishing, according to the concentration-compactness principle of P.L. Lions \cite{Leon1,Leon2}.

The existence and variational characterization of the ground state of the NLS equation have been considered only for the Neumann--Kirchhoff conditions and for the $\delta$-type conditions. If the ground state is attained, it admits the following properties \cite{AST15}.

\begin{proposition} 
Assume that $\Phi \in H^1_C(\GG)$ is a minimizer of the variational problem (\ref{minimization}). Then we have: 
		\begin{enumerate}
			\item There exists $\omega \in \mathbb{R}$ such that $\Phi$ is a weak solution of the stationary NLS equation (\ref{stationary}).
			\item For every vertex $v$, $\Phi$ satisfies  
			\begin{equation}
				\sum_{v \prec e} \frac{d \phi_e}{d x_e}(v) = \alpha(v)\phi_e(v)
				\label{eq: kirchhoff}
			\end{equation}
			\item \label{positivity} $\Phi$ is a classical solution on every edge of the graph $\GG$ and up to the change of sign, we may assume that $\Phi$ is real-valued and positive on $\GG$.
		\end{enumerate}
		\label{prop-EL}
	\end{proposition}

Besides the problem of controlling if the infimum $\mathcal{E}_{\mu}$ is attained 
for a given mass $\mu$, another problem is to relate the mass parameter 
$\mu$ in (\ref{minimization}) with the frequency parameter $\omega$ in (\ref{stationary}). Standing wave solutions of the stationary NLS equation 
(\ref{stationary}) for fixed $\omega$ may also appear as local minimizers or saddle points of the constrained energy $E(\Psi)$. If they are not global  minimizers, they are usually refer to as {\em bound states}.

\subsection{Stability of standing waves}

If $\Phi \in H^1_C(\GG)$ is the ground state of the variational problem 
(\ref{minimization}), it is expected to be stable in the time evolution 
of the NLS equation (\ref{tNLS}). However, if $\Phi \in H^1_C(\GG)$ is only the bound state, stability of this standing wave needs to be analyzed by 
different tools, e.g., by linearizing of the NLS equation at the neighborhood of $\Phi$.

Stability of standing waves is very similar to stability of equilibrium points or periodic orbits in finite-dimensional dynamical systems. 
In the particular case of the NLS dynamics, stability of standing waves 
needs to be understood in the orbital sense, as we explain next.

Recall that the NLS equation (\ref{tNLS}) on a metric graph $\GG$ is invariant with respect to the action of the group 
\begin{equation}
U(1) = \left\{ T(\theta) :\ \theta\in \RE\right\}, \quad T(\theta) := e^{-i\theta},
\end{equation} 
meaning that if $\Psi$ is a solution of the NLS equation, $T(\theta)\Psi$ is a solution as well. This symmetry is at the origin of existence of standing waves: a solution of the form $\Psi(x,t) = e^{-i\omega t} \Phi_{\omega}(x)$ is nothing but the action on $\Phi_\omega$ by the element of the group with $\theta=\omega t$, i.e. $\Psi  =T(\omega t) \Phi_\omega$. 
As a result, a standing wave cannot be stable in the usual sense of stability of equilibrium points. This is a general property of dynamical systems with symmetries, well known in the finite-dimensional case, where it gives rise to so called relative equilibria. 

Let us show why it is so. Fix a standing wave $\Psi(x,t) = e^{-i\omega t}\Phi_{\omega}(x)$ and consider a set of frequencies 
$\{ \omega_n \}_{n \in \mathbb{N}}$ such that 
$\omega_n\longrightarrow \omega$ as $n \to \infty$ and the corresponding standing waves $\Psi(x,t) = e^{-i \omega_n t} \Phi_{\omega_n}$. It can be shown that solutions of the stationary NLS equation (\ref{stationary}) are continuous with respect to the parameter $\omega$ such that $\| \Phi_{\omega_n} - \Phi_\omega \|_{H^1(\GG)} \longrightarrow 0$ as $n \to \infty$. Explicitly, for every $\delta>0$ there exists $N_\delta\in \mathbb N$ such that
for every $n\geq N_\delta$, we have 
$$
\|\Phi_{\omega_n} - \Phi_\omega \|_{H^1(\GG)} \leq \delta.
$$
On the other hand, we have along the time evolution that 
\begin{align*}
\| e^{-i \omega_n t} \Phi_{\omega_n} - e^{-i \omega t} \Phi_{\omega} \|_{H^1(\GG)} = & \|e^{-i\omega_n t}\Phi_{\omega_n} - e^{-i\omega_n t} \Phi_\omega + e^{-i\omega_n t} \Phi_\omega - e^{-i\omega t}\Phi_\omega \|_{H^1(\GG)} \\ \geq & 
\left|\ |e^{-i\omega_n t}-e^{-i\omega t}| \|\Phi_\omega\|_{H^1(\GG)} - \|\Phi_{\omega_n} - \Phi_\omega \|_{H^1(\GG)}\ \right | 
\end{align*}
and taking the supremum in time we obtain for fixed $n\geq N_{\delta}$ that 
\begin{align*}
\sup_{t\geq 0} \| e^{-i \omega_n t} \Phi_{\omega_n} - e^{-i \omega t} \Phi_{\omega} \|_{H^1(\GG)} \geq  2 \|\Phi_\omega\|_{H^1(\GG)} - \delta. 
\end{align*}
This means that solutions cannot be made as near as we want so that the given standing wave is not stable in the usual sense. Failure of stability is due to the fact that the phases $\theta = \omega_nt $ and $\theta = \omega t$ of the standing waves can assume any value drifting along the orbit of bound states that are initially nearby. 

Similar arguments exclude stability in the presence of other symmetry groups. The consequence of this fact is the need to update both the definition of stability and  of the criteria useful to assure stability of the standing waves. 
For the definition, the key is that instead of comparing the values of solutions at any fixed time, we compare the whole orbits and choose the nearest point along the orbit. We arrive to the notion of {\em orbital stability} and {\em orbital instability}. 

\begin{definition}
	\label{def-stab}
	The standing wave $\Psi(x,t) = e^{-i\omega t}\Phi_{\omega}(x)$ is orbitally stable if for every $\varepsilon > 0$ there exists $\delta > 0$ such that for every $\Psi_0 \in H^1_C(\GG)$ satisfying $\| \Psi_0 - \Phi_{\omega} \|_{H^1(\GG)} < \delta$, the unique solution $\Psi \in C([0,\infty),H^1_{C}(\GG))$ of the NLS equation (\ref{tNLS}) with $\Psi(\cdot,0) = \Psi_0$ satisfies 
	$$
	\inf_{\theta \in \mathbb{R}} \| \Psi(\cdot,t) - e^{-i \theta} \Phi_{\omega}  \|_{H^1(\GG)} < \varepsilon, \quad t > 0.
	$$
\end{definition}

\begin{definition}
	\label{def-unstable}
	The standing wave $\Psi(x,t) = e^{-i\omega t}\Phi_{\omega}(x)$ is orbitally unstable if there exists $\varepsilon > 0$ such that for every $\delta > 0$ there exists $\Psi_0 \in H^1_C(\GG)$ satisfying $\| \Psi_0 - \Phi_{\omega} \|_{H^1(\GG)} < \delta$ such that the unique solution $\Psi \in C([0,T),H^1_{C}(\GG))$ of the NLS equation (\ref{tNLS}) with $\Psi(\cdot,0) = \Psi_0$ and some $T > 0$ satisfies for some $t_0 \in (0,T)$
	$$
	\inf_{\theta \in \mathbb{R}} \| \Psi(\cdot,t_0) -  e^{-i \theta} \Phi_{\omega}  \|_{H^1(\GG)} \geq \varepsilon.
	$$
\end{definition}

Orbital stability in Definition \ref{def-stab} means that if 
the initial datum $\Psi_0  \in H^1_C(\GG)$ is close to the bound state $\Phi_{\omega}$ in the $H^1(\GG)$ norm, the unique solution $\Psi(\cdot,t) \in H^1_C(\GG)$ remain close to the orbit of the bound states $\{ T(\theta) \Phi_{\omega} \}_{\theta \in \mathbb{R}}$ for all $t > 0$. 
Orbital instability in Definition \ref{def-unstable} means that for any small distance from the bound state $\Phi_{\omega}$ in the $H^1(\GG)$ norm, there exists the initial datum $\Psi_0 \in H^1_C(\GG)$ such that the unique solution 
$\Psi(\cdot,t) \in H^1_C(\GG)$ departs from the orbit of the bound states $\{ T(\theta) \Phi_{\omega} \}_{\theta \in \mathbb{R}}$ at some time $t_0$ 
in some fixed distance $\varepsilon$.

The sufficient criteria ensuring orbital stability of solitary waves have been developed during 1970s and 1980s starting with the KdV equation \cite{Benj72,Bona75} and NLS equation \cite{VK}, and continued by  \cite{Grillakis,GSS,W2,W3} that consolidated a rather general and comprehensive theory, now applied in a wide range of examples including the NLS equation on metric graphs.

A Lyapunov function for stability of standing waves in the NLS equation (\ref{tNLS}) is constructed by using conserved quantities in the following action functional 
\begin{equation}
\label{S-def}
S_{\omega}(\Psi) := E(\Psi)-\omega M(\Psi).
\end{equation}
Differentiation yields $S_{\omega}'\Psi=0$, which is just the stationary NLS equation (\ref{stationary}) solved by $\Phi_{\omega}$. The idea for stability analysis is to derive the sufficient conditions of positivity (also called coercivity) of $S_\omega$  near $\Phi_{\omega}$ by using the second derivative test. 

If $\delta$ is small enough in Definition \ref{def-stab}, we set the perturbation $\eta =\alpha+i\beta$, explore conservation laws and phase invariance, and obtain 
\begin{align*}
\delta^2 & > S_{\omega}(\Psi_0)-S_{\omega}(\Phi_{\omega}) \\
& = S_{\omega}(\Psi(\cdot,t))-S_{\omega}(\Phi_{\omega})\\
& = S_{\omega}(e^{i\theta}\Psi(\cdot,t))-S_{\omega}(\Phi_{\omega}) \\
& = S_{\omega}(\Phi_{\omega}+\alpha+i\beta))-S_{\omega}(\Phi_{\omega})\\
& = \langle \mathcal{L}_+ \alpha, \alpha \rangle +  
\langle \mathcal{L}_- \beta,\beta \rangle + R(\alpha,\beta),
\end{align*}
where the second variation of $S_{\omega}$ is defined by two self-adjoint operators $\mathcal{L}_{\pm}$ in $L^2(\GG)$ given by 
\begin{eqnarray}
\label{L-plus}
\mathcal{L}_+ &:=& -\Delta_{\GG} + V(x) - (p+1)(2p+1) |\Phi_{\omega}(x)|^{2p} - \omega, \\
\mathcal{L}_- &:=& -\Delta_{\GG} + V(x) - (p+1) |\Phi_{\omega}(x)|^{2p} - \omega,
\label{L-minus}
\end{eqnarray}
with the dense domain $\mathcal{D}(\Delta_{\GG})$ 
and $R(\alpha,\beta)$ is the super-quadratic remainder term. 

If there existed positive $C_+\ , C_-$ such that  
\begin{equation}
\label{coercivity}
\langle \mathcal{L}_+ \alpha, \alpha \rangle \geq C_+ \|\alpha\|^2_{H^1(\GG)}, \qquad \langle \mathcal{L}_- \beta,\beta \rangle \geq C_-\|\beta\|^2_{H^1(\GG)}
\end{equation}
and the  remainder term $R$ is controlled in the $H^1(\GG)$ norm, 
orbital stability of Definition \ref{def-stab} follows with $\varepsilon = C \delta$ for some $C > 0$. 

The main issue for getting the coercivity bounds (\ref{coercivity}) is that $\mathcal L_-$ and $\mathcal L_+$ may admit zero and negative eigenvalues. 

The stationary NLS equation (\ref{stationary}) can be written as 
$\mathcal{L}_- \Phi_{\omega} = 0$ with $\Phi_{\omega} \in \mathcal{D}(\Delta_{\GG})$. Hence, $\Phi_\omega\in \Ker \mathcal L_-$. 
If $\Phi_{\omega}(x) > 0$ for every $x \in \GG$, which is true for the ground state due to the property (3) of Proposition \ref{prop-EL}, then the spectrum of $\mathcal{L}_-$ for the Neumann--Kirchhoff and $\delta$-type conditions is non-negative and the kernel is non-degenerate with $\Ker(\mathcal{L}_-) = {\rm span}(\Phi_{\omega})$  \cite{[ACFN14],CFN17}. The same properties of $\mathcal{L}_-$ were recently proven in \cite{Kur19} for more general boundary conditions on a metric graph $\GG$.

On the other hand, $\mathcal{L}_+$ cannot be non-negative because 
\begin{equation}
\label{negativity}
\langle \mathcal{L}_+ \Phi_{\omega}, \Phi_{\omega} \rangle = - 2p(p+1) \| \Phi_{\omega} \|^{2p+2}_{L^{2p+2}(\GG)} < 0.
\end{equation}
Therefore, $\mathcal{L}_+$ has at least one negative eigenvalue. The negative eigenvalue of $\mathcal{L}_+$ and the zero eigenvalue of $\mathcal{L}_-$ are potentially a source of instability if the constraint of fixed mass $M(\Psi) = \mu$ does not ensure strict positivity for $\mathcal L_-$ and $\mathcal L_+$. The fixed mass constraint induces the constrained space $L^2_c \subset L^2(\GG)$ defined by 
$$
L^2_c := \{u\in L^2(\GG) : \ \ \langle u,\Phi_\omega\rangle = 0\}.
$$ 
If ${\rm Ker}(\mathcal{L}_-) = {\rm span}(\Phi_{\omega})$ and the spectrum of $\mathcal{L}_-$ is non-negative, then $\mathcal{L}_-$ is strictly positive in $L^2_c$. Controlling positivity of $\mathcal{L}_+$ in $L^2_c$ gives the orbital stability or instability of standing waves in the time evolution 
of the NLS equation (\ref{tNLS}) according to the following result \cite{Grillakis,GSS,VK}.

\begin{proposition}
	Assume that $\Phi_{\omega} \in \mathcal{D}(\Delta_{\GG})$ is a positive real solution of the stationary NLS equation (\ref{stationary}) and $\omega$ does not belong to the spectrum of $-\Delta_{\GG} + V(x)$ in $L^2(\GG)$. Then,
\begin{itemize}
	\item[(i)] If ${\rm Ker}(\mathcal{L}_-) = {\rm span}(\Phi_{\omega})$ and $\mathcal{L}_+$ has exactly one negative eigenvalue with the rest of its spectrum in $L^2(\GG)$ being strictly positive, then $\Phi_{\omega}$ is orbitally stable in $H^1_C(\GG)$ if the following slope condition holds:
	\begin{equation}\label{VK}
	\frac{d}{d\omega}\|\Phi_\omega\|^2<0
	\end{equation}
	while it is orbitally unstable if $\frac{d}{d\omega}\|\Phi_\omega\|^2>0$. 	
	\item[(ii)] If $\mathcal{L}_+$ has two or more negative eigenvalues, then $\Phi_{\omega}$ is orbitally unstable in $H^1_C(\GG)$ 		
\end{itemize}
\label{prop-stability}
\end{proposition}

\begin{remark} 
	$\Phi_\omega$ is a negative vector for $\mathcal L_+$ due to (\ref{negativity}), but it is not the eigenvector of $\mathcal{L}_+$. The condition \eqref{VK} assures strict positivity of $\mathcal L_+$ on $L^2_c$. Nothing can be said on stability of standing waves without further analysis if $\frac{d}{d\omega}\|\Phi_\omega\|^2 = 0$.
\end{remark}

\begin{remark}
	The number of negative eigenvalues of $\mathcal L_+$ is called the {\em Morse index} of the linearization. Proposition \ref{prop-stability} shows that the Morse index is a key object in the analysis of stability.
\end{remark}

\begin{remark}
	When the bound state $\Phi_{\omega}$ is a ground state of the 
	variational problem (\ref{minimization}), and more generally when the bound state is a local minimizer of a constrained variational problem with exactly one constraint, then the linearization can have at most one negative direction so that $\mathcal{L}_+$ has exactly one negative eigenvalue due to (\ref{negativity}). When the bound state $\Phi_{\omega}$ is a saddle point of a constrained variational problem, then the linearization may have more negative directions so that $\mathcal{L}_+$ may have more than one negative eigenvalues. 
\end{remark}

The first proof of orbital stability for the ground state of the NLS equation in $\mathbb R^n$ was given by Cazenave and Lions  \cite{CL82}. The proof does not use the expansion of the Lyapunov function (\ref{S-def}) but relies on properties of minimizing sequences and the concentration-compactness theorem. The proof is valid if the ground state in the minimization problem 
(\ref{minimization}) is unique. If the ground state is not unique, one can only obtain the stability of the set of the ground states, which is a weaker result.  Uniqueness of the ground states for the NLS equation on metric graphs was recently discussed in \cite{DST20} and can be proven for some simplest graphs such as the tadpole graph or the graph with a terminal edge, see Figure \ref{fig-graphs} below. Uniqueness of the ground state on the metric graphs with the external potential $V$ can also be proven 
for the branch of ground states bifurcating from the lowest eigenvalue of the associated self-adjoint operator \cite{CFN17}, as in Proposition \ref{prop-bifurcation} below.  

Besides the variational problem (\ref{minimization}), other variational characterizations of bound states have been used for the NLS equation on a metric graph $\GG$. One can consider the variational problem 
\begin{equation}
\label{infB}
\mathcal{B}_{\omega} = \inf_{\Psi \in H^1_C(\GG)} \left\{ B_{\omega}(\Psi) : \quad \| \Psi \|_{L^{2p+2}(\GG)} = 1 \right\}, 
\end{equation}
where 
\begin{equation}
\label{B-def}
B_{\omega}(\Psi) := \| \nabla \Psi \|_{L^2(\GG)}^2 + \int_{\GG} V(x) |\Psi|^2 - \omega \| \Psi \|^2_{L^2(\GG)}
\end{equation}
and $\omega < 0$. The variational problem (\ref{infB}) gives generally a larger set of standing waves with profiles $\Phi_{\omega}$ compared to the set of ground states in the variational problem (\ref{minimization}). 

Versions of the variational problem (\ref{infB}) arise in the determination of
the best constant of the Sobolev inequality, which is equivalent to the Gagliardo--Nirenberg inequality
in $\RE^n$ (see, e.g., \cite{Agueh06, Agueh08, DELL14}).
However, the minimizer of the variational problem (\ref{infB}) on a metric graph $\GG$ does not give the best constant in the Gagliardo--Nirenberg inequality \cite{NP20}.

Another variational problem is the minimization of the action functional 
on the Nehari manifold, 
\begin{equation}
\label{Nehari}
\mathcal{N}_{\omega} := \inf_{\Psi \in H^1_C(\GG) \backslash \{0\}} \left\{ S_{\omega}(\Psi) : \quad B_{\omega}(\Psi) = (p+1) \| \Psi \|^{2p+2}_{L^{2p+2}(\GG) } \right\},
\end{equation}
where $S_{\omega}$ is given by (\ref{S-def}) and $B_{\omega}$ is given by (\ref{B-def}). This approach was used in \cite{FOO2008} for the $\delta$-type potential on the line $\mathbb{R}$ and generalized in \cite{[ACFN14]} for the star graph with the $\delta$-type conditions. More recently, the variational problems on the Nehari manifold were analyzed in \cite{Pankov19,Pankov18}.

Relation between the variational problems (\ref{infB}) and \eqref{Nehari} 
was studied in \cite{NP20} for $p = 2$ and for the tadpole graph. 
It was shown that although the minimizers of the two problems do not generally coincide, if the minimizers	satisfy the same monotonicity properties on the tadpole graph (see Proposition \ref{global-existence} below), then the minimizer of one variational problem is at least a local minimizer of the other variational problem.

\subsection{Example of star graphs} 

The NLS equation on star graphs with $N$ unbounded edges was the first example studied in the context of the existence of ground states. In \cite{[ACFN14], acfn-aihp} (see also \cite{N14} for an overview) bound states were classified in the case of attractive and repulsive $\delta$-type conditions at the vertex with the parameter $\alpha$ in (\ref{BC-delta}). In this simple case, the bound states are obtained by gluing pieces of solitons of the line, and the classification is pictorially described by using the number of tails and bumps of the resulting bound state. 

The standing wave $\Phi^j_\omega$ is indexed by the number $j$ of bumps and it turns out that there are less bumps than tails for $\alpha<0$ and more bumps than tails for $\alpha>0$. The self-adjoint operators $\mathcal L_\pm$ have been studied around the $N$-tail state (no bumps) in \cite{[ACFN14], acfn-aihp}, where $\Phi^0_{\omega}$ is the ground state for small masses in the minimization problem (\ref{minimization}) and also the the minimizer of the action on the Nehari manifold (\ref{Nehari}) for $\alpha < 0$. Later, spectral properties of the self-adjoint operators $\mathcal L_\pm$ for the bound states $\Phi^j_\omega$ with $j \geq 1$ (at least one bump), and also for $\alpha > 0$ have been studied with different techniques in \cite{AngGol2} and more comprehensively in \cite{K19}. 

We denote by $n(\mathcal L_\pm)$  the Morse index (number of negative eigenvalues counted with their multiplicity) and by $z(\mathcal L_\pm):=\text{dim Ker}\ {\mathcal L_\pm}$. Let $\Phi^j_\omega$ be the bound state for $p>0$, 
$$
j=0,\cdots, \frac{N-1}{2}, \quad \omega<-\frac{\alpha^2}{(N-2j)^2}.
$$
It was found in \cite{K19} that 
\begin{itemize} 
	\item The absolutely continuous spectrum of $\mathcal L_\pm$  is located on $[|\omega|,+\infty)$,
	\item The discrete spectrum of $L_\pm$ in $(-\infty, |\omega|)$ is finite,
	\item $z(\mathcal L_-)=1, \ \ \ \ \ \ n(\mathcal L_-)=0$
	\item $z(\mathcal L_+)=0, \ \ \ \ \ \ n(\mathcal L_+)=
	\begin{cases}
	j+1 \ \ \ \ \ \text{if} \ \ \alpha<0 \\ N-j\ \ \ \ \text{if} \ \ \alpha>0
	\end{cases}	$
	\item $\frac{d}{d\omega}||\Phi_\omega^0||^2_2<0 \ \ \ \ \omega<-\frac{\alpha^2}{N^2}\, \ \text{and} \ \alpha<0  $
\end{itemize}
Proposition \ref{prop-stability} applies so that $\Phi_\omega^0$ is orbitally stable for $\omega<-\frac{\alpha^2}{N^2}$ and $\alpha < 0$, and $\Phi_\omega^j$ is orbitally unstable for any $j \geq 1$ (at least one bump) if $\alpha < 0$ and for any $j \geq 0$ if $\alpha > 0$.

The Neumann--Kirchhoff conditions ($\alpha=0$) support the bound states made by gluing half-solitons on the half-lines of the star graph for odd $N$ (obtaining a single degenerate bound state invariant by permutation of edges), or by placing shifted solitons on half of the lines of the star graph for $N$ even (see again \cite{[ACFN14], N14}). We will show in Section \ref{sec-4} that there is no ground state in this case. However, the bound states could in principle be linearly or orbitally stable. It was proved in \cite{KP1,KP2} that the bound states are unstable both for odd and even $N$. 

Other classes of boundary conditions have been considered in the last years, allowing for discontinuous wave-functions at vertices, as in the boundary conditions (\ref{BC-delta-prime}) and (\ref{NK-generalized}).
Stability and instability of standing waves on star graphs in the case of a $\delta'$ vertex have been considered in \cite{PG20,G22}. The case of generalized Kirchhoff boundary conditions has been studied in \cite{KP2,KGP} when the strength of nonlinearity is suitably adapted on any edge to allow for the existence of shifted states. It is shown that these shifted translating states are orbitally unstable.

\section{Variational methods for the ground state}
\label{sec-4}

In order to construct the ground state of the variational problem (\ref{minimization}), we first note that the energy functional 
$E(\Psi)$ is not bounded from below, irrespectively of the structure of the metric graph $\GG$. To understand this, we can consider Neumann--Kirchhoff conditions and no potential $V(x)$. Let us fix $\Psi \in H^1_C (\GG)$ and compute $E(\lambda \Psi)$ at the scaled function $\lambda \Psi$ with parameter $\lambda > 0$. 
For every $p > 0$ we have
$$ 
E(\lambda \Psi) =  {\lambda^2}  \| \Psi' \|^2_{L^2 (\GG)} -  \lambda^{2p+2}
\| \Psi \|^{2p+2}_{L^{2p+2} (\GG)} \ \to \ - \infty \quad \quad \text{as } \ \lambda \to + \infty.
$$
This implies that the energy is unbounded from below 
and it does not possess an infimum without a constraint. 
The situation is however different if the mass constraint 
$\| \Psi \|_{L^2 (\GG)}^2 = \mu$ is imposed in the variational problem (\ref{minimization}).

Recall the Gagliardo-Nirenberg inequality (\ref{gagliardo2}) for the graphs 
with finitely many edges and at least one half-line. After taking into account that the mass is fixed, we obtain 
$$
\| \Psi \|^{2p+2}_{L^{2p+2}(\GG)} \leq \ C_p \mu^{\frac{p+2}{2}} \| \Psi' \|^p_{L^2(\GG)}. 
$$
For every $p > 0$ we have the following lower bound on the energy
\begin{equation}
\label{lower-bound}
E(\Psi) \ \geq \  \| \Psi' \|^2_{L^2(\GG)} - C_p \mu^{\f{p+2}{2}} \| \Psi' \|^p_{L^2(\GG)},
\end{equation}
which allows us to conclude the following.
\begin{itemize}
\item If $0 < p < 2$, then the inequality (\ref{lower-bound}) implies that $E(\Psi)$ is bounded from below. The infimum always exists and can be approached by a minimizing sequence. However, the minimizing sequence may not converge to a minimizer $\Phi \in H^1_{C}(\GG)$ if it is not attained.

\item If $p = 2$ \emph{(critical power)},
then the inequality (\ref{lower-bound}) implies that 
$$
E(\Psi) \ \geq \ \left( 1 - C_2 \mu^{2} \right)  \| \Psi' \|^2_{L^2(\GG)},
$$
so that there exists a \emph{critical mass} $\mu_\GG := C_2^{-1/2}$ below which all states have \emph{positive} energy $E(\Psi)$. The infimum exists for $\mu \in (0,\mu_{\GG}]$.

\item If $p > 2$, then the inequality (\ref{lower-bound}) is not conclusive. 
In this case, one can use the following scaling transformation 
for the half-line: $\Psi(x) \mapsto \sqrt{\lambda} \Psi(\lambda x) =: \Psi_{\lambda}(x)$, $x > 0$. This transformation preserves the mass $M(\Psi_{\lambda}) = M(\Psi) = \mu$ and changes the energy:
$$
E(\Psi_{\lambda}) = \lambda^2 \| \Psi' \|^2_{L^2} - \lambda^p \| \Psi \|^{2p+2}_{L^{2p+2}} \to -\infty \quad \mbox{\rm as} \quad  \lambda \to +\infty.
$$
Hence, the infimum does not exist for any $\mu > 0$.
\end{itemize}

A similar analysis and the same conclusions hold true also for compact graphs, making use of the corresponding Gagliardo-Nirenberg inequality (\ref{gagliardo3}). 
Due to the conclusions above, we will only consider the bound states in the subcritical and critical cases for $p \in (0,2]$.

\begin{figure}[h!]
	\begin{center}
		\begin{tikzpicture}[xscale = 2.0, yscale=2.0]			
		\draw (-2.5,0) -- (2.5,0);
		\draw [black, thick] plot [domain=-2.5:2.5, smooth] (\x, {2.5/(exp(-2*\x)+ exp(2*\x))});
		\fill (-2.5,0) node[left] {$\infty$};
		\fill (2.5,0) node[right] {$\infty$};
		\fill (0.8, 1) node[right] {$\phi_{\mu}$};
		\end{tikzpicture}
	\end{center}
	\caption{NLS soliton (\ref{soliton-mu}) on the line $\mathbb{R}$}
	\label{fig-soliton}
\end{figure}
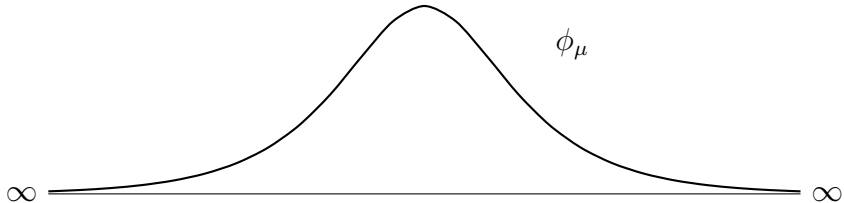

For $p \in (0,2)$, the ground state of the NLS equation are well-known in the cases of line $\mathbb{R}$ and half-line $\mathbb{R}_+$. The unique solution of the stationary NLS equation 
\begin{equation}
\label{soliton-NLS}
-\phi_{\omega}'' - (p+1) \phi_{\omega}^{2p+1} = \omega \phi_{\omega}.
\end{equation}
is given up to a translation by the NLS soliton
\begin{equation}
\label{soliton}
\phi_\omega(x) = |\omega|^{\frac{1}{2p}}{\rm sech}^{\frac{1}{p}}(p\sqrt{|\omega|} x), \quad x \in \mathbb{R}.
\end{equation}
The mass $\mu$ of the NLS soliton can be calculated as a function of frequency $\omega$:
\begin{equation}
\label{soliton-mass}
\mu = C_p |\omega|^{\frac{1}{p}-\frac{1}{2}}, \quad 
C_p := \int_{\mathbb{R}} {\rm sech}^{\frac{2}{p}}(px) dx.
\end{equation}
Inverting this relation for $p \in (0,2)$ allows us to parametrize the NLS soliton (\ref{soliton}) by $\mu$:
\begin{equation}
\label{soliton-mu}
\phi_\mu (x) \ = \ c_p \mu^{\f 1 {2-p}} \hbox{\mbox{sech}}^{\f 1 {p}} (c_p^p
\mu^{\f {p} {2-p}} x).
\end{equation}
where $c_p$ is a positive constant that depends on $p \in (0,2)$ but not on $\mu$. This solution is shown on Figure \ref{fig-soliton}.

The soliton $\phi_{\mu}$ up to a translation in $x$ is \emph{the ground state} of the variational problem (\ref{minimization}) for $\GG= \mathbb{R}$ and every $\mu > 0$ if $p \in (0,2)$. Similarly, the half-soliton for the double mass $\phi_{2\mu}(x)$ for $x > 0$ 
is \emph{the ground state} of the variational problem (\ref{minimization}) for $\GG= \mathbb{R}_+$ and every $\mu > 0$ if $p \in (0,2)$. Notice that the translational symmetry is broken by the Neumann condition at $x = 0$ so that the ground state on $\mathbb{R}_+$ is uniquely defined.

In the critical case $p=2$, the mass $\mu$ of the NLS soliton (\ref{soliton-mass}) is independent on $\omega$. The constant values of the mass on the half-line $\mathbb{R}_+$ and the full line $\mathbb{R}$ are given by 
\begin{equation}
\label{mass-critical}
\mu_{\mathbb{R}_+} =  \frac{\pi}{4}, \quad 
\mu_{\mathbb{R}}  =\frac{\pi}{2}.
\end{equation}
Consequently, if $p = 2$, then the ground state of the variational problem (\ref{minimization}) on $\GG = \mathbb{R}$ and $\GG= \mathbb{R}_+$ exists only for 
$\mu = \mu_{\mathbb{R}}$ and $\mu = \mu_{\mathbb{R}_+}$ respectively.

\subsection{Compact graphs} 

Ground states on compact graphs have been treated in  \cite{CDS18,Dovetta18}. The subcritical case is simple because of the compact embedding of $H^1_C(\GG)$ into $L^p(\GG)$ for every $p \geq 2$. As a result, a minimizing sequence $\{ \Phi_n \}_{n \in \mathbb{N}}$ satisfying the mass constraint $M(\Phi_n) = \mu$ converges strongly in  $H^1_C(\GG)$ and its limit $\Phi \in H^1_C(\GG)$ also satisfies the mass constraint $M(\Phi) = \mu$.
Hence $\Phi$ is the ground state of the variational problem (\ref{minimization}) 
and it is attained irrespectively of details of the compact graph $\GG$. 

It has been proved in \cite{CDS18} that for $p \in (0,2]$ there exists the critical mass $\mu_*$ such that the unique constant solution of the stationary NLS equation (\ref{stationary}) with $V(x) = 0$ is the ground state for $\mu \in (0,\mu_*)$,  whereas the ground state is not a constant solution for $\mu > \mu_*$. 

A detailed analysis of the ground state and other bound states on the special example of a compact graph given by {\em the dumbbell graph} has been given in \cite{marzuola} and \cite{G19}.

\subsection{Unbounded graphs} 

We now consider unbounded metric graphs in the case of Neumann--Kirchhoff conditions and no external potential $V(x)$. The lack of compactness of $\GG$ makes the convergence of a minimizing sequence a non trivial issue, especially in the critical case $p = 2$. 

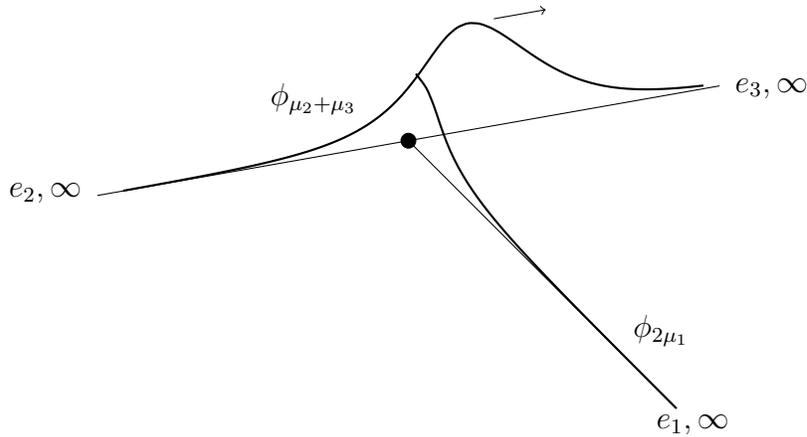
\begin{figure}[htb]
\begin{center}
		\begin{tikzpicture}[xscale= 1.4,yscale=1.4]			
			\draw [rotate around = {10:(0,0)}] (-3,0) -- (3,0);
			\draw  (0,0) -- (2.5,-2.5);
			\draw  [rotate around = {10:(0,0)}] [black, thick] plot [domain= -0.55:2.1, smooth] (\x +0.75, { 2/(exp(-2*\x)+ exp(2*\x))});
			\draw  [rotate around = {10:(0,0)}] [black, thick] plot [domain=-3.5:-0.55, smooth] (\x +0.75, { 2/(exp(-2*\x)+ exp(2*\x))});
			\draw [rotate around = {315:(0,0)}] [black, thick] plot [domain=0:4, smooth] (\x -0.4, {1/(exp(-2*\x)+ exp(2*\x))});
			\draw plot [mark=*] coordinates {(0,0)};
			\draw [rotate around = {10:(0,0)}] [->] (1,1) -- (1.5,1);
			\fill (-3,-0.5) node[left] {$e_2,\infty$};
			\fill (3, 0.5) node[right] {$e_3,\infty$};
			\fill (2.7, -2.7) node[] {$e_1, \infty$};
			\fill (-0.9, 0.4) node[] {$\phi_{\mu_2+\mu_3}$};
			\fill (2.4, -1.8) node[] {$\phi_{2\mu_1}$};
		\end{tikzpicture}
	\end{center}
		\caption{Vanishing of the minimizing sequence on the star graph with three edges}
\label{fig-star-three}
\end{figure}

What can go wrong in the convergence of a minimizing sequence can be understood from the following example (originally discussed in \cite{[ACFN2]}).
Let $\GG$ be a star graph with three infinite edges. For any function $\Psi \in H^1_C(\GG)$ with given mass $\mu$, it is always possible to find a function $\Phi \in H^1_C(\GG)$ with the same mass but lower energy shown on Figure \ref{fig-star-three}.

The original function $\Psi$ has masses $\mu_1$, $\mu_2$, and $\mu_3$ in its components on $e_1$, $e_2$, and $e_3$, respectively. Put half of the soliton $\phi_{2\mu_1}$ on $e_1$ (recall that $\phi_{2\mu_1}$ has the same mass $\mu_1$) and a soliton of mass $\mu_2+\mu_3$ on the union of the two half-lines $e_2 \cup e_3$. The Neumann--Kirchhoff boundary conditions are satisfied. Since $\phi_{\mu_2+\mu_3}$ is the ground state on the line for the mass $\mu_2+\mu_3$, and $\phi_{2\mu_1}$ is the ground state on the half-line for a mass $\mu_1$, we have $E(\Psi) \geq  E(\Phi)$. However, $\Phi$ is not a ground state. An easy calculation shows that for $p\in (0,2)$, $\inf{E}(\Phi)$ is attained for $\mu_1 = 0$. As $\mu_1$ decreases, the soliton $\phi_{\mu_2+\mu_3}$ on $e_2 \cup e_3$ moves away from the vertex of the star graph $\GG$ and runs away to $\infty$ taking all the mass while the half-soliton on $e_1$ disappears. However, the limit of this $\Phi$ is the zero function in $H^1_C(\GG)$.

This behavior of the minimizing sequences is typical on unbounded star graphs: they move all the mass along the half lines to infinity and vanish near the vertex points. The bound states given by the half solitons on the star graphs are saddle points of the constrained minimization problem (\ref{minimization}) \cite{[ACFN2]}. They are shown to be unstable in the time evolution of the 
NLS equation for $p \in (0,2)$ \cite{KP1}. Shifted states 
can be constructed for even number of half lines or for generalized Neumann--Kirchhoff conditions and most of them are linearly unstable \cite{KP2}. Even if the shifted states are linearly stable, they are nonlinearly unstable \cite{KGP}.

A thorough analysis of the ground state on unbounded graphs $\GG$ with finitely many edges and the Neumann--Kirchhoff conditions has been given in a series of papers by R. Adami, E. Serra and P. Tilli \cite{AST15, AST16, AST17,AST19}. The following proposition gives the topological and metric conditions on $\GG$ that guarantee non-existence of the ground state. 

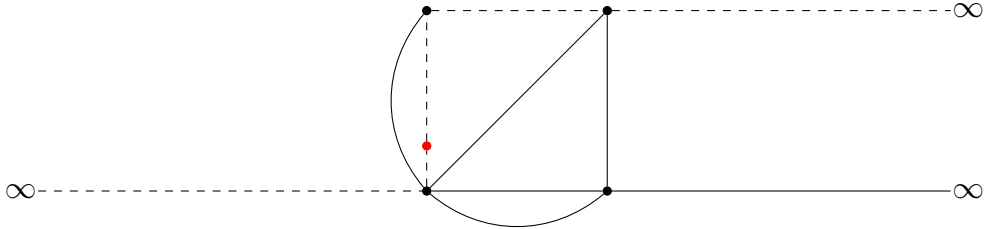
\begin{figure}[htb!]
\begin{center}
\begin{tikzpicture}[xscale= 1.2,yscale=1.2]
\node at (-6.5,0) [infinito](0) {${\infty}$};
\node at (-2,0) [nodo,black] (5) {};
\node at (0,0) [nodo,black] (4) {};
\node at (-2,0.5) [nodo,red] (9) {};
\node at (-2,2) [nodo,black] (2) {};
\node at (0,2) [nodo,black] (3) {};
\node at (4,0) [infinito] (8) {${\infty}$};
\node at (4,2) [infinito] (7) {${\infty}$};


\draw [-,black,dashed] (0) -- (5);
\draw [-,black] (5) -- (4);
\draw [-,black] (5) to [out=-40,in=-140] (4);
\draw [-,black] (5) to [out=130,in=-130] (2);
\draw [-,black] (5) -- (3);
\draw [-,black,dashed] (5) -- (9);
\draw [-,black,dashed] (9) -- (2);
\draw [-,black,dashed] (2) -- (3);
\draw [-,black,dashed] (3) -- (7);
\draw [-,black] (4) -- (3);
\draw [-,black] (4) -- (8);
\end{tikzpicture}
\end{center}
\caption{The dashed line represents the trail from the red point to infinity along two different paths.}
\label{fig-path}
\end{figure}

\begin{figure}[htb!]
\begin{center}
\begin{tikzpicture}[xscale= 2.0,yscale=1.8]
%

\node at (-2.5,0) [infinito](0) {$\infty$};
\node at (0.5,0) [nodo,black] (5) {};
\node at (0.5,0.6) [nodo,black] (1) {};
\node at (3.5,0) [infinito] (8) {$\infty$};

\draw [-,black] (0) -- (5);
\draw [-,black] (5) -- (8);
\draw [-,black] (0.5,0.3) circle (0.3) ;
\draw [-,black] (0.5,0.9) circle (0.3) ;

\end{tikzpicture}
\end{center}
\caption{A tower graph with two bubbles.}
\label{fig-bubble}
\end{figure}
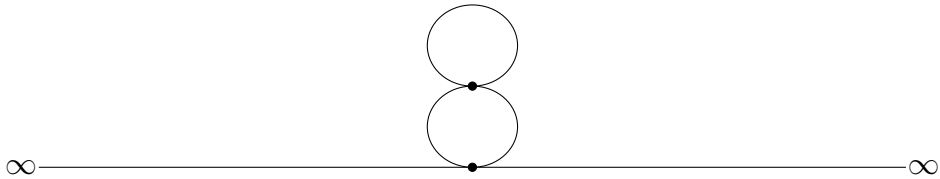

\begin{proposition}
\label{AST-condition-H-proposition}
Let $0<p<2$. Assume that every point of the graph $\GG$ lies on a trail containing two different half-lines, that is, one can reach infinity from every point of the graph through two disjoint paths as shown on Figure \ref{fig-path}. Then, the ground state does not exist unless $\GG$ is isometric to a bubble tower graph shown on Figure \ref{fig-bubble}.
\end{proposition}

\begin{remark} One can easily understand why ground states can exist on the bubble tower graphs. Fold a line on any bubble tower graph and then place the NLS soliton defined on the line. Neumann--Kirchhoff conditions are respected by placing the maximum of the soliton on top of the bubble tower and exploiting symmetry of the soliton. The folded NLS soliton is the ground state. It is unique because the translation invariance of the line has been broken by the vertices.
\end{remark}

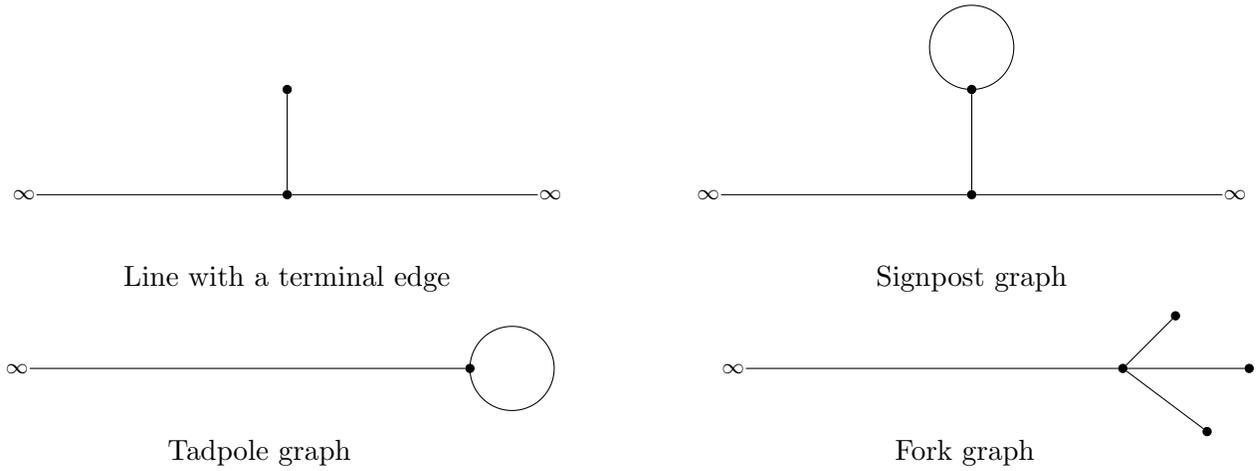
\begin{figure}[htb!]  
	\begin{center}
		\begin{tikzpicture}[scale= 1.4]
		\node at (-3.5,0) [infinito]  (1) {$\scriptstyle\infty$};
		\node at (-1,0) [nodo] (2) {};
		\node at (1.5,0) [infinito]  (3) {$\scriptstyle\infty$};
		\node at (-1,1) [nodo] (4) {};
		\draw [-] (1) -- (2) ;
		\draw [-] (2) -- (3) ;
		\draw [-] (2) -- (4) ;
		\node at (-1,-.8)  [minimum size=0pt] (10) {\small{Line with a {terminal edge}}};
		
		\node at (3,0) [infinito]  (5) {$\scriptstyle\infty$}; 
		\node at (5.5,0) [nodo] (6) {};
		\node at (8,0) [infinito]  (7) {$\scriptstyle\infty$};
		\node at (5.5,1) [nodo] (8) {};
		\draw [-] (5) -- (6) ;
		\draw [-] (6) -- (7) ;
		\draw [-] (6) -- (8) ;
		\draw(5.5,1.4) circle (0.4);
		\node at (5.5,-.8)  [minimum size=0pt] (11) {\small{Signpost graph}};
		\end{tikzpicture}
	\end{center} 
	\begin{center}
		\begin{tikzpicture}[scale= 1.4]
		\node at (-2.5,0) [infinito]  (1) {$\scriptstyle\infty$};
		\node at (1.8,0) [nodo] (2) {};
		\draw [-] (1) -- (2) ;
		\draw(2.2,0) circle (0.4);
		\node at (-0.2,-.8)  [minimum size=0pt] (10) {\small{Tadpole graph}};
		
		\node at (4.3,0) [infinito]  (5) {$\scriptstyle\infty$}; 
		\node at (8,0) [nodo] (6) {};
		\node at (8.5,.5) [nodo]  (7) {};
		\node at (8.8,-.6) [nodo] (8) {};
		\node at (9.2,0) [nodo] (9) {};
		\draw [-] (5) -- (6) ;
		\draw [-] (6) -- (7) ;
		\draw [-] (6) -- (8) ;
		\draw [-] (6) -- (9) ;
		\node at (6.5,-.8)  [minimum size=0pt] (11) {\small{Fork graph}};
		\end{tikzpicture}
	\end{center}
	\caption{Exampes of graphs where the ground state exists.}
	\label{fig-graphs}
\end{figure}

There are many metric graphs that escape the non-existence result of Proposition \ref{AST-condition-H-proposition}. The list includes graphs with the terminal (pendant) edge, the signpost graph, the tadpole graph, and the fork graph shown on Figure \ref{fig-graphs}. In fact, any graph with just one half-line violates the condition of Proposition \ref{AST-condition-H-proposition} 
and the ground state may exist on such graphs with just one half-line. 
On the other hand, simple metric graphs such as the double bridge graph (a circle with two half-lines) meet the hypotheses of Proposition \ref{AST-condition-H-proposition} and fail to trap a ground state.

The main ideas involved in the proof of Proposition \ref{AST-condition-H-proposition} are given as follows.

\begin{itemize}
	\item {\em Rearrangements.} This is an extension of a classical technique that, given any function in $H^1_C(\GG)$ allows to build another function which has the same $L^p(\GG)$ norms but smaller kinetic energy. There are two types of rearrangements useful in the analysis of ground states on metric graphs: 
	\begin{itemize}
		\item Decreasing rearrangements which associate to any function $\Psi\in H^1_C(\GG)$ a second function $\Psi^*$ on $H^1(\mathbb{R}_+)$ with the same $L^p(\mathbb{R}_+)$ norm and a lower kinetic energy;
		\item Symmetric rearrangements which do the same but for a function $\hat \Psi\in H^1(\RE)$.
	\end{itemize}

	\item {\em Comparison with the half-soliton and the full-soliton.} The following inequality holds for fixed mass $M(\Psi)=\mu$:
	\begin{equation}
	\label{inequalities-energy}
	E(\phi_{2\mu})=\inf_{\Psi\in H^1(\mathbb{R}_+)}E(\Psi) \leq \inf_{\Psi\in H^1_C(\GG)}E(\Psi) \leq \inf_{\Psi\in H^1(\RE)} E(\Psi) = E(\phi_{\mu}).
	\end{equation}
	The first part of this inequality is a consequence of monotonic rearrangement and states that no function on the metric graph can do better than the half-soliton on the half-line. The second part is a nontrivial consequence of a property of minimizing sequences.
	
	\item {\em Sufficient condition for existence of the ground state.} If there exists $\Phi \in H^1_C(\GG)$ such that 
	\begin{equation}
	\label{smaller-energy}
	E(\Phi) \leq \inf_{\Psi\in H^1(\RE)} E(\Psi) = E(\phi_\mu)
	\end{equation} 
	then the ground state exists. This is exactly the condition that prevents the minimizing sequence to runaway to infinity and to vanish inside $\GG$.
\end{itemize}

The ground state on the tadpole graph in the critical case $p = 2$ was studied in \cite{NP20} by using the variational problem (\ref{infB}). 
	It was shown that the ground state exists for $\mu \in (\mu_{\mathbb{R}_+},\mu_{\mathbb{R}}]$ and two bound states 
	coexist for some $\mu > \mu_{\mathbb{R}}$: one is a local minimizer 
	and the other one is a saddle point of the variational problem (\ref{minimization}). Local minimizers of the variational problem (\ref{minimization}) in the absence of the ground state are studied in \cite{PSV19}. 
	
Several papers study bound states (not necessarily ground states) making use of other analytical techniques. Number-theoretic properties for existence of standing waves on the double bridge graph were discovered and studied in \cite{NRS19}. Canonical perturbation theory is used to classify standing waves in \cite{GW1,GW2}. Existence of standing waves in the presence of generalized Neumann--Kirchhoff conditions is studied in  \cite{sabirov_ea13, sobirov_ea10}, with an attempt to establish a connection with integrable Hamiltonian systems. The integrable properties of the system of NLS equations on edges of the star graph were studied in \cite{caudrelier}.

\subsection{$\delta$-type conditions and external potentials}

The $\delta$-type conditions at one or more vertices or the presence of external potentials may give rise to a negative eigenvalue of
the linear operator associated to the quadratic part of the energy. In this case, a ground state always exist in the subcritical \cite{CFN17},
critical \cite{C18}, and supercritical \cite{Ardila18} cases. First examples of the ground state for the $\delta$-type conditions were considered in \cite{acfn-aihp, ACFN16}. We only sketch the subcritical case $p \in (0,2)$. 

Assume that $\GG$ is a connected graph with a finite number of edges and it is composed by at  least one half-line attached to a compact core. 
Let $V = V_+ - V_-$ with $V_\pm\geq 0$, $V_+\in L^1(\GG) + L^{\infty}(\GG)$, and $V_-\in L^r(\GG)$ for some $r\in [1,1+1/\mu]$.  The quadratic part of the energy is given by 
$$
E(\Psi)=\| \Psi ' \|^2_{L^2(\GG)} + \sum_{v\in V} \al(v) |\Psi(v)|^2 +\langle\Psi,V\Psi\rangle,
$$
which is generated by the self adjoint operator $H := -\Delta_{\GG} + V$ with the form domain $\mathcal{E}(H) = H^1_C(\GG)$. The spectrum of $H$ is defined in $L^2(\GG)$ with the domain $\mathcal{D}(H)$ given by the $\delta$-type conditions. Assume that the spectrum $\sigma(H)$ has the infimum at $-E_0 < 0$, which is an isolated eigenvalue of $H$.

\begin{proposition}
Let $0<p<2$. If the above assumptions hold true, then 
\[
-\infty <  \mathcal{E}_{\mu} \leq - E_0 \mu 
\]
for any $\mu>0$.  Moreover, there exists $\mu^\ast>0$ such that for $0< \mu<\mu^\ast$  the infimum $\mathcal{E}_{\mu}$ is attained, i.e., the ground state exists and it is orbitally stable.  
\label{prop-bifurcation}
\end{proposition}

\begin{remark}
If $-E_0$ is a \emph{simple eigenvalue} of $H$, one can use bifurcation theory to find a candidate for the ground state. It turns out that the ground state of Proposition \ref{prop-bifurcation} bifurcates from the linear eigenstate associated to the eigenvalue $-E_0$.
\end{remark}

\begin{remark} A sufficient condition to have an isolated eigenvalue of the linear operator $H$ is 
\begin{equation}
\int_\GG V  + \sum_{v\in V} \alpha(v) < 0
\end{equation}
which is a generalization of the condition on the line $\mathbb{R}$.
Simplicity of the ground state is a consequence of the resolvent (or semigroup) of $H$ being positivity improving (see \cite{CFN17}).
\end{remark}

\subsection{Localized nonlinearities}

A natural question is to understand what happens when the nonlinearity affects only a compact portion of the metric graph $\GG$. Namely, one can consider a metric graph where the dynamics is given by the NLS equation on a compact part $\mathcal K$ 
and by the linear Schr\"{o}dinger equation on the remaining unbounded part $\GG \backslash \mathcal{K}$. We call this setting as a metric graph with localized nonlinearity. Such models are physically meaningful, because it is reasonable to consider the wave interactions to be asymptotically linear. 

The simplest graph with the localized nonlinearity is given by the line, 
for which the nonlinearity is present on a bounded segment $[-L,L]$ with some $L > 0$ which is connected to two exterior half-lines $(-\infty,-L] \cup [L,\infty)$. When the segment length reduces to zero, that is, when $L \to 0$, the ruling equation is the linear Schr\"odinger equation on the line, with no bound states. In the opposite limiting case, when the segment length increases to infinity, that is, when $L \to \infty$, there is only a ground state (up to symmetries) coinciding with the soliton, and not any other bound states. Therefore, the intermediate case is rather non-trivial. 

The problem was addressed for some metric graphs in \cite{GSD11}, where 
the existence of many bound states in the presence of a rather complex pattern of scattering resonances was shown by means of numerical approximations. Analysis of ground states and more generally bound states for a metric graph with localized nonlinearity has been treated in the papers \cite{ST1,ST2, T16} in the subcritical case $0 < p < 2$ and in \cite{DT19} in the critical case $p=2$. 

A first unexpected result which makes localized nonlinearities quite different from the standard nonlinearities concerns existence of the ground state (see \cite{T16}). Let the energy functional $E : H^1_C(\GG)\to\RE$ be given by 
\[
E(u)=\frac12\|u'\|_{L^2(\GG)}^2-\frac{1}{2p+2}\|u\|_{L^{2p+2}(\K)}^{2p+2}
\]
and as before we consider minimization of $E(u)$ under the fixed constant 
mass $M(u) = \mu$ for $\mu > 0$. Then, 
$-\infty < \mathcal{E}_{\mu} \leq 0$ 
and the infimum is attained if it is strictly negative. Moreover, we have the following:
\begin{itemize}
	\item If $p\in(0,1)$, then there exists the ground state for every $\mu>0$ 
	\medskip
	\item If $p\in [1,2)$, then there exist two constants $0 < \mu_2 < \mu_1$ such that:\\
	--  the ground state exists for every $\mu>\mu_1$;\\
	-- the ground state does not exist for every $\mu<\mu_2$.	
\end{itemize}
The physically important cubic NLS equation $(p=1)$ is the critical case since there exists the gap in the set of masses admitting a ground state for $p > 1$ 
and no gap exists for $p < 1$. Existence of bound states with higher energy levels than the ground state energy has also been studied for localized nonlinearities  \cite{ST1, ST2}.

\section{Bound states via period function}

Bound states including the ground state, if it exists, can be analyzed by other methods, e.g., by studying orbits on the phase plane, since each component 
of $\Phi \in H^1_C(\GG)$ satisfies the second-order differential equation on each edge of the general graph $\GG$. Let us look again for real-valued solutions of the stationary NLS equation (\ref{stationary}) but rewritten without the potential $V(x)$ as
\begin{equation}
\label{stationary-no-potential}
- \Delta_{\GG} \Phi_{\omega} - (p+1) |\Phi_{\omega}|^{2p} \Phi_{\omega} = \omega \Phi_{\omega}.
\end{equation} 
The following class of functions define the {\em edge-localized states} 
on a general graph $\GG$:
\begin{itemize}
\item $\Phi_{\omega}(x) > 0$ for every $x \in \GG$;
\item on each bounded edge $e \in E$, there is at most one local critical point of $\phi_{\omega}(x)$ (either maximum or minimum) inside the edge;
\item on each unbounded edge $e \in E$, the function 
$\phi_{\omega}(x)$ is monotonically decreasing and has exponential decay to $0$ as $x \to +\infty$.
\end{itemize}

\begin{figure}[htbp] 
	\centering
	\includegraphics[width=1.8in, height = 1.5in]{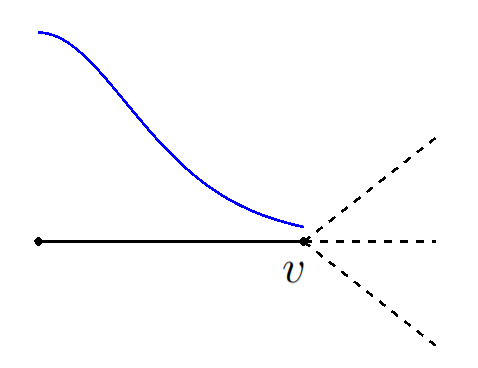} 
	\hspace{0.5cm}
	\includegraphics[width=1.8in, height = 1.5in]{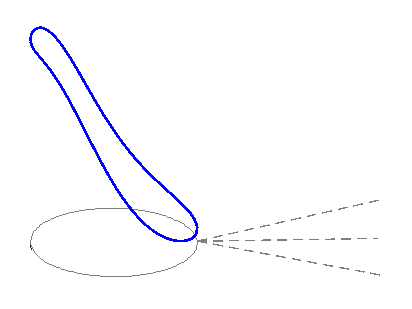}
	\hspace{0.5cm}
	\includegraphics[width=1.8in, height = 1.5in]{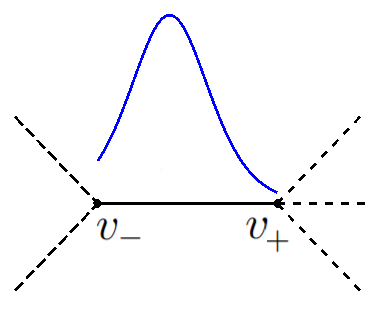}
	\caption{Schematic illustration of a half-pulse on a pendant (left) 
		and a pulse on a looping edge (middle) and on an internal edge (right). }
	\label{fig-states}
\end{figure}

On each bounded edge $e$ of the graph $\GG$ we point out the case when $\phi_{\omega}$ has a local maximum inside the edge, and call such a profile {\em a pulse}. If the maximum happens at one vertex of the edge $e \in \GG$, 
then we call it {\em a half-pulse}. The examples of a pulse and a half-pulse on different types of edges are given in Figure \ref{fig-states}. 

Figure \ref{fig-states} does not give an exhaustive list of all possible edge-localized states as one can think of cases when, for example, a local maximum is strictly inside of the pendant edge, half-pulses are formed at a vertex connecting several edges, or two maxima are located at both ends of the edge $e$. The corresponding bound states correspond to larger values 
of energy at fixed mass and we will not consider such bound states here.

We shall distinguish edge-localized states based on the number of pulse profiles present. If a pulse profile occurs on a single edge of the graph, we call the corresponding bound state {\em a single-pulse state} 
and if several pulse profiles occur on multiple edges of the graph, we call it {\em a multi-pulse state}.

Existence of edge-localized states with pulse profiles on certain edges was confirmed analytically and numerically for several simple graph models, e.g. the tadpole graph \cite{AST16,cfn15,DST20,NP20,nps}, the dumbbell graph \cite{G19,marzuola}, and the flower graph \cite{KMPZ}. More general results were obtained in {\em the limit of large mass} in \cite{AST19,BMP} and \cite{KP21}, where the existence of single-pulse states and multi-pulse states was proven for a general graph, see Section \ref{sec-multi-pulse}. 

Here we focus on a novel analytical approach where pulses on bounded edges are obtained by using properties of the {\em period function} for second-order differential equations. The period function is typically used for analysis of existence of periodic solutions to nonlinear evolution equations \cite{Vill1, GV15} as well as their spectral stability \cite{GeyerSAPM,GP17}. Recently, the period function was used for parts of the periodic solutions in the context of metric graphs in \cite{KMPZ,NP20} and for solitary waves with compact singular heads in \cite{PKR}.

\subsection{Integral curves for the pulse solutions}
\label{subsection-dynamical-system}

A pulse profile might appear generally on any bounded edge $e$ of a graph $\GG$. We parametrize the edge $e$ by an interval $[-\ell_1, \ell_2]$. The pulse profile on $e$ satisfies the stationary NLS equation \eqref{stationary-no-potential} on the interval, 
\begin{equation}
\label{pulse-profile-ode}
-\phi''(x) - (p+1) \phi(x)^{2p+1} = \omega \phi(x), \quad 
x \in (-\ell_1, \ell_2).
\end{equation}
The parametrization of $e$ is chosen in such a way that $x=0$ corresponds to the maximum of the pulse, this imposes the condition $\phi'(0) = 0$.  If $e$ is a looping edge, then $\ell_1 = \ell_2$ follows by the symmetry of the Neumann--Kirchhoff conditions. If $e$ is a pendant (terminal) edge, then we set $\ell_1 = 0$ and consider the half-pulse with the maximum at $x = 0$. 
If $e$ is an internal edge, then it is typical that $\ell_1 \neq \ell_2$ .
The exact representation of the pulse on the edge $e$ is determined by information from the remainder of the graph $\GG \backslash \{e\}$.

The parameter $\omega < 0$ in (\ref{pulse-profile-ode}) can be eliminated by 
the scaling transformation 
\begin{equation}
\label{scaling-transform}
\omega : = -\epsilon^2, \quad \phi(x) = \epsilon^{\frac{1}{p}} u(\epsilon x),
\end{equation}
where $\epsilon > 0$. The stationary NLS equation (\ref{pulse-profile-ode}) transforms into  the normalized equation 
\begin{equation}
\label{stat-NLS-halfedge} 
-u''(z) + u(z) - (p+1) u(z)^{2p+1} = 0, \quad z \in (-\epsilon \ell_1, \epsilon \ell_2),
\end{equation}
where $z = \epsilon x$ and $\epsilon$ determines the new length of the interval.

The scaling transformation (\ref{scaling-transform}) does not change the derivative condition at the origin where the pulse maximum is located, hence $u'(0) = 0$ is still true. This implies that the solution is an even function of $z$ and it makes sense to partition the interval $[-\epsilon \ell_1, \epsilon \ell_2]$ into two intervals $[-\epsilon \ell_1, 0]$ and $[0, \epsilon \ell_2]$. On each of these intervals the pulse solution is monotonic towards the vertices of the edge. The analysis on both intervals is very similar, so we further focus on the properties of solutions restricted to the interval 
$[0, \epsilon \ell]$ subject to the boundary conditions 
\begin{equation}
\label{endpoint-values-of-u}
z = 0 : \quad (u, u') = (\mathfrak{p}_+, 0), \quad \text{and} \quad 
z = \epsilon \ell : \quad (u,u') = (\mathfrak{p}, -\mathfrak{q}),
\end{equation}
where $\mathfrak{p}_+ > \mathfrak{p} > 0$ and $\mathfrak{q} > 0$ due to positivity and monotonicity assumptions on $u$ restricted to the interval $[0, \epsilon \ell]$. 

Let us denote $v := u'$. The solution $(u(z), v(z))$ stays for all $z \in (0, \epsilon \ell)$ on the invariant level curve given by 
\begin{equation}
\label{curve-invariant}
E_\beta := \{(u,v) : \quad v^2-u^2+u^{2p+2} = \beta \}, \quad 
\beta \geq \beta_p := -\frac{p}{(p+1)^{\frac{p+1}{p}}},
\end{equation}
where $\beta_p$ is the minimal value of the function $f(u) := -u^2 + u^{2p+2}$.

The choice $\beta = \beta_p$ gives the constant solution $u(z) = u_p := (p+1)^{-\frac{1}{2p}}$ for all $z \in [0, \epsilon \ell]$. Since we are looking only for pulse profiles on the edge $e$, the constant solution is neglected 
and $\beta > \beta_p$.

\begin{figure}[htbp!] 
	\centering
	\includegraphics[width=4in, height = 3in]{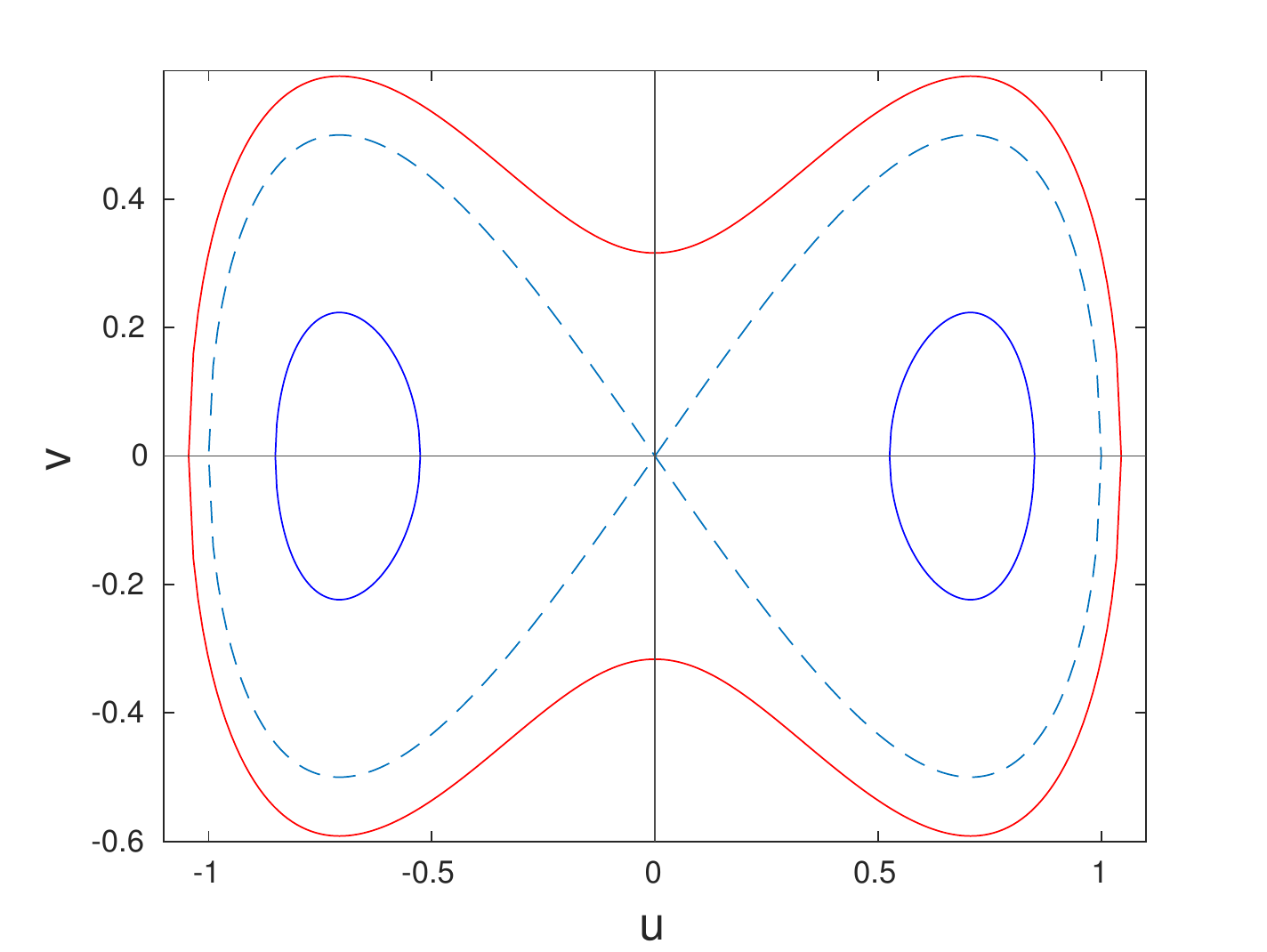}
	\caption{Phase portrait of the differential equation $-u'' + u - 2 u^3 = 0$ on the $(u,v)$-plane given by the level curves $E_{\beta}$.}
	\label{fig-phase-portrait}
\end{figure}

Figure \ref{fig-phase-portrait} shows the phase portrait on the $(u,v)$-plane given by the level curves $E_\beta$. Each level curve $E_{\beta}$ represents a solution to the differential equation \eqref{stat-NLS-halfedge} for $p = 1$. 

We note the following:
\begin{itemize}
\item $E_0$ at the level $\beta = 0$ represents two homoclinic orbits (the blue dashed lines),
and the solution on the positive homoclinic orbit is represented by the shifted NLS soliton 
\begin{equation}
\label{NLS-soliton}
u(z) = {\rm sech}^{\frac{1}{p}}p(z+a),
\end{equation}
where $a \in \mathbb{R}$ is an arbitrary translation parameter. 

\item $E_\beta$ with $\beta \in (\beta_p, 0)$ corresponds to either strictly positive or strictly negative periodic orbit (the blue solid lines).

\item $E_\beta$ with $\beta > 0$ corresponds to a sign-indefinite periodic orbit outside the homoclinic loops (the red solid line).
\end{itemize}

In some explicit computations (see, e.g., \cite{BMP,marzuola} for $p = 1$ 
and \cite{NP20} for $p = 2$), it is useful to recall the exact solutions 
of the differential equation \eqref{stat-NLS-halfedge}. For $p = 1$, the dnoidal Jacobi elliptic solutions represent the sign-definite periodic orbits and 
the cnoidal Jacobi elliptic solutions represent the sign-indefinite periodic orbits 
in Figure \ref{fig-phase-portrait}. 

Up to translation, writing a dnoidal Jacobi elliptic solutions as 
$u(z) = A\, {\rm dn}(Bz; k)$ with elliptic modulus $k \in (0,1)$ and constants $A, B$ gives $v(z) = -ABk^2 {\rm sn}(Bz; k) {\rm cn}(Bz; k)$. Substituting $u$ and $v$ to the relation $v^2 - u^2 + u^4 = \beta$ and using identities for the Jacobi elliptic functions yield the exact solution:
\begin{equation}
\label{dnoidal}
u_{\rm dn}(z) = \frac{1}{\sqrt{2-k^2}} {\rm dn} \left(\frac{x}{\sqrt{2-k^2}}; k\right), \quad \beta = \frac{k^2-1}{(2-k^2)^2}.
\end{equation}
Since $\beta \in \left(-\frac{1}{4}, 0\right)$, the dnoidal solution 
(\ref{dnoidal}) corresponds to a level curve $E_{\beta}$ inside the positive homoclinic orbit.

Similar computations for cnoidal waves imply that 
\begin{equation}
\label{cnoidal}
u_{\rm cn}(z) = \frac{k}{\sqrt{2k^2-1}} {\rm cn} \left(\frac{x}{\sqrt{2k^2-1}}; k\right), \quad \beta = \frac{(1-k^2)k^2}{(2k^2-1)^2}.
\end{equation}
Since $\beta > 0$, the cnoidal solution corresponds to a level curve $E_\beta$ outside the homoclinic orbits. Note that the parameter $k$ takes values in $(\frac{1}{\sqrt 2}, 1]$ and that the cnoidal solution $u_{\rm cn}(z)$ can be obtained from the dnoidal solution $u_{\rm dn}(z)$ by the Jacobi real transformation 
$$
{\rm dn}(x;k) = {\rm cn}(k x;k^{-1}).
$$

The main disadvantage of Jacobi elliptic functions is that they can only be applied to the cases $p = 1$ and $p = 2$. In comparison, the method relying 
on the period function can be applied to the differential equation (\ref{stat-NLS-halfedge}) for every $p > 0$.

\begin{figure}[htbp] 
	\centering
	\includegraphics[width=3.5in, height = 2.5in]{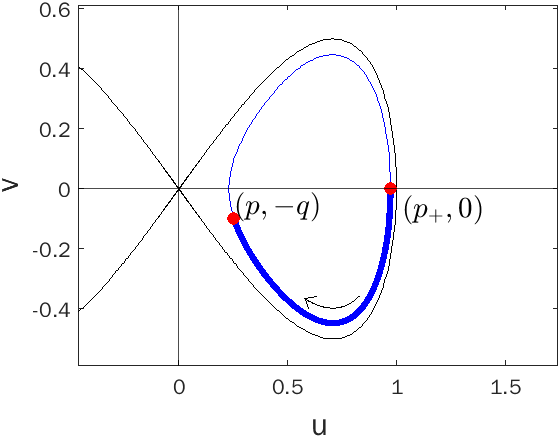}
	\caption{Phase plane $(u,v)$ for the differential equation $-u'' + u - 2 u^3 = 0$ showing the homoclinic loop for $\beta = 0$ and the integral curve $E_{\beta}$ for $\beta \in (-\frac{1}{4},0)$. The part of the integral curve between points $(\mathfrak{p}_+,0)$ and $(\mathfrak{p},-\mathfrak{q})$ corresponds to the unique decreasing solution with the boundary conditions \eqref{endpoint-values-of-u}.}
	\label{fig-nls-soliton-like-phase-plane}
\end{figure}

Let us now introduce {\em the period function} for the periodic orbit 
on the level curve $E_{\beta}$. Given a point $(\mathfrak{p},\mathfrak{q})$ on the level curve $E_{\beta}$, we define the period function as follows:
\begin{equation}
\label{period}
T_+(\mathfrak{p},\mathfrak{q}) := \int_{\mathfrak{p}}^{\mathfrak{p}_+} \frac{du}{\sqrt{\beta + A(u)}}, \quad \text{where} \quad A(u) := u^2-u^{2p+2}.
\end{equation}
The period function gives the $z$-length of the solution $u(z)$ obtained along the level curve $E_{\beta}$ between the two points (\ref{endpoint-values-of-u}). 
Figure \ref{fig-nls-soliton-like-phase-plane} shows the two points and the part of the invariant curve between these points on the phase plane $(u,v)$. 
The point $(\mathfrak{p}_+,0)$ represents a turning point for the periodic orbit defined by
$\beta + A(\mathfrak{p}_+) = 0$. Since $(\mathfrak{p},\mathfrak{q})$ and $(\mathfrak{p},-\mathfrak{q})$ also belong to the level curve $E_{\beta}$, they are related by $\beta + A(\mathfrak{p}) = \mathfrak{q}^2$.

Since the pulse profile $u$ should be located on the interval $[0, \epsilon \ell]$, the period function defines another relation between $\mathfrak{p}$ and $\mathfrak{q}$:
\begin{equation}
\label{stat-NLS-alternative}
T_+(\mathfrak{p},\mathfrak{q}) = \epsilon \ell.
\end{equation}
The pulse solution $u$ can be constructed using the following steps: for every point $(\tilde{\mathfrak{p}}, - \tilde{\mathfrak{q}})$ on the same level curve $E_{\beta}$ between points $(\mathfrak{p}_+, 0)$ and $(\mathfrak{p},-\mathfrak{q})$, we have $\tilde z = T_+(\tilde{\mathfrak{p}}, \tilde{\mathfrak{q}}) \in (0,\epsilon \ell)$ such that $u(\tilde z) = \tilde{\mathfrak{p}}$ and $u'(\tilde z) = -\tilde{\mathfrak{q}}$.

\subsection{Properties of the period function}

The period function $T_+(\mathfrak{p},\mathfrak{q})$ defined in \eqref{period} satisfies certain monotonicity properties which are useful in analysis of pulses on edges of a metric graph $\GG$. These results were obtained in \cite{KMPZ,KP21} in the particular cases of cubic nonlinearity, which are reviewed here.  

Recall that $A(u) := u^2 - u^4$ attains a minimum at $u = \pm \mathfrak{p}_*$, where $\mathfrak{p}_* = \frac{1}{\sqrt 2}$. Properties of $T_+(\mathfrak{p},\mathfrak{q})$ depend whether $\mathfrak{p} \in (0, \mathfrak{p}_*]$ or $\mathfrak{p} \in (\mathfrak{p}_*, 1)$. The following two propositions were proven in Section 3.3 in \cite{KMPZ}.

\begin{proposition}
\label{monotone-on-left}
For every $\mathfrak{p} \in (0, \mathfrak{p}_*]$, $T_+(\mathfrak{p},\mathfrak{q})$ is a monotonically decreasing function of $\mathfrak{q}$ in $(0,\infty)$.
\end{proposition}

\begin{proposition}
\label{nonmonotone-on-right}
For every $\mathfrak{p} \in (\mathfrak{p}_*, 1)$,
$T_+(\mathfrak{p},\mathfrak{q})$ is a non-monotone function of $\mathfrak{q}$ in $(0,\infty)$ satisfying 
$T_+(\mathfrak{p},\mathfrak{q}) \to 0$ as $\mathfrak{q} \to 0$ and $\mathfrak{q} \to \infty$.
\end{proposition}

\begin{figure}[htbp] 
\centering
\includegraphics[width=3in, height = 2in]{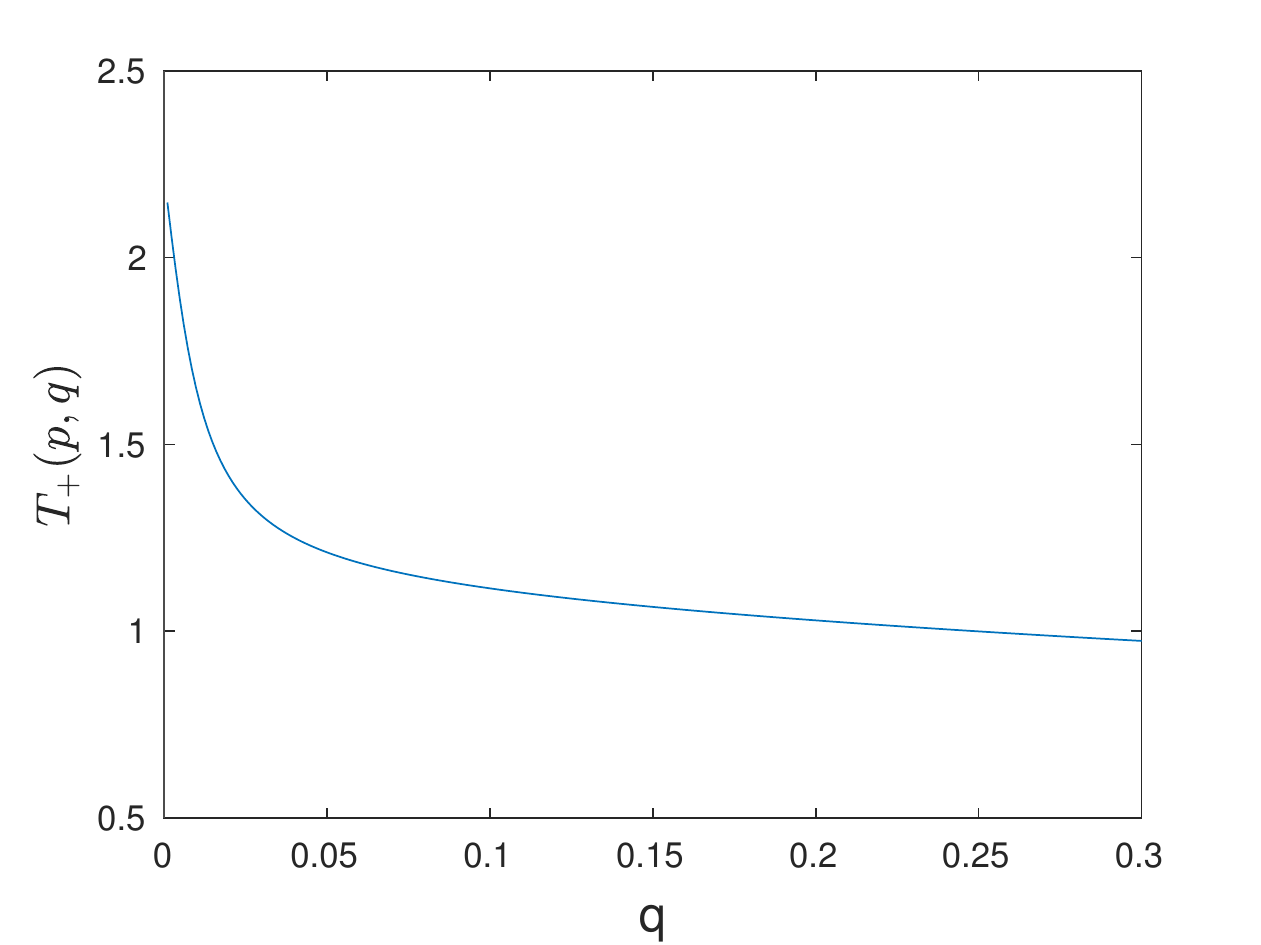}
\includegraphics[width=3in, height = 2in]{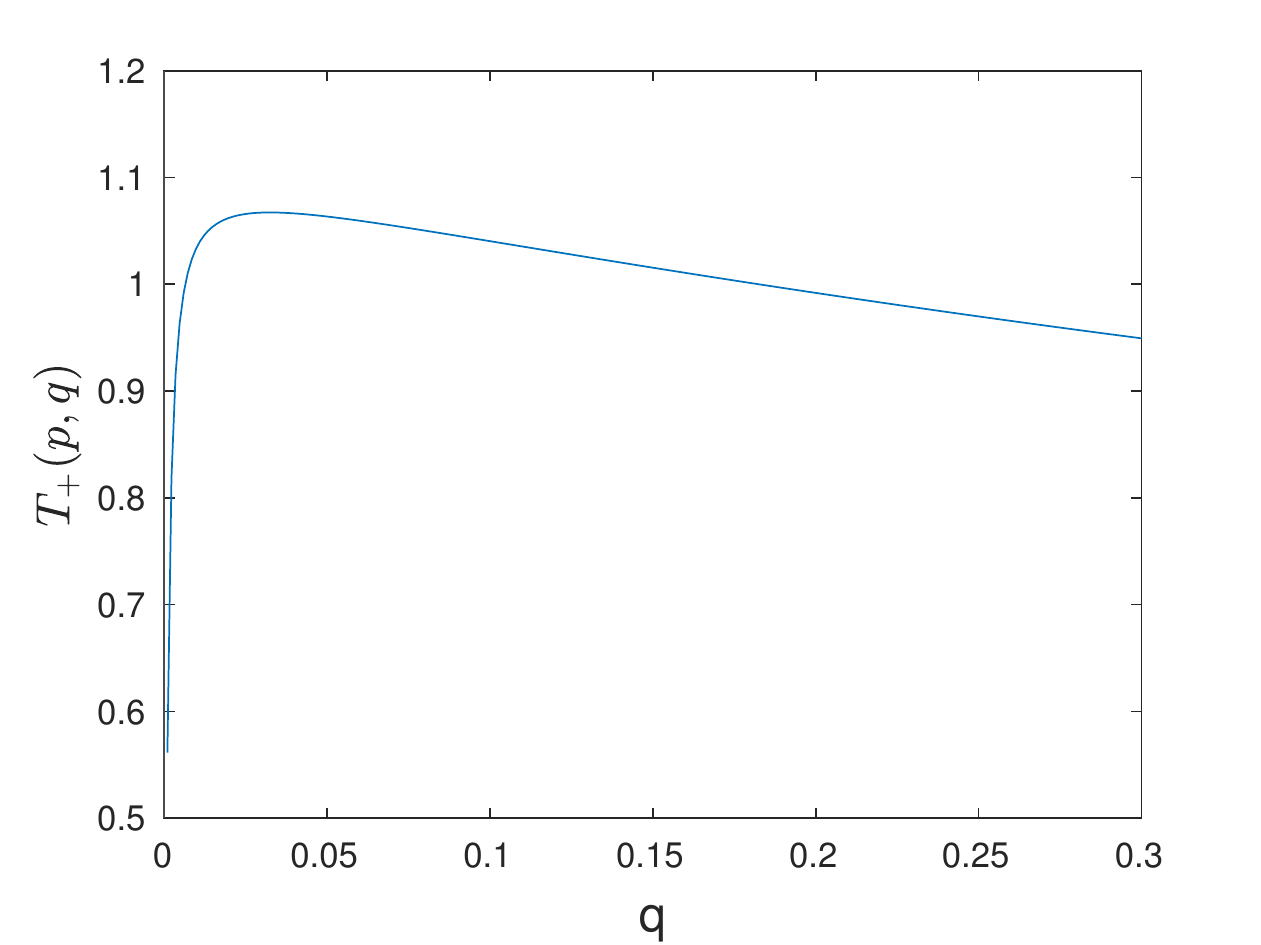}
\caption{The dependence of $T_+(\mathfrak{p},\mathfrak{q})$ versus $\mathfrak{q}$ for $\mathfrak{p} \in (0,\mathfrak{p}_*)$ (left) and $\mathfrak{p} \in (\mathfrak{p}_*, 1)$ (right).}
\label{period-function-mon-fig}
\end{figure}

Figure \ref{period-function-mon-fig} illustrates the monotonicity result of Proposition \ref{monotone-on-left} (left) and the non-monotonicity result of Proposition \ref{nonmonotone-on-right} (right). 

In the case of non-monotonicity, Figure \ref{period-function-mon-fig} (right) suggests that there exists exactly one critical point of $T_+(\mathfrak{p},\mathfrak{q})$ with respect to $\mathfrak{q}$ for a fixed value of $\mathfrak{p} \in (\mathfrak{p}_*, 1)$, and it is the maximum point. 
Let us denote the corresponding value of $\mathfrak{q}$ by $\qmax(\mathfrak{p})$. This maximum point exists for every $\mathfrak{p} \in (\mathfrak{p}_*, 1)$, and since the period function $T_+(\mathfrak{p},\mathfrak{q})$ is a continuous differentiable function of $(\mathfrak{p},\mathfrak{q})$ in $(0,1) \times (0,\infty)$,
the value $\qmax$ is a continuous function of $\mathfrak{p}$. 

\begin{figure}[htbp]
	\includegraphics[width=3in, height = 2in]{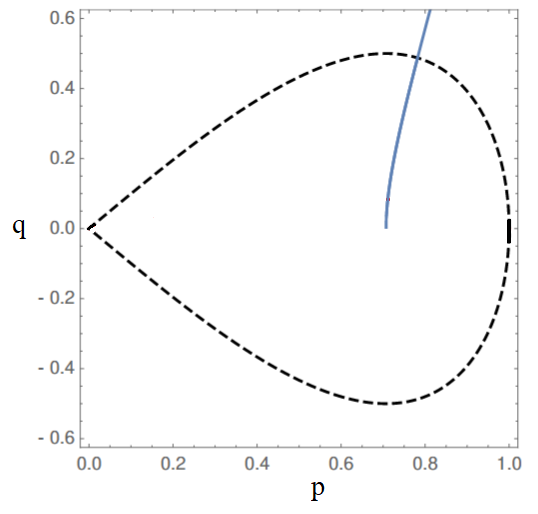}
	\caption{Dependence of $\qmax$ (the maximum of $T_+(\mathfrak{p},\mathfrak{q})$ in $\mathfrak{q}$) versus $\mathfrak{p}$ is shown by the blue solid line. The black dashed line shows the homoclinic orbit.}
	\label{fig:tmax}
\end{figure}

Figure \ref{fig:tmax} shows that the value $\qmax$ (the blue solid line) is monotonically increasing function of $\mathfrak{p}$. The black dashed line displays the homoclinic orbit at the energy level $E_0$ at $\beta = 0$ so that 
$\qmax$ crosses the homoclinic orbit at the unique value of $\mathfrak{p}_{**} \approx 0.782$. This implies that, for every $\mathfrak{p} \in (\mathfrak{p}_{**}, 1)$, the period function $T_+(\mathfrak{p},\mathfrak{q})$ is a monotonically increasing function of $\mathfrak{q}$ inside the homoclinic orbit. The following proposition was also proven in \cite{KMPZ}.

\begin{proposition}
\label{behaviour-qmax}
There exists $\mathfrak{p}_{**} \in (\mathfrak{p}_*, 1)$ such that for every $\mathfrak{p} \in (\mathfrak{p}_*, \mathfrak{p}_{**})$, there is exactly one critical (maximum) point of $T_+(\mathfrak{p},\mathfrak{q})$ in $\mathfrak{q}$
inside $(0, \sqrt{A(\mathfrak{p})})$, where $A(\mathfrak{p}) := \mathfrak{p}^2 - \mathfrak{p}^4$. For $\mathfrak{p} \in (\mathfrak{p}_{**}, 1)$, $T_+(\mathfrak{p},\mathfrak{q})$ is monotonically increasing in
$\mathfrak{q}$ inside $(0, \sqrt{A(\mathfrak{p})})$.
\end{proposition}

Finally, the asymptotic representation
$$
T_+(\mathfrak{p},\mathfrak{q}) \sim -\ln(\mathfrak{p} + \mathfrak{q}), 
\quad \mbox{\rm as} \;\; \mathfrak{p},\mathfrak{q} \to 0
$$ 
proven in \cite{BMP} (see also \cite{KP21}) suggests the following property, which is important for the study of pulses in the limit of large $\epsilon$ in $\omega = -\epsilon^2$.

\begin{proposition}
	\label{period-monotonicity-in-the-limit}
	There is a small $\delta>0$ such that
	\begin{equation*}
	\partial_{\mathfrak{p}} T_+(\mathfrak{p},\mathfrak{q}) < 0 \quad \text{and} \quad \partial_{\mathfrak{q}} T_+(\mathfrak{p},\mathfrak{q}) < 0
	\end{equation*}
	for every $\mathfrak{p} \in (0, \delta)$ and $\mathfrak{q} \in (0, \delta)$. 
\end{proposition}

\subsection{Example of the tadpole graph}

Here we show how the period function $T_+(\mathfrak{p},\mathfrak{q})$ can be applied to obtain precise analytical results on existence of single-pulse states on a tadpole graph with a single loop (see Figure \ref{fig-graphs}). The arguments can be extended to the flower graph with $N$ loops in \cite{KMPZ}. 

We parameterize the unbounded edge of the tadpole graph $\GG$ by $[0, \infty)$ and the loop by $[-\pi, \pi]$. After the scaling transformation (\ref{scaling-transform}), the stationary NLS equation (\ref{stat-NLS-halfedge})
with the cubic nonlinearity can be written as the following system
of differential equations:
\begin{equation}
\label{NLS-scaled-tadpole}
\left\{ \begin{array}{l} -u_1''(z) + u_1(z) - 2 u_1(z)^3 = 0, \quad z \in (-\pi \epsilon,\pi \epsilon),\\
-u_0''(z) + u_0(z) - 2 u_0(z)^3 = 0, \quad z > 0,\\
u_1(\pi \epsilon) = u_1(-\pi \epsilon) = u_0(0), \\
 u_1'(\pi \epsilon) - u_1'(-\pi \epsilon) = u_0'(0).
\end{array} \right.
\end{equation}
The only dependence
of system (\ref{NLS-scaled-tadpole}) on $\epsilon$ is due to the length of the interval $[-\pi \epsilon,\pi \epsilon]$. 

We are looking for positive decaying solution to equation 
$-u_0''(z) + u_0(z) - 2 u_0(z)^3 = 0$ on the half-line, which is 
the shifted NLS soliton $u_0(z) = {\rm sech}(z + a)$, 
where $a \in \mathbb{R}$ is an arbitrary positive translation parameter. 
In general, if one is interested in other bound states, then $a$ can be chosen to be negative, which leads to non-monotone $u_0$  on $[0,\infty)$.

After eliminating $u_0(z) = {\rm sech}(z+a)$, we obtain 
the closed boundary-value problem for a single looping edge $[-\pi \epsilon,\pi \epsilon]$. However, since the pulse is symmetric on the looping edge, 
the boundary-value problem can be reduced on the half-interval:
\begin{equation}
\left\{ \begin{array}{l} -u_1''(z) + u_1(z) - 2 u_1(z)^3 = 0, \quad z \in (0, \pi \epsilon), \\
u_1(0) = \mathfrak{p}_+, \quad u_1(\pi \epsilon) = \mathfrak{p}, \\
u_1'(0) = 0, \quad  u_1'(\pi \epsilon) = -\frac{1}{2} \sqrt{A(\mathfrak{p})},
\end{array} \right.
\label{bvp-u1-tadpole}
\end{equation}
where $\mathfrak{p} = u_0(0) = \sech (a) \in (0,1)$ is a free parameter 
of the boundary-value problem, whereas $\mathfrak{p}_+$ is computed from the 
turning point satisfying 
$$
 - A(\mathfrak{p}_+) = \mathfrak{q}^2 - A(\mathfrak{p}) = -\frac{3}{4} A(\mathfrak{p})
$$ 
with $\mathfrak{q} := \frac{1}{2} \sqrt{A(\mathfrak{p})}$  and $A(\mathfrak{p}) := \mathfrak{p}^2 - \mathfrak{p}^4$ since $u_0'(0) = -\sqrt{A(\mathfrak{p})}$ is determined on the homoclinic orbit $E_0$ with $\beta = 0$. 

To solve the boundary-value problem (\ref{bvp-u1-tadpole}), we use the period function $T_+(\mathfrak{p},\mathfrak{q})$ introduced in \eqref{period}. On the interval $[0, \pi \epsilon]$, the system \eqref{bvp-u1-tadpole} is equivalent to \eqref{stat-NLS-halfedge} and \eqref{endpoint-values-of-u} with $\ell = \pi$ and  $p=1$. The existence of a single-pulse edge-localized state on $\GG$ depends on the existence of a root to the  nonlinear equation \eqref{stat-NLS-alternative}, which becomes
\begin{equation}
\label{tadpole-period}
T_+\left( \mathfrak{p},\frac{1}{2} \sqrt{A(\mathfrak{p})}\right) = \epsilon \pi.
\end{equation}
It was proven in \cite{KMPZ} that the mapping $\mathfrak{p} \to T_+\left( \mathfrak{p},\frac{1}{2} \sqrt{A(\mathfrak{p})}\right)$ is monotonically decreasing for every $\mathfrak{p} \in (0,1)$ and has the asymptotic limits: 
$$
\lim_{\mathfrak{p} \to 0} T_+\left( \mathfrak{p},\frac{1}{2} \sqrt{A(\mathfrak{p})}\right) = \infty, \qquad 
\lim_{\mathfrak{p} \to 1} T_+\left( \mathfrak{p},\frac{1}{2} \sqrt{A(\mathfrak{p})}\right) = 0. 
$$
This means that for every positive $\epsilon$ there is an unique root $\mathfrak{p}$ to the nonlinear equation \eqref{tadpole-period} so that the following result holds. 

\begin{proposition}
\label{global-existence}
For every $\epsilon>0$, there exists only one ingle-pulse edge-localized solution to  the stationary NLS equation \eqref{NLS-scaled-tadpole} on the tadpole graph. Such solution is symmetric on the loop parameterized
by $[-\pi \epsilon,\pi \epsilon ]$, and is monotonically decreasing on the tail $[0,\infty)$.
\end{proposition}

\begin{figure}[htbp] 
   \centering
  \includegraphics[width=2.8in, height = 2in]{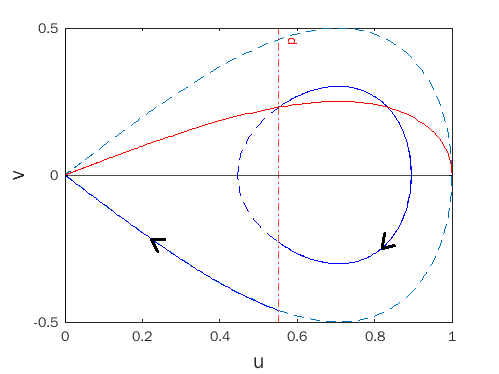}
    \includegraphics[width=2.8in, height = 2.1in]{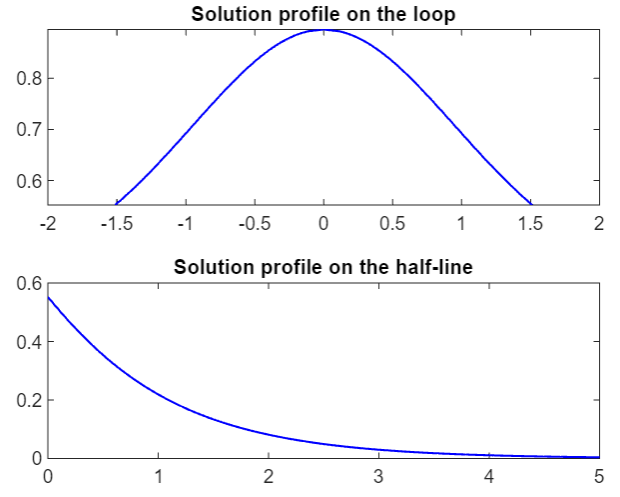}
   \caption{Left: Geometric construction of the single-pulse state on the phase plane for the tadpole graph. Right: Profile of solutions on the loop and on the half-line.}
   \label{fig-plane-tadpole}
\end{figure}

Figure \ref{fig-plane-tadpole} shows a geometric construction of solutions to the boundary-value problem (\ref{NLS-scaled-tadpole}) on the phase plane $(u,v)$.
The dashed line represents the homoclinic orbit at $E_0 = 0$ with the solid part depicting the shifted NLS soliton $u_0(z) = \sech (z+a)$. The dashed-dotted
vertical line depicts the value of $\mathfrak{p} = u_0(0) = {\rm sech}(a)$. The level curve $E(\mathfrak{p},\frac{1}{2} \sqrt{A(\mathfrak{p})}) = v^2 - A(u)$  is shown by the dashed line, whereas the solid part depicts a suitable solution to the boundary-value problem (\ref{bvp-u1-tadpole}).

\begin{figure}[htbp] 
	\centering
	\includegraphics[width=2.5in, height = 2.3in]{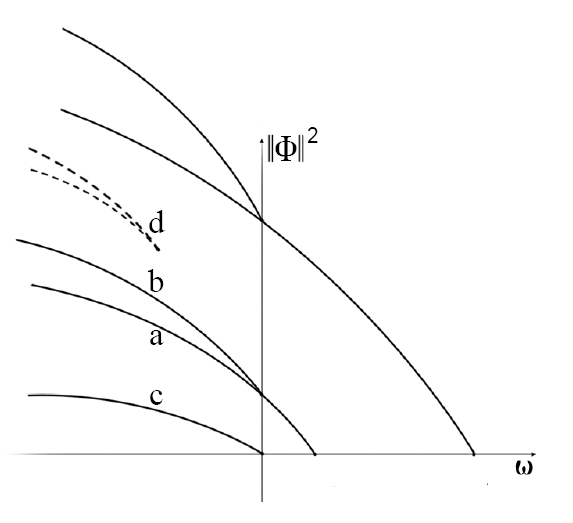}
	\caption{Bifurcation diagram for bound states on the tadpole graph. The curve $c$ corresponds to a family of single-pulse states for $\omega<0$ which bifurcate from the trivial solution at $\omega = 0$. The rest of the curves in the diagram correspond to sign-indefinite bound states.}
	\label{fig-bifurcation}
\end{figure}

First results about bound states of the cubic NLS equation on the tadpole graph were obtained in \cite{cfn15} by using the Jacobi elliptic functions (\ref{dnoidal}) and (\ref{cnoidal}) on the looping edge. Parameters of the Jacobi elliptic functions and their translations were found 
from the boundary conditions in the boundary-value problem (\ref{NLS-scaled-tadpole}). Bifurcation diagram on Figure \ref{fig-bifurcation} 
shows all branches of bound states. Branch $c$ for $\omega < 0$ 
corresponds to the single-pulse state of Proposition \ref{global-existence}.

The branch $a$ corresponds to a sign-indefinite solution in the looping edge with zero 
tail $u_0 = 0$. This and similar solutions are continued from $\omega < 0$ to $\omega > 0$ 
up to the points of their bifurcations from the linear boundary-value problem. At $\omega = 0$, the branch $a$ splits with a new branch $b$ which exists for $\omega < 0$ and consists 
of a positive exponentially decaying tail $u_0$ and a sign-indefinite 
solution on the looping edge. Bifurcations of such solutions for $\omega < 0$ were studied in \cite{nps} for the general power nonlinearities by using the Lyapunov--Schmidt reduction method.

Finally, there exist pairs of branches $d$, for which the solution in the looping edge consists of several pulses, whereas the exponentially decaying tails $u_0$ are either monotone or non-monotone. These excited states 
are usually not interesting in applications because they are typically unstable in the time evolution of the NLS equation \cite{nps}.

Similar classifications of bound states were obtained from explicit analysis of the stationary NLS equations with cubic nonlinearity on the dumbell graph in \cite{G19,marzuola}, the truncated star graph \cite{Band}, and the periodic  graph \cite{PS17}. Positive bound states on the tadpole graph 
in the case of quintic nonlinearity were considered with the period 
function in \cite{NP20}.

Stability of the ground state of the NLS equation on the tadpole graph with Neumann--Kirchhoff condition has been studied for the subcritical and critical case in \cite{nps} and \cite{NP20}, respectively. Branch $c$ on Figure \ref{fig-bifurcation} was shown to be orbitally stable in \cite{nps} 
for $0 < p < 2$. Branch $b$ was shown to be orbitally unstable in \cite{nps}, however, no conclusion on orbital stability of branch $a$ was made. In the critical case $p=2$, orbital stability of the ground state 
of the variational problem (\ref{infB}) depends on the values of $\omega$  \cite{NP20}. It was shown that there exists $\omega_* < 0$ such that 
the ground state is orbitally stable for $\omega \in (\omega_*,0)$ and orbitally unstable for $\omega \in (-\infty,\omega_*)$ \cite{NP20}.
Moreover, the ground state of the variational problem (\ref{infB}) along the solution branch undergoes a transition from the ground state of the variational problem (\ref{minimization}), to a local constrained minimum, and eventually to a saddle point.

\section{Bound states in the limit of large mass}
\label{sec-multi-pulse}

The question of existence of edge-localized states to the stationary NLS equation becomes more complicated when we consider a general graph $\GG$ 
where many pulses can be placed on different bounded edges. 
Although each pulse still corresponds to a level curve $E_\beta$, 
the value of $\beta$ can be different between different bounded edges. 
Even in the case of cubic nonlinearity, when each pulse is expressed 
by the dnoidal or cnoidal Jacobi elliptic functions, it is hard to adjust 
their parameters to satisfy the boundary conditions on different vertices of the general graph $\GG$. 

Another method can be used to overcome this difficulty in the limit 
of large negative $\omega$, which in the case of subcritical nonlinearities 
transfers to the limit of large mass $M(\Phi) = \mu$ in the variational 
problem (\ref{minimization}). The positive edge-localized states were considered in the limit of large mass in \cite{AST19,BMP,D20,KP21,KS20}.

\subsection{Single-pulse states}

Single-pulse states were constructed in \cite{AST19} by solving the energy minimization problem 
\begin{equation}
\label{min-another}
E(\Phi) = \inf\limits_{u \in V_\mu(\GG)} E(\Psi)
\end{equation} 
in the space of functions of fixed mass $M(\Psi) = \mu$ which additionally attain their maximum on a particular bounded edge $e$ of the graph $\GG$:
\begin{equation}
\label{set}
V_\mu(\GG) = \left\{ \Psi \in H^1_C(\GG): \quad M(\Psi) = \mu \;\; \text{and} \;\; \| \Psi \|_{L^\infty(\GG)} = \| \psi_e \|_{L^\infty(e)}   \right\}.
\end{equation}
It was proven in \cite{AST19} for the subcritical powers $0 < p < 2$ 
and for unbounded graphs $\GG$ that for sufficiently large values of $\mu$, there exists a positive solution $\Phi \in V_{\mu}(\GG)$ to the energy minimization problem (\ref{min-another}) and (\ref{set}). By construction, $\Phi$ is a single-pulse state as it reaches its maximal value on the edge $e$ only, and such $\Phi$ can be found for every chosen edge $e$. For example, if $\GG$ has $k$ bounded edges, then for sufficiently large $\mu$, $\GG$ admits $k$ single-pulse states of sufficiently large mass $\mu$ each localized on a single edge of $\GG$. Each bound state $\Phi \in V_{\mu}(\GG)$ is a local minimizer of energy $E(\Psi)$ under the fixed mass constraint $M(\Psi) = \mu$ such that {\em the Morse index}, which is the number of negative eigenvalues of the linearized operator $\mathcal{L}_+$ introduced in (\ref{L-plus}), is equal to {\em one}. Consequently, the bound state $\Phi \in V_{\mu}(\GG)$ is orbitally stable in the time evolution, according to Proposition \ref{prop-stability}(i), because the slope condition $\frac{d}{d\omega} \| \Phi \|^2_{L^2(\GG)} < 0$ is always satisfied for sufficiently large negative $\omega$ for the subcritical powers $0 < p < 2$. 

Existence of single-pulse states was proven in \cite{BMP} for both bounded and unbounded graphs $\GG$ but only in the case of cubic nonlinearity. 
The idea explored in \cite{BMP} relies on the {\em Dirichlet-to-Neumann (DtN) mappings} constructed separately for two parts of the graph $\GG$: a single bounded edge $e$ and the rest of the graph $\GG \setminus \{e\}$. The single-pulse states are then obtained by studying the intersections of the DtN mappings corresponding to the two parts of a graph $\GG$. Moreover, the same asymptotic method enables classification of the single-pulse states of the least energy at fixed mass $M(\Phi) = \mu$ by comparing the exponentially small terms in the expansion of $\mu := Q(\Phi)$ in $\epsilon := \sqrt{|\omega|}$. 
The following result was obtained in \cite{BMP}.

\begin{proposition}
If $e$ is either a looping or internal 
edge of length $2 \ell$ or a pendant of length $\ell$ as on Figure \ref{fig-states}, then for sufficiently large $\epsilon$, there exists 
the single-pulse edge-localized state $\Phi \in H^1_C(\GG)$ such that the component $\phi_e$ on edge $e$ attains a single local maximum on $e$, monotone from its maximum to the vertices of $e$, and 
concentrated on $e$ in the following sense
\begin{equation}
\label{eq:L2_proportion}
\frac{\left\|\Phi \right\|_{L^2(e)}}
{\left\|\Phi\right\|_{L^2(\Gamma)}}
\geq 1-Ce^{-2\epsilon\ell},
\end{equation}
where the constant $C$ is independent of $\epsilon$. The profile $\Phi$ has no internal maxima on the remainder of graph $\Gamma \backslash \{e\}$. 
The state of the least energy at fixed (large) mass $\mu$ localizes on 
the following edges of the graph $\Gamma$:
\begin{itemize}
	\item[(i)] The longest among pendants; in the case of a tie, the pendant incident to fewest edges.
	\item[(ii)]  If (i) is void, the shortest among loops incident to a single edge.
	\item[(iii)] If (i)--(ii) are void, a loop incident to two edges.
	\item[(iv)] If (i)--(iii) are void, the longest edge among the following: a looping edge incident to  more than two edges or an internal edge
	incident to more than one other edge.
\end{itemize}
\label{prop-single-pulse}
\end{proposition}

Additional results related to the single-pulse states in the limit of large mass were obtained in \cite{D20}, where bounded graphs with the pendant edges were considered and convergence of the edge-localized states to the half-solitons was proven. Multi-pulse states were also studied in \cite{D20}, all pulses localize at the terminal vertices. 

Single-pulse states were considered in \cite{KS20} by recasting the existence problem to the semi-classical limit of an elliptic problem. It was proven in \cite{KS20} that the location of the positive single-pulse state with a single maximum as the state of the least energy at fixed mass is determined by the longest pendant edge of a bounded graph or the longest internal edge if no pendant and looping edges are present. These results are included in more general results of Proposition \ref{prop-single-pulse}.

\subsection{Multi-pulse states} 

Multi-pulse states were classified for the stationary NLS equation with the cubic nonlinearity in the limit of large mass under some consistency assumptions in \cite{KP21} based on the properties of the DtN mappings proven in \cite{BMP}. 

In order to illustrate these results, we consider the cubic NLS equation on bounded edges of a metric graph:
\begin{equation}
\label{nls-scaled}
-\phi_e''(z) + \phi_e(z) - 2 \phi_e(z)^3 = 0, \quad 
z \in e, 
\end{equation}
where $\phi_e$ is the corresponding component of $\Phi$ on the edge $e$ 
and the scaling transformation (\ref{scaling-transform}) with $p = 1$ has already been used. 
The original graph $\GG = \{ E, V \}$ is transformed to the 
$\epsilon$-scaled graph $\GG_{\epsilon} = \{ E_{\epsilon}, V \}$ 
for which every bounded edge $e \in E$ of length $\ell_e$ 
transforms to the edge $e_{\epsilon} \in E_{\epsilon}$ of length $\epsilon \ell_e$ but the unbounded edge $e \in E$ remains the same as $e \in E_{\epsilon}$. 

As we are interested in the multi-pulse edge-localized states on $\GG$, we partition the graph into two parts: one part containing $N$ bounded edges $E_N := \{ e_1,e_2,\dots, e_N\}$ which will contain the pulse profiles, while the other part is the complement $\GG \setminus E_N$ which will have small components. Similar notations are used for $\GG_{\epsilon}$, $E_{N,\epsilon} \subset E_{\epsilon}$, and $e_{j,\epsilon} \in E_{N,\epsilon}$.

Vertices in the set $E_{N,\epsilon}$ which are connected to $\Gamma_{\epsilon} \backslash E_{N,\epsilon}$ are declared as the boundary vertices 
and the following {\em Dirichlet data}
$\vec{p} = (p_1,p_2, \ldots, p_{|B|})$ and {\em Neumann data}
$\vec{q} = (q_1,q_2, \ldots, q_{|B|})$ are introduced as follows:
\begin{equation}
\label{Neumann_data}
p_j := u_{e \sim v_j}(v_j), \quad 
q_j := \sum_{e \sim v_j} \partial u_e(v_j), \quad v_j \in V_B,
\end{equation}
where $V_B$ is the set of $|B|$ boundary vertices,
the derivatives $\partial$ are directed away from $\Gamma_{\epsilon} \backslash E_{N,\epsilon}$ and $e \sim v_j$ lists all edges $e \in \Gamma_{\epsilon} \backslash E_{N,\epsilon}$ incident to the vertex $v_j$. Pulse solutions 
on $E_{N,\epsilon}$ are required to satisfy
\begin{equation}
\label{bound-up}
\sup_{z \in e_{j,\epsilon}} |u_j(z)| > \frac{1}{\sqrt{2}}, \quad 
e_{j,\epsilon} \in E_{N,\epsilon}
\end{equation}
whereas the small remainder terms on $\GG_{\epsilon} \backslash E_{N,\epsilon}$ 
are required to satisfy
\begin{equation}
\label{bound-down}
\sup_{z \in e_{j,\epsilon}} |u_j(z)| < \frac{1}{\sqrt{2}}, \quad 
e_{j,\epsilon} \in \Gamma_{\epsilon} \backslash E_{N,\epsilon},
\end{equation}
where $\frac{1}{\sqrt{2}}$ is the 
constant solution of the differential equations in (\ref{nls-scaled}).

The following two results were proven in \cite{BMP} (Theorem 2.9 and Lemma 2.12) and have been used for analysis in \cite{KP21}. The first result concerns to the elliptic estimates of small solutions to the stationary NLS equation (\ref{nls-scaled}) on the remainder graph $\Gamma_{\epsilon} \backslash E_{N,\epsilon}$. The second result gives estimates of the large solution 
to the stationary NLS equation (\ref{nls-scaled}) on just one edge 
$e \in E_{N,\epsilon}$.

\begin{proposition}
	\label{lemma-1}
	There exist $C_0 > 0$, $p_0 > 0$, and $\epsilon_0>0$ such that for every
	$\vec{p}$ with $\| \vec{p}\| < p_0$ and every
	$\epsilon > \epsilon_0$, there exists a solution $U \in H^2_{\rm NK}(\Gamma_{\epsilon} \backslash E_{N,\epsilon})$
	to the stationary NLS equation (\ref{nls-scaled}) 
	on $\Gamma_{\epsilon} \backslash E_{N,\epsilon}$ subject 
	to the Dirichlet data on $V_B$ which
	is unique among functions satisfying \eqref{bound-down}.
	The solution satisfies the estimate
	\begin{equation}
	\label{eq:nlin_DTN_solution}
	\| U \|_{H^2(\Gamma_{\epsilon} \backslash E_{N,\epsilon})} \leq C_0 \|\vec{p}\|,
	\end{equation}
	while its Neumann data satisfies
	\begin{equation}
	\label{eq:nlin_DTN_value}
	|q_j - D_j p_j| \leq C_0 \left( \|\vec{p} \| e^{-\epsilon \ell_{\rm min}} +
	\|\vec{p}\|^3 \right),	\qquad 1 \leq j \leq |B|,
	\end{equation}
	where $D_j$ is the degree of the $j$-th boundary vertex 
	in $\Gamma_{\epsilon} \backslash E_{N,\epsilon}$ and $\ell_{\rm min}$
	is the length of the shortest edge in $\Gamma \backslash E_N$.
\end{proposition}

\begin{proposition}
	\label{lemma-2}
	There exist $C_0 > 0$, $p_0 > 0$, and $\epsilon_0>0$ such that for every
	$p \in (0,p_0)$ and every $\epsilon > \epsilon_0$, there exists a real solution $u \in H^2(0,\epsilon \ell)$ to the differential equation 
	$$
	-u'' + u - 2 u^3 = 0, \quad 0 < z < \epsilon \ell,
	$$ 
	satisfying $u'(0) = 0$ and $u(\epsilon \ell) = p$, which is unique among positive and decreasing functions satisfying \eqref{bound-up}.	The solution satisfies $u'(\epsilon \ell) < 0$ and
	\begin{equation}
	\left| u'(\epsilon \ell) - u(\epsilon \ell) + 4 e^{-\epsilon \ell} \right|
	\leq C_0 \epsilon e^{-3 \epsilon \ell}.
	\label{eq:single-bump-DtN}
	\end{equation}
\end{proposition}

For each boundary vertex $v_j \in V_B$ with $1 \leq j \leq |B|$, 
we use the Dirichlet data $p_j$ as the unknown variable
and write the flux boundary condition to determine $p_j$. 
The main advantage of this method is that the value of $p_j$ 
can be found independently from the other boundary vertices under some consistency assumptions. For simplicity 
of presentation, let us only consider the case when 
the bounded edges in $E_{N}$ are only represented by 
the looping edges.

By Proposition \ref{lemma-1}, the Neumann data at the boundary vertex $v_j$ directed away from $\Gamma_{\epsilon} \backslash E_{N,\epsilon}$ is 
\begin{equation}
q_j^{(1)} = D_j p_j + \mathcal{O}(\| \vec{p} \| e^{-\epsilon \ell_{\rm min}} 
+ \| \vec{p} \|^3).
\end{equation}
By Proposition \ref{lemma-2}, the Neumann data at the same boundary vertex $v_j$ directed away from $E_{N,\epsilon}$ is 
\begin{eqnarray}
q_j^{(2)} = 2 L_j p_j - 8 \sum_{e \in E_{j,\epsilon}} e^{-\epsilon \ell_e}  
+ \mathcal{O}(\epsilon e^{-3 \epsilon \ell_{\rm j, min}}),
\end{eqnarray}
where $E_{j,\epsilon}$ is a subset of $E_{N,\epsilon}$ of looping edges 
connected to $v_j$ and $\ell_{j,{\rm min}}$ is the minimal half-length of the looping edges in $E_{j,\epsilon}$.  

The flux boundary condition gives $q_j^{(1)} + q_j^{(2)} = 0$, which becomes the implicit equation on $p_j$ with the following explicit solution:
\begin{eqnarray}
p_j = \frac{8}{Z_j}  \sum_{e \in E_{j,\epsilon}} e^{-\epsilon \ell_e} 
+ \mathcal{O}\left( \| \vec{p} \| e^{-\epsilon \ell_{\rm min}} + 
+ \| \vec{p} \|^3 + \epsilon e^{-3 \epsilon \ell_{\rm j, min}} \right),
\label{solution-pj}
\end{eqnarray}
where $Z_j$ is the total degree of the vertex $v_j$. 

\begin{remark}
The same solution (\ref{solution-pj}) works trivially if 
the vertex $v$ is not a boundary vertex in $V_B$
but an interior vertex between edges in the set $E_{N,\epsilon}$. Therefore, 
the Dirichlet data at the interior vertices are settled trivially, 
and it is only required to identify the proper solutions 
for the Dirichlet data at the boundary vertices. 
\end{remark}

Although the boundary conditions are satisfied 
for all vertices if $p_j$ is defined by (\ref{solution-pj}) 
for $1 \leq j \leq |B|$, one needs to verify that 
$\| \vec{p}\|$ in the error terms in (\ref{solution-pj})
is smaller than the leading-order terms in (\ref{solution-pj}) 
as $\epsilon \to \infty$ for each $j$. This is true if 
$$
\| \vec{p} \| e^{-\epsilon \ell_{\rm min}} \ll e^{-\epsilon \ell_{j,{\rm min}}}, 
\quad \| \vec{p} \|^3 \ll e^{-\epsilon \ell_{j,{\rm min}}}, \quad 
1 \leq j \leq |B|, 
$$
which gives the following restrictions on the lengths of edges in the set $E_N$:
	\begin{equation}
\label{constraint-on-length-1}
\max_{1 \leq j \leq |B|} \ell_{j,{\rm min}} - \min_{1 \leq j \leq |B|} \ell_{j,{\rm min}} < \ell_{\rm min}
\end{equation}
and
\begin{equation}
\label{constraint-on-length-2}
\max_{1 \leq j \leq |B|} \ell_{j,{\rm min}} <	3 \min_{1 \leq j \leq |B|} \ell_{j,{\rm min}},
\end{equation}
where $\ell_{\rm min}$ is the length of the shortest edge 
in $\Gamma \backslash E_N$. If $N = 1$, $|B| = 0$, $|B| = 1$, 
or $|B| \geq 2$ with all $\ell_{j,\rm min}$ being equal, 
the consistency assumptions (\ref{constraint-on-length-1}) and (\ref{constraint-on-length-2}) are trivially satisfied. 

With the same approach, pendant edges can be included into consideration due to the symmetry of solutions on the looping edges.  
However, for internal edges, one need to control the center of symmetry 
for pulses on the internal edges, like on Figure \ref{fig-graphs}. 
This leads to more complicated solutions of the flux boundary conditions 
and to additional restrictions that each internal edge in the set $E_N$ 
is assumed to have no common vertices with other internal edges and its half-length is strictly minimal to the half-lengths of looping edges adjustent to the two vertices of the internal edge \cite{KP21}. 

Morse index of the $N$-pulse positive state constructed from the DtN mappings 
has been computed in \cite{KP21} by using Sturm's Oscillation Theorem. 
Note that Sturm's Oscillation Theorem has been 
previously used for the star graphs in  \cite{KP2} and for flower graphs in \cite{KMPZ}. The following result was obtained in \cite{KP21}.

\begin{proposition}
Assume that $E_N$ consists of pendant and looping edges, the lengths of which 
satisfy the constraints (\ref{constraint-on-length-1}) and (\ref{constraint-on-length-2}). For  sufficiently large $\epsilon$, 
there exists the $N$-pulse edge-localized state $\Phi \in H^1_C(\GG)$ such that 
the profile $\Phi$ is concentrated on $E_N$ in the following sense 
\begin{equation}
\label{eq:L2}
\frac{\left\|\Phi \right\|_{L^2(\Gamma \backslash E_N)}}
{\left\|\Phi\right\|_{L^2(E_N)}}
\leq  C e^{-\epsilon \ell_N},
\end{equation}
where the constant $C$ is independent of $\epsilon$ 	
and $\ell_N$ is the minimum of $\{ \ell_{j,{\rm min}}\}_{j=1}^N$. 
The profile $\Phi$ has no internal maxima on the remainder of graph $\Gamma \backslash E_N$. Morse index 
of the positive $N$-pulse bound state $\Phi$ is exactly $N$.
\label{prop-multi-pulse}
\end{proposition}

The implication of Proposition \ref{prop-multi-pulse} 
to the time evolution of perturbations to 
the $N$-pulse edge-localized states is that these states are orbitally 
unstable under the NLS time flow if $N \geq 2$ according to Proposition \ref{prop-stability} (ii). In agreement with the variational 
characterization of the single-pulse states on unbounded graphs in \cite{AST19},  the single-pulse states with $N = 1$ 
are orbitally stable according to Proposition \ref{prop-stability} (i) 
because the slope condition $\frac{d}{d\omega} \| \Phi \|^2_{L^2(\GG)} < 0$ is always satisfied for sufficiently large negative $\omega$. Consequently, multi-pulse states with $N \geq 2$ 
cannot be the ground state of the variational problem (\ref{minimization}) because their Morse index exceeds {\em one}. It is still unclear (see discussion in \cite{BMP}) if the ground state can only be realized among the single-pulse states of Fig. \ref{fig-states} in the large-mass limit.

\subsection{Examples of the flower and dumbbell graphs}

Here we show how the construction of multi-pulse edge-localized states 
in the limit of large mass corresponds to the complete study of positive 
bound states on the flower and dumbbell graphs. The flower graph 
was considered in \cite{KMPZ}. The dumbbell graph was considered 
in \cite{G19,marzuola}. These examples and their generalizations 
were also reviewed with some additional numerical approximations in \cite{BMP,KP21}.

\begin{figure}[htbp]
	\centering
	\includegraphics[width=4in, height = 2in]{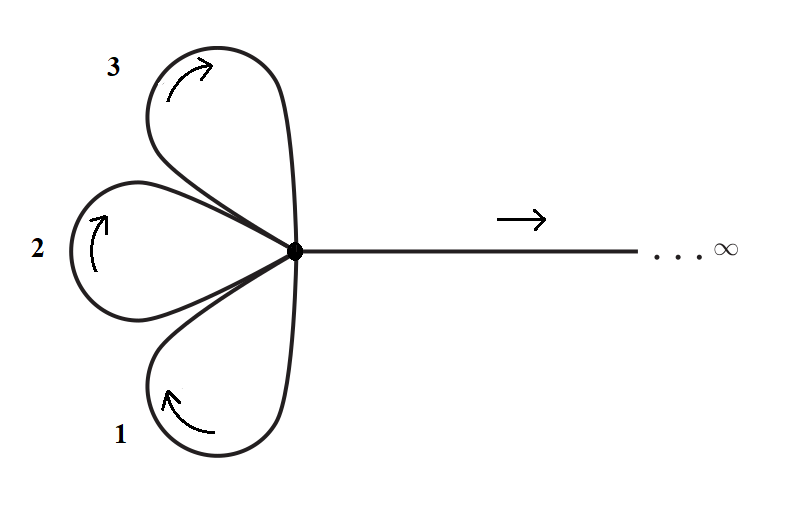}
	\caption{A flower graph with three loops. }
	\label{fig-flower}
\end{figure}

Let us consider the flower graph with $N = 3$ loops shown on Figure \ref{fig-flower}. 
Each $j$-th loop is parametrized by a segment $[-\ell_j, \ell_j]$ of length $2 \ell_j$  and the unbounded edge is parametrized by a half-line $[0, \infty)$. The consistency assumptions (\ref{constraint-on-length-1}) and (\ref{constraint-on-length-2}) are trivially satisfied 
because $|B| = 1$. Therefore, by Proposition \ref{prop-multi-pulse}, one can get single-pulse, double-pulse, 
and triple-pulse solutions in the limit of large mass. Moreover, 
the Morse index of the $N$-pulse state is $N$. 

\begin{figure}[htbp] 
	\centering
	\includegraphics[width=3in, height = 2in]{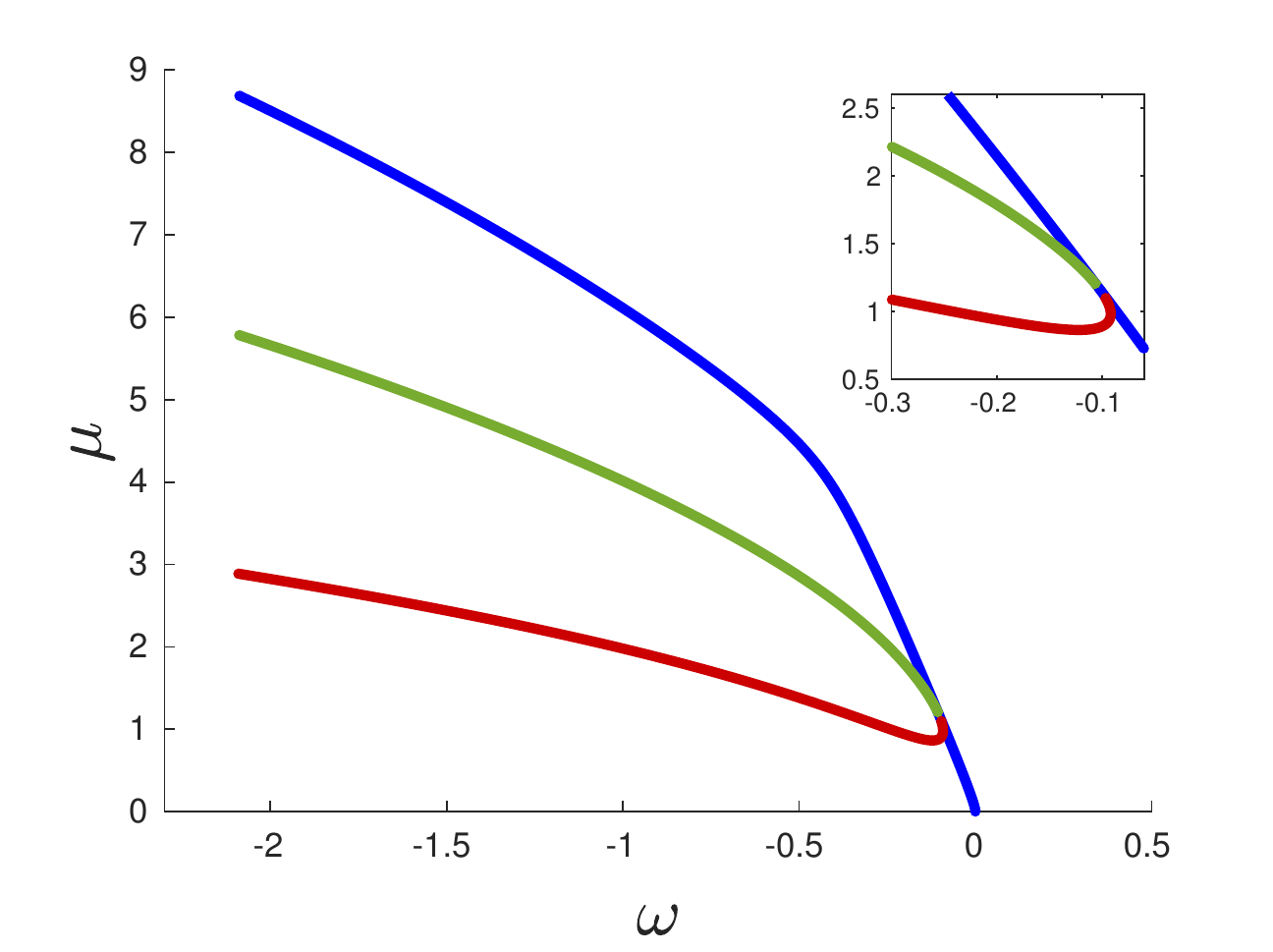}
	\includegraphics[width=3in, height = 2in]{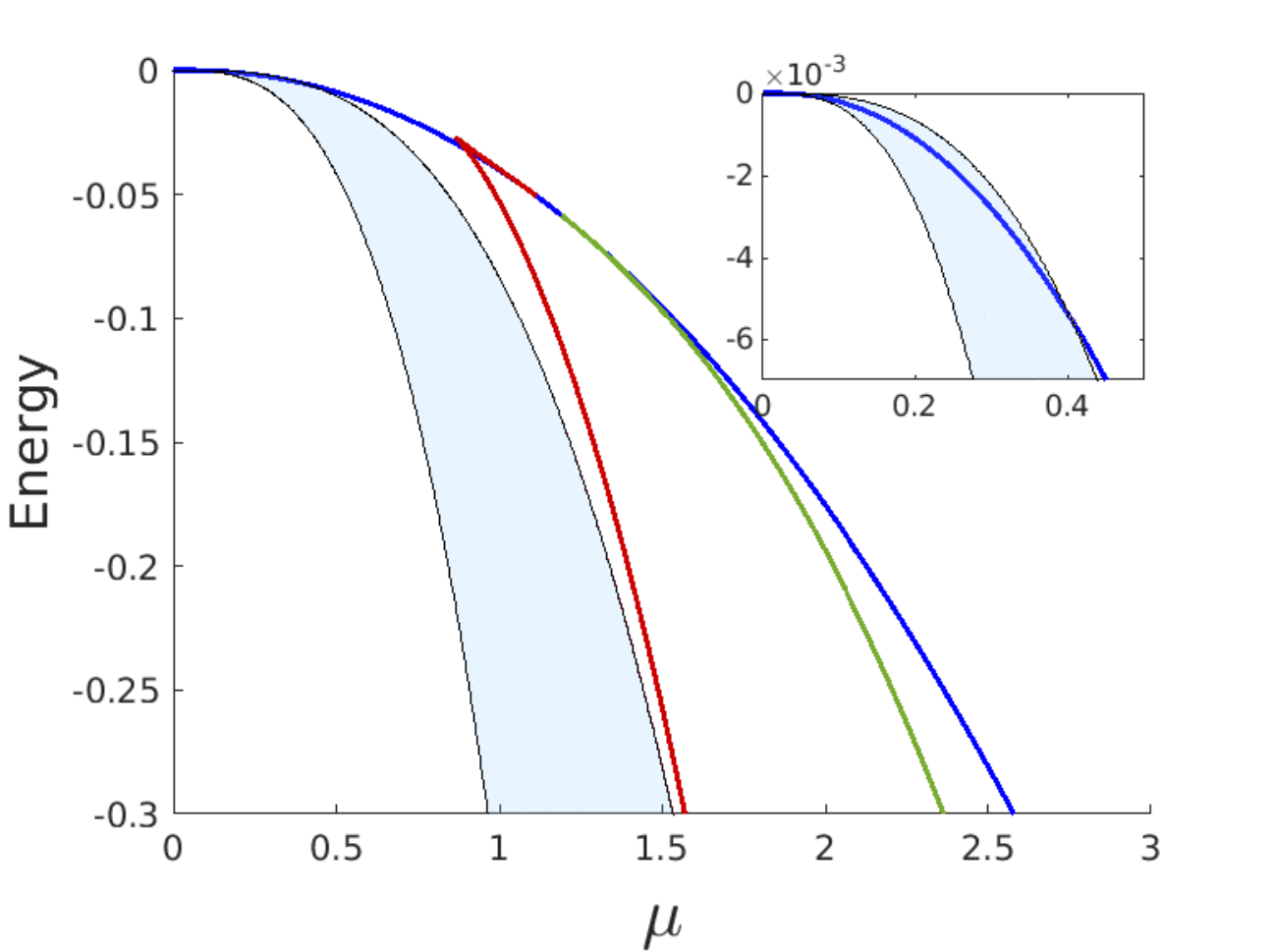}
	\caption{Bifurcation diagram of multi-pulse edge-localized states 
		on the parameter plane $(\omega,\mu)$ (left) 
		and on the mass--energy plane (right) for the flower graph 
		with $N = 3$ loops.  }
	\label{fig-flower-states}
\end{figure}

This conclusion coincides with 
Theorems 1,2,3 of \cite{KMPZ} and is illustrated on Figure \ref{fig-flower-states}.
The blue line shows the three-pulse state with Morse index equal to {\em one} before the bifurcation point (smaller values of mass) 
and to {\em three} after the bifurcation point (larger values of mass).
Therefore, the three-pulse state is orbitally stable before bifurcation 
and unstable after bifurcation. 

The green line depicts the double-pulse state with two components having larger amplitudes than the third one. It only exists for masses larger than 
the mass at the bifurcation point and it has Morse index equal 
to {\em two}, hence it is orbitally unstable. 

The red line shows a single-pulse state with one component having larger amplitude than the other components. 
It bifurcates to smaller masses than the mass at the bifurcation point but 
the solution branch turns at the fold bifurcation and 
then extends to the limit of large mass. The single-pulse state 
has Morse index is equal to {\em two} near the bifurcation point 
and to {\em one} after turning at the fold bifurcation. 
It is orbitally stable for large negative $\omega$.

The energy-mass diagram (right panel of Fig. \ref{fig-flower-states}) 
shows that the three-pulse solution for smaller mass belongs to the shaded area where minimizers of the variational problem 
(\ref{minimization}) exist, see inequalities in (\ref{inequalities-energy}). However, once the solution curve leaves the shaded area for larger mass, it never comes back. Hence, for sufficiently large mass, the infimum of the variational problem (\ref{minimization}) is not attained at the 
flower graph with $N = 3$ (in fact, for every $N \geq 2$)  
even at the single-pulse solution given by the red curve. 
Consequently, it is a local rather than global minimizer of the constrained variational problem (\ref{minimization}). No ground state 
exists for the flower graph with $N \geq 2$ in the large-mass limit, 
even thought the flower graph escapes the condition of Proposition \ref{AST-condition-H-proposition} because it has only one half-line.

\begin{figure}[htbp] 
	\centering
	\includegraphics[width=4in, height = 2in]{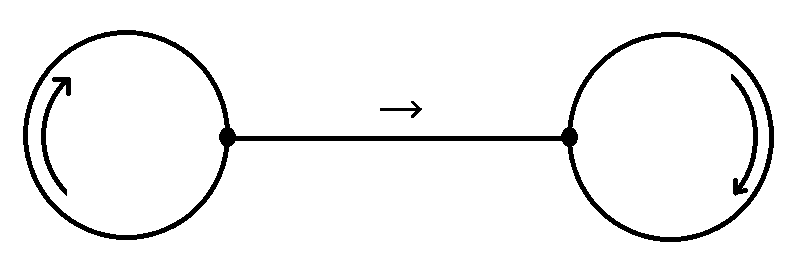}
	\caption{A dumbbell graph. }
	\label{fig-dumbbell}
\end{figure}

Let us now consider a dumbbell graph shown in Figure \ref{fig-dumbbell}.  The dumbbell graph consists of two loops of lengths $2\ell_-$ and $2 \ell_+$ connected by an internal edge of length $2 \ell_0$ at two vertices. The numerical approximations of the positive multi-pulse states were performed by using the Quantum Graph Package \cite{QGP} for $\ell_- = \ell_+ = \pi$, $\ell_0 = 2 < \ell_-$, and fixed $\omega = -4$.

\begin{figure}[htbp] 
	\centering
	\includegraphics[width=3.25in, height = 2.75in]{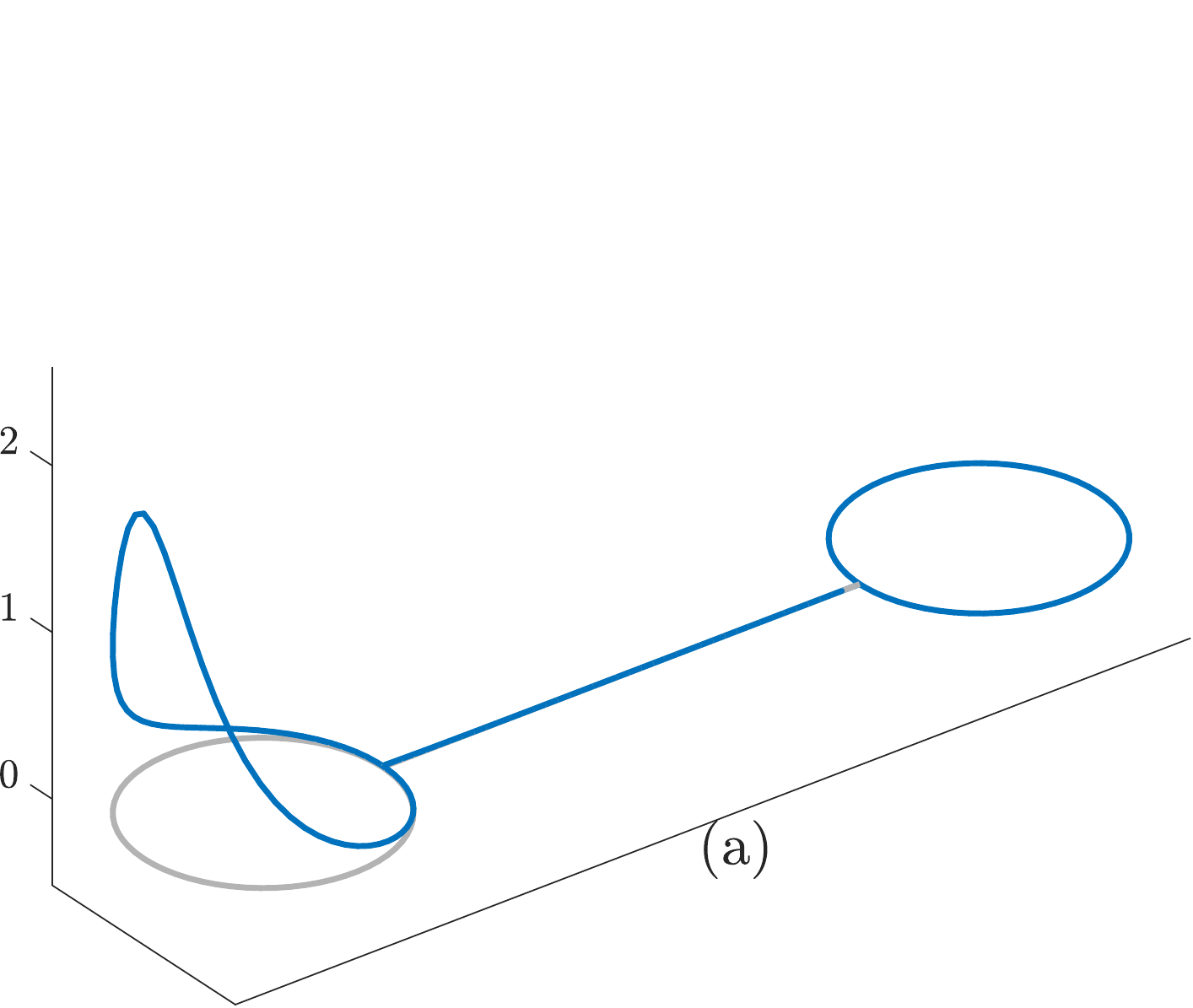}
	\includegraphics[width=3.25in, height = 2.75in]{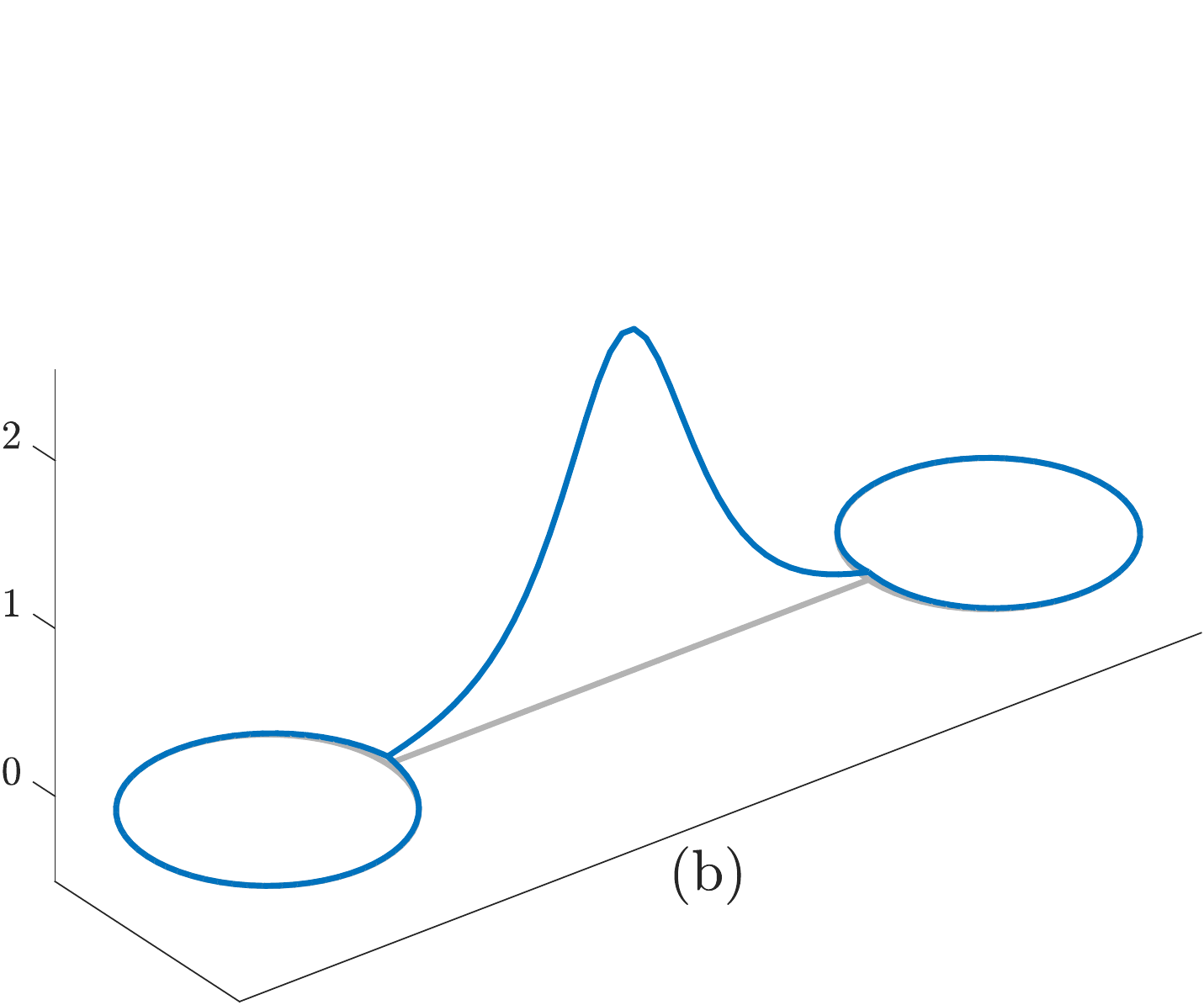} \\
	\includegraphics[width=3.25in, height = 2.75in]{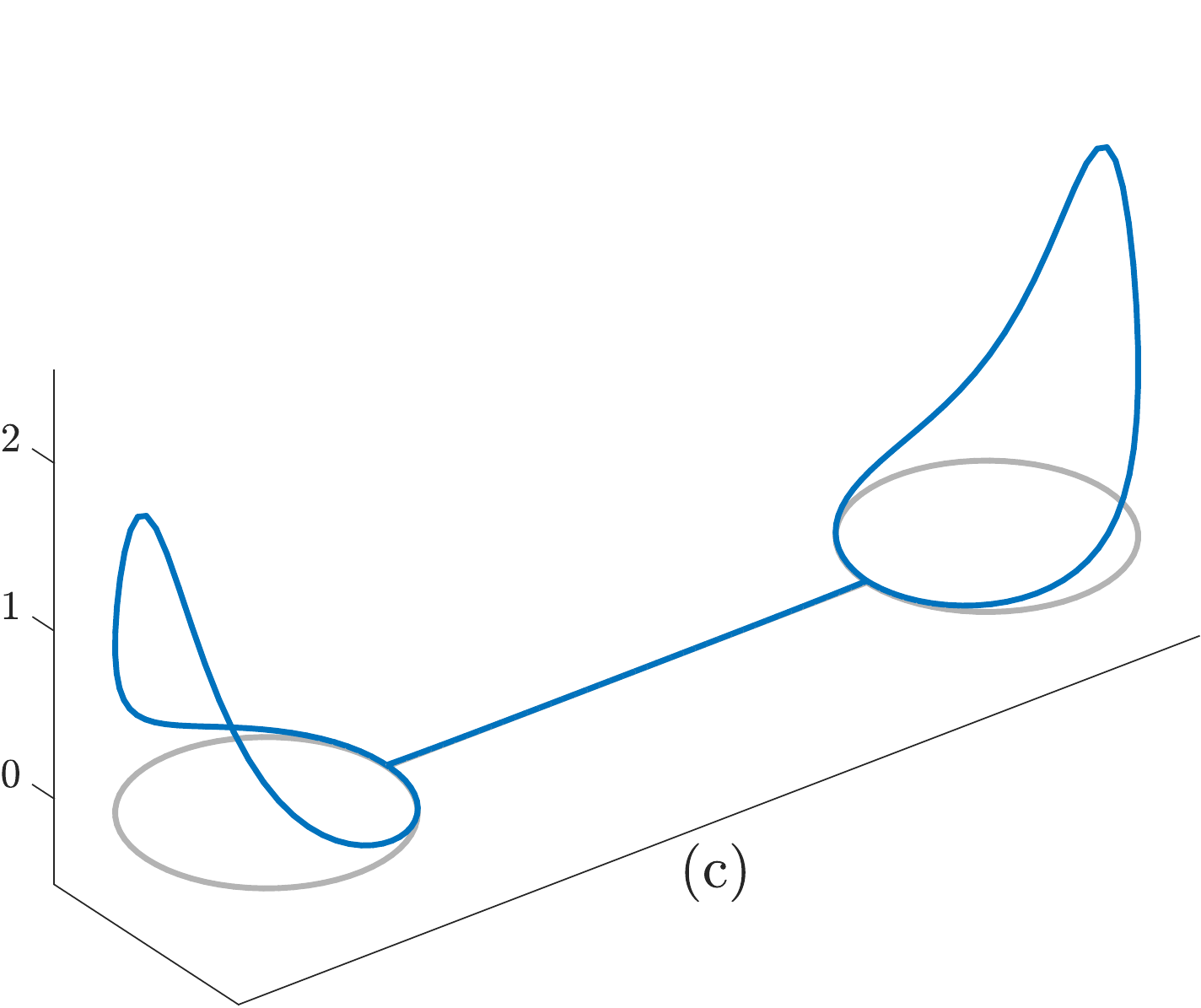}
	\includegraphics[width=3.25in, height = 2.75in]{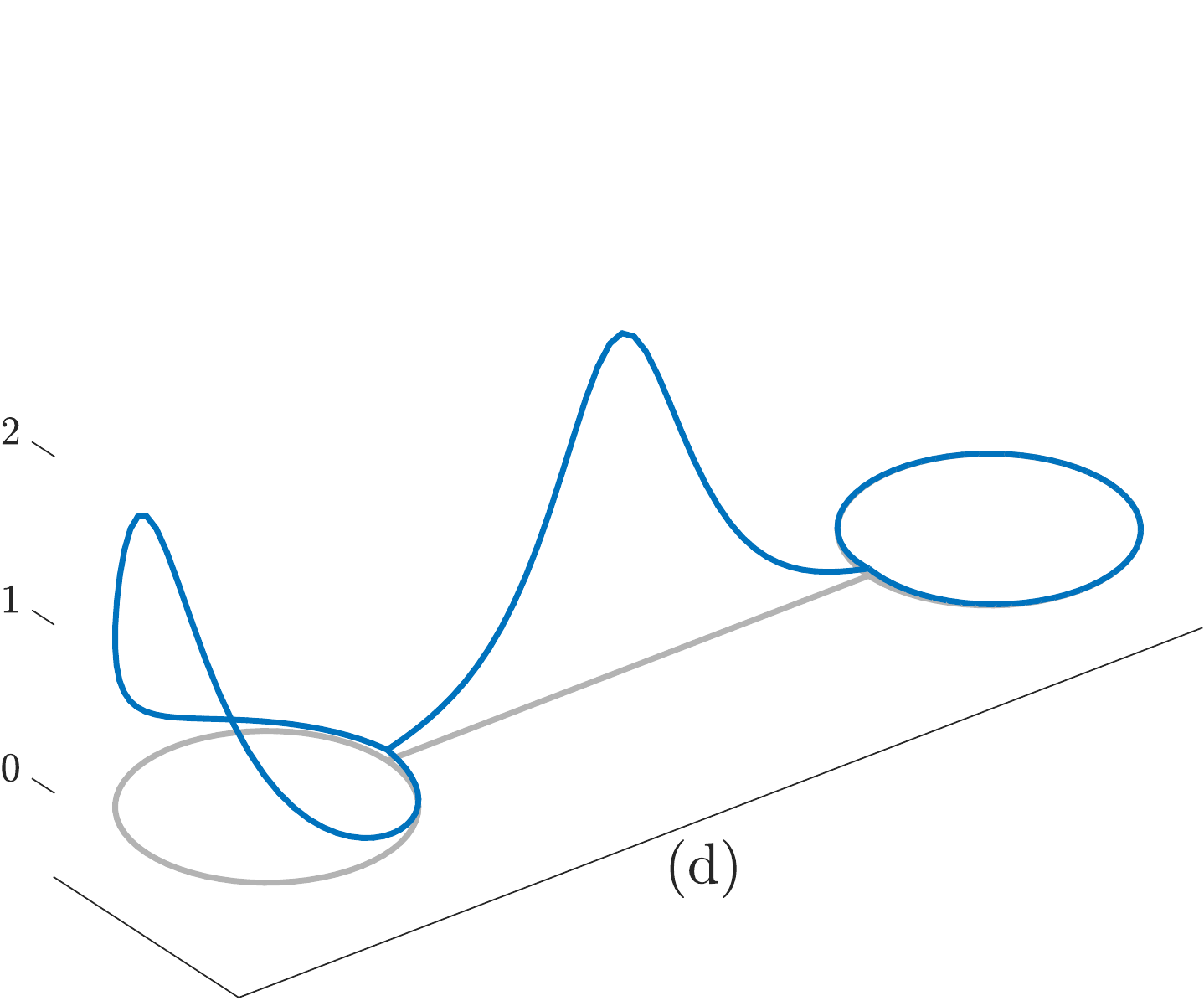}\\
	\includegraphics[width=3.25in, height = 2.75in]{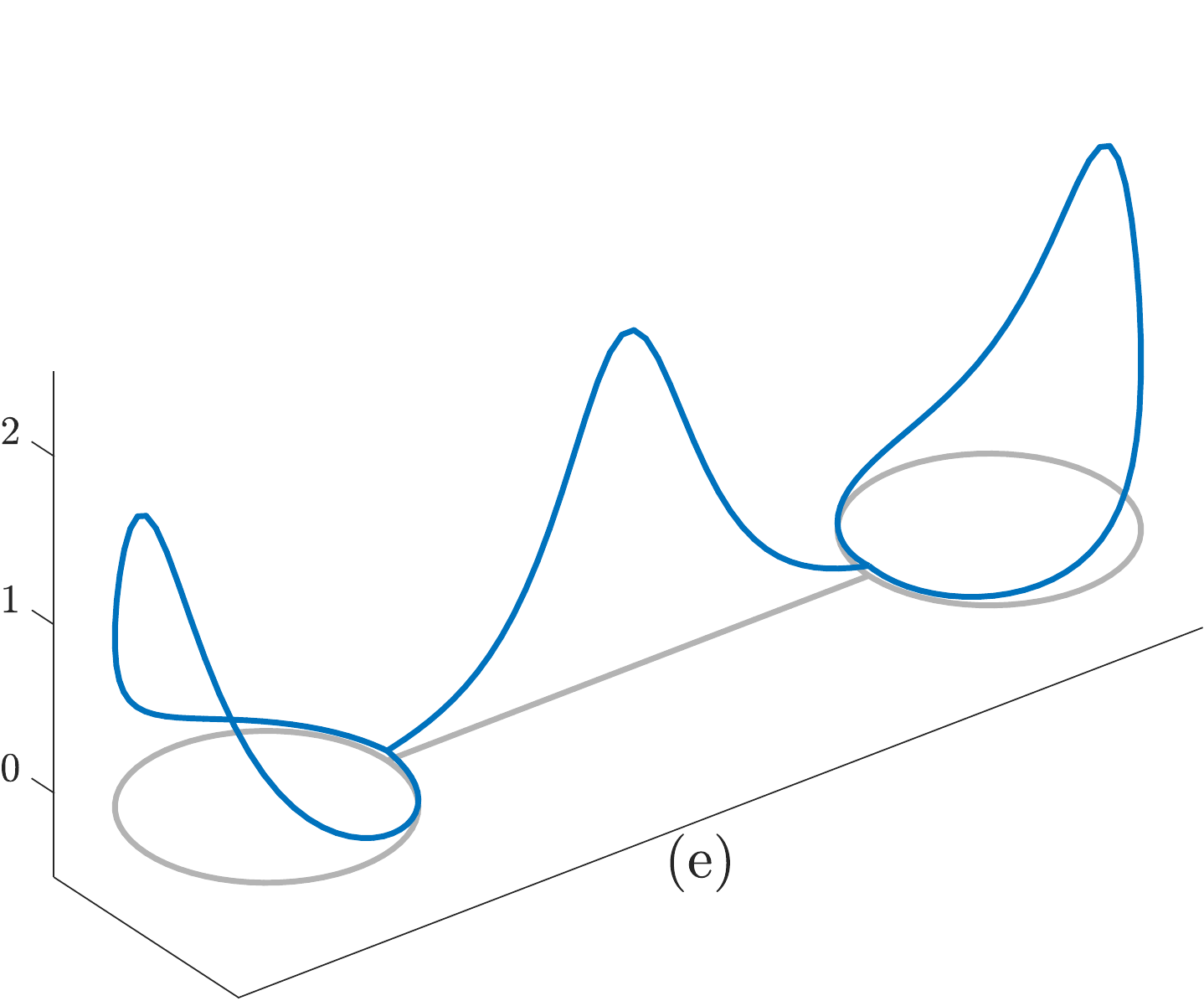}
	\caption{Five edge-localized states on the dumbbell graph with  $\ell_- = \ell_+ = \pi$, $\ell_0 = 2$, and $\omega = -4$.}
	\label{fig-states-dumbbell}
\end{figure}

The single-pulse state can be placed at any of the three edges in the limit of large mass. When it is placed on the loop, its Morse index is equal to {\em one}. It was shown in \cite{marzuola} (and also follows from the variational method in \cite{AST19}) that when it is placed on the internal edge, its Morse index is also equal to {\em one}. Figure \ref{fig-dumbbell} (a,b) show the two single-pulse states and we have confirmed that Morse index of each state is equal to {\em one}, hence, these states are orbitally stable for large negative $\omega$.

The double-pulse state can be placed at the two loops if $\ell_+ < 2 \ell_0 + \ell_-$ and $\ell_+ < 3 \ell_-$, which follow from the conditions (\ref{constraint-on-length-1}) and (\ref{constraint-on-length-2}) respectively. The constraints are satisfied if $\ell_- = \ell_+$. Morse index of this double-pulse state is equal to {\em two}, hence it is orbitally unstable. Figure \ref{fig-dumbbell} (c) shows this double-pulse state and we have confirmed that its Morse index is equal to {\em two}.
	
The double-pulse state can also be placed at one loop and the internal edge if $\ell_0 < \ell_- = \ell_+$. Figure \ref{fig-dumbbell} (d) shows this two-pulse state and we found that its Morse index is equal to {\em three}, hence it is orbitally unstable. This example shows that Proposition \ref{prop-multi-pulse} stating that the Morse index of $N$-pulse state is $N$ only holds for the multi-pulse states 
supported in the looping and pendant edges and may generally fail when some of the pulses are supported on internal edges. 
	
The three-pulse state can be placed at all three edges if $\ell_0 < \ell_- \leq \ell_+$. Figure \ref{fig-dumbbell} (e) shows this three-pulse state and we have found that its Morse index is equal to {\em five}, hence it is orbitally unstable. Again, Morse index 
exceeds the number of pulses of this three-pulse state.

\section{Other nonlinear dispersive wave models on metric graphs}
\label{sec-7}

Here we review other important nonlinear dispersive wave models 
which have been recently considered on metric graphs. 
The models include the Dirac equation, the Klein--Gordon equation, and the Korteweg--de Vries (KdV) equation.

\subsection{Dirac equation}

The nonlinear Dirac (NLD) equation with power type nonlinearity 
in one spatial dimension is written in the form:
\begin{equation}
\label{eqdirac}
i\f{\partial \Psi}{\partial t}=\mathsf{D}_m \Psi-|\Psi|^{2p+1}\Psi, \quad 
\mathsf{D}_m = -i c\sigma_{1} \frac{\partial\Psi}{\partial x} + m c^{2} \sigma_{3} \Psi,
\end{equation}
where $m\geq0$ is the mass, $c>0$ is the relativistic parameter, and $p > 0$ defines the nonlinearity power. The Dirac operator $\mathsf{D}_m$ is expressed in terms of the \emph{Pauli matrices} $\sigma_1$ and $\sigma_3$ given by 
\begin{equation}
\label{eq-pauli}
\sigma_1:=\begin{pmatrix}
0 & 1 \\
1 & 0
\end{pmatrix}
\qquad\text{and}\qquad
\sigma_3:=\begin{pmatrix}
1 & 0 \\
0 & -1
\end{pmatrix}.
\end{equation}
The NLD equation (\ref{eqdirac}) is an effective model for materials with Dirac points such as \emph{graphene} and \emph{germanene} \cite{WB-SB14} (see \cite{FW14} for a rigorous analysis of materials with Dirac points).
It has also been used for Bose--Einstein condensation in periodic traps \cite{HC09,PG07}. In some recent papers the NLD equation has been considered as an effective equation for BEC in branched quasi-unidimensional domains and optical fibers, when relativistic effects cannot be neglected (see \cite{p,SBMK18,TLB14}). 

A rigorous analysis of standing waves for the NLD equation on noncompact metric graphs has been recently reported in \cite{BCT19,BCT21}.
As for any differential equation on a metric graph one has to choose boundary conditions at vertices of the graph, which amounts to choose a specific self-adjoint realization of the Dirac operator \cite{BT90}. The preferred choice is given by the \emph{Kirchhoff} vertex conditions. 

Let us denote the Dirac operator with Kirchhoff vertex conditions by ${\mathsf D}_{\mathcal G}$.  The action of ${\mathsf D}_{\mathcal G}$ is defined on each edge as the action of the Dirac operator $\mathsf{D}_m$:
\begin{equation}
\label{eq-Dstand}
{\DD}_{\mathcal G} \psi_{e}=-i c\sigma_{1}\psi_e'+mc^{2}\sigma_{3}\psi_e,\qquad\forall e\in E.
\end{equation}
The domain of ${\mathsf D}_{\mathcal G}$ is defined by 
\[
\mathcal{D}(\DD_{\mathcal G}):=\{\psi=(\phi,\chi)^T\in H^1(\mathcal G):\text{\eqref{eq-cont_bis} and \eqref{eq-kirch_bis} are satisfied}\},
\]
where the Kirchhoff vertex conditions are 
\begin{gather}
\label{eq-cont_bis} \ \ \ \ \ \ \ \ \ \ \ \ \ \ \ \ \ \ \ \ \ \ \ \ \ \ \ \ \ \phi_e(v)=\phi_f(v),\qquad\forall e,f\succ v,\qquad\forall v\in V,\\[.2cm]
\label{eq-kirch_bis} \sum_{e\succ v}\chi_e(v)=\sum_{e \leftarrow v} \chi_e (v)  - \sum_{e \rightarrow v}
\chi_e (v) = 0,\qquad\forall v\in V.
\end{gather}

Let us add few remarks. Firstly, squaring the operator $\DD_{\mathcal G}$ and applying it to spinors with the nonzero first component, one obtains the Laplacian $\Delta_{\GG}$ with the Neumann--Kirchhoff conditions plus zero-order corrections, see \cite{BCT21}. Secondly, in application to the materials with Dirac points, the NLD equation is an effective model equation, and some caution should be used to interpretate the parameter $c>0$.

The NLD equation on the metric graph $\GG$ is given by
\begin{equation}
\label{eq-NLDtime}
i\f{\partial \Psi}{\partial t}=\DD_{\mathcal G}\Psi-|\Psi|^{2p+1}\Psi.
\end{equation}
One can also consider the model with localized nonlinearity, i.e.
\begin{equation}
\label{eq-NLDtimeconc}
i\f{\partial \Psi}{\partial t}=\DD_{\mathcal G}\Psi-\chi_{\K}|\Psi|^{2p+1}\Psi,
\end{equation}
where $\chi_{\K}$ is the indicator function of the compact part $\K$ of the metric graph $\GG$, where the nonlinearity is localized.
Bound states of the NLD equation \eqref{eq-NLDtime} are the spinors $\psi$ such that
\begin{equation}
\label{eq-NLD}
-i c\sigma_{1}\psi_e'+mc^{2}\sigma_{3}\psi_e-|\psi_e|^{2p+1}\psi_e=\omega \psi_e,\qquad \forall e\in E,
\end{equation}
where components of $\psi$ satisfies the Kirchhoff vertex conditions (\ref{eq-cont_bis}) and (\ref{eq-kirch_bis}). The wave function $\Psi(t,x)=e^{-i\omega t}\psi(x)$ is the standing wave solution of the NLD equation (\ref{eq-NLDtime}). The same definition obviously hold true for the model (\ref{eq-NLDtimeconc}), where the nonlinearity is concentrated, one can just add a factor $\chi_{\K}$ in front of the nonlinearity. 

In order to study the bound states of the stationary NLD equation (\ref{eq-NLD}), one should take into account the important fact that spectrum of $\DD_{\mathcal G}$ is given by (see \cite[Appendix A]{BCT19})
\[
\sigma(\DD_{\mathcal G})=(-\infty,-mc^2]\cup[mc^2,+\infty)
\]
As a consequence, the Dirac operator with the Kirchhoff vertex conditions has the gap  $\RE\backslash\sigma(\DD_{\mathcal G})=(-mc^2,mc^2)$. Bound states 
correspond to the values of $\omega$ in the gap $(-m c^2, mc^2)$. 

The NLD equation (\ref{eq-NLDtimeconc}) with localized nonlinearity has been considered in \cite{BCT19} and it is proved that for every $\omega\in(-mc^2,mc^2)$ there exist infinitely many distinct standing wave solutions. Moreover it is shown that for a sequence  $\{c_n\}_{n \in \mathbb{N}}$ of values of $c$ such that $c_n\to+\infty$, there exist another sequence $\{ \omega_n \}_{n \in \mathbb{N}}$ of values of $\omega$ such that  $\omega_n\to+\infty$ with corresponding bound states $\{\psi_n\}_{n \in \mathbb{N}}$ converging to a certain spinor of the form $(u,0)^T$, where $u \in \mathcal{D}(\Delta_{\GG})$ is a bound state satisfying the stationary NLS equation with concentrated nonlinearity:
\begin{equation}
\label{eq-NLSconc}
-\frac{1}{2m}u_e''-\chi_{\K}|u_e|^{p-2}u_e=\lambda u_e,\qquad\forall e\in E,
\end{equation}
Thus, the NLS equation on a metric graph $\GG$ becomes the non-relativistic limit of the NLD equation on the same metric graph $\GG$, for localized nonlinearities. 

The situation is technically more involved in the case of the NLD equation \eqref{eq-NLDtime}, due of the lack of compactness. Only the case of star graphs has been treated in \cite{BCT21} by means of bifurcation theory and not making use of variational methods. The analysis shows that starting from any real solution $u \in \mathcal{D}(\Delta_{\GG})$ of the NLS equation on the star graph with Neumann--Kirchhoff conditions,
\begin{equation}
\label{eq-NLS}
-\frac{1}{2m}u_e''-|u_e|^{p-2}u_e=\lambda u_e,\qquad\forall e\in E,
\end{equation}
one can define a branch of bound states of the stationary NLD equation \eqref{eq-NLD} for $\omega$ sufficiently near to $mc^2$ (and in the gap). In particular, there exists at least one branch of bound states of the NLD equation (\ref{eq-NLD})  bifurcating near the threshold $\omega=mc^2$. Moreover, one should notice that for a star graph with an odd number of edges there is a unique real solution of the stationary NLS equation \eqref{eq-NLS}, but for star graph with an even number of half lines there are infinitely many real solutions grouped in continuous families \cite{[ACFN14],KP1,KP2}.

\subsection{Nonlinear Klein--Gordon equation}

	The square of the Dirac operator is the Klein-Gordon operator, and the Klein-Gordon equation
	\begin{equation}
	\label{NKG}
	\frac{\partial^2\psi}{\partial t^2}=-\Delta\psi +m^2\psi + \lambda|\psi|^{2\sigma}\psi
	\end{equation}
	is another well-known dispersive partial differential equation. Notice that the Klein--Gordon equation (\ref{NKG}) is a nonlinear {\em wave equation}, second order in time. Recent analysis of the nonlinear Klein--Gordon equation (\ref{NKG}) on star graphs can be found in \cite{Gol21}.
	
	A variation on the nonlinear Klein--Gordon equation is the sine-Gordon equation
	\begin{equation}
	\label{sine-Gordon}
	\frac{\partial^2\psi}{\partial t^2}=-\Delta\psi + m^2\sin \psi.
	\end{equation}
This equation appears in several contexts, from field theory to differential geometry and DNA filaments dynamics. The sine-Gordon equation (\ref{sine-Gordon}) on a star graph has been recently used to model Josephson junctions in tri-crystal boundaries (see \cite{CD2014,Susanto} and  \cite{AngPlaza21-1,AngPlaza21-2,AngPlaza21-3}).

\subsection{Korteweg--de Vries equation} 

In its original form the Korteweg de Vries (KdV) equation reads
\begin{equation}\label{eq:kdv}
\frac{\partial \eta}{\partial t}=\frac{3}{2}\sqrt{\frac{g}{\ell}}\left(\frac{\sigma}{3}\frac{\partial^3 \eta}{\partial x^3} +\frac{2\gamma}{3}\frac{\partial \eta}{\partial x}+ \eta\frac{\partial \eta}{\partial x}\right)
\end{equation}
where the unknown $\eta$ is the elevation of the {\it shallow water} surface with respect to its average depth $\ell$ in a shallow canal. This dispersive equation has been the first where soliton solutions were discovered. In this case the soliton is just a translation with constant velocity of a given profile function.
Renaming coefficients and scaling variables one obtains the KdV equation with parameter $\alpha$ in front of the third derivative and $\beta$ in front of the first derivative
\begin{align}\label{KdV}
\frac{\partial u}{\partial t} &=\alpha \frac{\partial^3 u}{\partial x^3} +\beta \frac{\partial u}{\partial x} +u \frac{\partial u}{\partial x}.
\end{align}

The linear part of the KdV equation was reviewed in Section 1.4. 
When it is considered on the line, one usually takes $\alpha<0$ 
and set  $\beta=0$ due to the Galilean transformation. However, the Galilean transformation is not generally possible on the half-line or on metric graphs 
and the properties of the linear KdV equation depends on the signs of $\alpha$ and $\beta$.

Very little is known about the behavior of solutions of the KdV equation (\ref{KdV}) on metric graphs $\GG$. In \cite{Cav18}, well-posedness of the initial-value problem for the KdV equation on a star graph with three edges and special boundary conditions was considered in Sobolev spaces of low regularity. It was shown that the unique solution exists but in a very weak sense, in particular, it is not even twice differentiable in spatial coordinate. Further results in this direction were obtained in \cite{AC18,CCV20} for star graphs with bounded edges in the context of arterial flows.

Stability of standing waves for the KdV equation on a balanced star graph was considered in \cite{AC21}, based on \cite{MNS18}. The linear instability was proven for a special class of boundary conditions. Nothing is known about existence of stable bound states of the KdV equation on a star graph.

Finally, we mention other nonlinear equations related to fluid flow, in particular with application to arterial flow. The Benjamin--Bona--Mahony equation (or the regularized long wave equation) was considered on finite trees in \cite{bona}. Traveling waves of the same equation on more general metric graphs were considered in  \cite{MR14}. Shallow water equations in the context of river flow at forks were considered in \cite{Caputo}. Starting with a two-dimensional Boussinesq model in a forked channels region, a reduced one-dimensional equation on a metric graph was deduced in \cite{NS15}, with suitable boundary conditions at vertex. It was also shown in \cite{NS15} that the reduced model supports propagation of solitary waves.

\section{Further directions}

We have reviewed nonlinear evolutionary models on quantum graphs. Existence and stability of standing waves in the nonlinear Schr\"{o}dinger equation was studied with different analytical techniques such as the variational method, the period function, and the Dirichlet--to--Neumann mappings. Many fundamental results have already obtained for metric graphs both with bounded and unbounded edges under the assumption that the number of edges and vertices is finite. 
Extensions of these results to other evolutionary models such as the nonlinear Dirac equation, the nonlinear Klein--Gordon equation, and the Korteweg--de Vries equation have been recently considered with some preliminary results. 

In the end of this review, we would like to mention the recent work 
on the metric graphs with the unbounded number of edges. 

The NLS equation was considered on periodic graphs in one direction, such as a periodic chain with the fundamental cell consisting of a single loop and a single internal edge. For such graphs, homogenization to the NLS and NLD equations with constant coefficients 
was considered in \cite{p}. Standing wave solutions on the same periodic graph were classified by using dynamical system methods \cite{PS17} 
and variational methods \cite{Pankov18}. Existence of the ground state at every mass was proven for the periodic graph in \cite{Dovetta-per}. Symmetry 
of the ground state was clarified in the limit of large mass in \cite{BMP}.

In the case of the periodic graphs in two dimensions, an interesting dimensional crossover was found in \cite{ADST} in the context of the NLS equation with power nonlinearity. It was shown that if the power is below the cubic nonlinearity, ground states exist for every value of the mass, while 
if the power is between the cubic (including) and quintic (excluding) nonlinearities, ground states exist if and only if the mass exceeds a threshold value that depends on the power. This dimensional crossover is related 
to the coexistence of the one-dimensional and two-dimensional Sobolev inequalities, leading to a generalized Gagliardo–Nirenberg inequality for the doubly-periodic metric graphs. Failure of the dimensional crossover was investigated in \cite{DT21} for the doubly periodic graphs with compact or non-compact defects. 

Standing waves of the NLS equation on various metric graphs have approximated numerically with different techniques. Petviashvili's iteration method has been explored in \cite{marzuola} and \cite{BMP}. Other iterative and shooting methods were used in \cite{G19,KGP} and resulted in the Quantum Graph Package \cite{QGP}. The gradient method was applied to the stationary NLS equation in \cite{Besse}. 

Overall, standing waves of the nonlinear evolutionary models on quantum (metric) graphs is a rapidly growing area with many interesting developments.

\vspace{1cm}

{\bf Acknowledgments.} D. Noja acknowledges for funding the EC grant IPaDEGAN (MSCA-RISE-778010). D.E. Pelinovsky acknowledges the support of the NSERC Discovery grant.

\end{document}